\documentclass[oneside,final, letterpaper]{ucr}

\usepackage{amsmath, amssymb, amsfonts, mathrsfs, adjustbox,xurl}
\usepackage{tocloft}
\usepackage{fouriernc}

\usepackage{textcomp}
\newcommand{\textapprox}{\raisebox{0.5ex}{\texttildelow}}

\usepackage[toc, page]{appendix}
\usepackage[all,2cell,cmtip]{xy}
\usepackage[dvipsnames]{xcolor}
\usepackage{proof}
\usepackage[numbers]{natbib}
\definecolor{darkgreen}{rgb}{0,0.45,0}
\usepackage[colorlinks, citecolor=darkgreen, urlcolor=blue, final, hyperindex, pagebackref, linktoc = page, linkcolor = blue]{hyperref}

\usepackage{amsthm}

\usepackage[capitalize]{cleveref}

\theoremstyle{plain}
\newtheorem{thm}{Theorem}
\newtheorem{prop}[thm]{Proposition}
\newtheorem{cor}[thm]{Corollary}
\newtheorem{lem}[thm]{Lemma}

\theoremstyle{definition}
\newtheorem{defn}[thm]{Definition}
\newtheorem{expl}[thm]{Example}

\theoremstyle{remark}
\newtheorem*{rmk}{Remark}

\newcommand{\Z}{\mathbb{Z}}
\newcommand{\N}{\mathbb{N}}

\newcommand{\CMC}{\mathsf{CMC}}

\newcommand{\Set}{\mathsf{Set}}

\newcommand{\Petri}{\mathsf{Petri}}

\newcommand{\SSMC}{\mathsf{SSMC}}
\newcommand{\Csp}{\mathsf{Csp}}
\newcommand{\CMon}{\mathsf{CommMon}}
\newcommand{\Cat}{\mathsf{Cat}}

\newcommand{\X}{\mathsf{X}}

\newcommand{\Q}{\mathsf{Q}}
\newcommand{\sA}{\mathsf{A}}
\newcommand{\Cospan}{\lC \mathbf{sp}}
\newcommand{\Open}{\mathsf{Open}} 
\newcommand{\Span}{\mathsf{Span}}
\newcommand{\dRel}{\mathbb{R}\mathbf{el}}

\newcommand{\Hom}{\mathrm{Hom}}
\newcommand{\maps}{\colon}

\newcommand{\beq}{\begin{equation}}
\newcommand{\eeq}{\end{equation}}

\newcommand{\define}[1]{{\bf \boldmath{#1}}}

\newcommand{\Prof}{\mathsf{Prof}}

\newcommand{\lA}{\ensuremath{\mathbb{A}}}
\newcommand{\lB}{\ensuremath{\mathbb{B}}}
\newcommand{\lC}{\ensuremath{\mathbb{C}}}
\newcommand{\lD}{\ensuremath{\mathbb{D}}}

\newcommand{\law}[1]{\mathsf{#1}}
\newcommand{\Net}[1]{{#1}\text{-}\mathsf{Net}}

\newcommand{\Mod}{\mathsf{Mod}}
\newcommand{\Law}{\mathsf{Law}}
\newcommand{\PreNet}{\mathsf{PreNet}}
\newcommand{\cat}[1]{#1\text{-}\mathsf{Cat}}
\newcommand{\Grph}{\mathsf{Grph}}
\newcommand{\A}{^\bullet\mathrm{A}}
\newcommand{\B}{^{\bullet}\mathrm{B}}

\newcommand{\backb}{\mathrm{B}^{\bullet}}

\newcommand{\backa}{\mathrm{A}^{\bullet}}

\newcommand{\CAT}{\mathsf{CAT}}
\newcommand{\Mnd}{\mathsf{Mnd}}

\DeclareSymbolFont{bbold}{U}{bbold}{m}{n}
\DeclareSymbolFontAlphabet{\mathbbb}{bbold}

\DeclareMathOperator{\colim}{\mathrm{colim}}
\newcommand{\ten}{\otimes}
\newcommand{\Mor}{\mathrm{Mor}}
\newcommand{\Ob}{\mathrm{Ob}}

\newcommand{\po}{\ar@{}[dr]|{\text{\pigpenfont R}}}
\newcommand{\pb}{\ar@{}[dr]|{\text{\pigpenfont J}}}

\let\maps\colon

\newcommand{\RCat}{\mathsf{RCat}}
\newcommand{\RMat}{\mathsf{RMat}}

\newcommand{\op}{\mathrm{op}}

\newcommand{\Mat}{\mathsf{Mat}}

\newcommand{\Mon}{\mathsf{Mon}}

\usepackage{tikz}
\usepackage{tikz-cd}
\definecolor{lblue}{rgb}{0,250,255}
\tikzstyle{species}=[circle,fill=yellow,draw=black,scale=1.15]
\tikzstyle{inarrow}=[->, >=stealth, shorten >=.03cm,line width=1.5]
\tikzstyle{empty}=[circle,fill=none, draw=none]
\tikzstyle{inputdot}=[circle,fill=black,draw=black, scale=.25]
\tikzstyle{inputarrow}=[->,draw=purple, shorten >=.05cm]
\tikzstyle{simple}=[-,draw=black,line width=1.000]

\usetikzlibrary{backgrounds,circuits,circuits.ee.IEC,shapes,fit,matrix}

\tikzstyle{simple}=[-,line width=2.000]
\tikzstyle{arrow}=[-,postaction={decorate},decoration={markings,mark=at position .5 with {\arrow{>}}},line width=1.100]
\pgfdeclarelayer{edgelayer}
\pgfdeclarelayer{nodelayer}
\pgfsetlayers{edgelayer,nodelayer,main}

\tikzstyle{none}=[inner sep=0pt]

\definecolor{lblue}{rgb}{0,250,255}
\tikzstyle{species}=[circle,fill=yellow,draw=black,scale=1.15]
\tikzstyle{transition}=[rectangle,fill=lblue,draw=black,scale=1.15]
\tikzstyle{inarrow}=[->, >=stealth, shorten >=.03cm,line width=.8]
\tikzstyle{empty}=[circle,fill=none, draw=none]
\tikzstyle{inputdot}=[circle,fill=darkgreen,draw=Black, scale=.50]
\tikzstyle{node}=[circle,fill=blue,draw=Black, scale=0.8]
\tikzstyle{inputarrow}=[->,draw=purple, shorten >=.05cm]
\tikzstyle{simple}=[-,draw=black,line width=1.000]

\usetikzlibrary{arrows,shapes,snakes,automata,backgrounds,petri}
\tikzstyle{place}=[circle,thick,draw=blue!75,fill=blue!20,minimum size=6mm]
\tikzstyle{red place}=[place,draw=red!75,fill=red!20]
\tikzstyle{transition}=[rectangle,thick,draw=black!75,
  			  fill=black!20,minimum size=4mm]

\tikzstyle{dot}=[circle,fill=black,draw=black, scale=.4]

\begin{document}
\ssp
\title{Composing Behaviors of Networks\\ \today}
\author{Jade Edenstar Master}
\degreemonth{June}
\degreeyear{2020}
\degree{Doctor of Philosophy}
\chair{Dr. John C.\ Baez}
\othermembers{
Dr. Wee Liang Gan\\
Dr. Jose Gonzalez}
\numberofmembers{3}
\field{Mathematics}
\campus{Riverside}

\maketitle
\copyrightpage{}
\approvalpage{}

\degreesemester{Spring}

\begin{frontmatter}

\begin{acknowledgements}
I wouldn't have made it this far or even been the person I am today without the support of my friends and family. Allison I love you and feel so lucky to have you in my life. Thank you to my parents Eliza and Steve. I love you both so much thank you for everything. Thank you to my grandmother Babette. You have always been my biggest supporter and you inspire me to live my best life. Thank you to my friends for putting up with me. Thank you to Joyo, Mia, Alice, Casey, Maria, Ellie, and everyone on the server. I love you all. Thank you to everyone I know on the internet, I hope we meet in person someday. Thank you to all my haters, you pushed me to be a better person. 

Thank you to my advisor John Baez, I've learned so much from you and would not be where I am now without your dedication and hard work. Thank you to Christian Williams and Joe Moeller for sharing camaraderie and inspiration. Thank you to Daniel Cicala, Kenny Courser, Christina Vasilokopolou and Brandon Coya. Thank you to Anastasios Stefanou, Todd Trimble, Joshua Meyers, Sarah Rovner-Frydman, Zans Mihejez, Fabrizio Genovese, Morgan Rogers, Jelle Herold, Valeria De Paiva, Sarah Griffith, Alex Pokorny, Rany Tith, Evan Patterson, Arquimedes Canedo, Mike Shulman, Jules Hedges, and Martha Lewis. Math is created by communities not individuals and this thesis is no exception. 

Thank you to all the people who produced the things I need to survive. Thank you to the Tongva tribe, this work was done on their land.

The material from Chapter \ref{QNetchap} consists of work from ``Petri nets based on Lawvere theories" \cite{master_2020}. The material from Chapter \ref{openQNets} generalizes the work of ``Open Petri nets" joint with John Baez \cite{open}. Chapter \ref{enriched} and Chapter \ref{compalgpathchap} consist of work from ``The open algebraic path problem" \cite{openalgpathproblem}.

\end{acknowledgements} 

\begin{dedication}
    \null\vfil
    {\large
    \begin{center}
 To Allison, grow your wings and fly.
    \end{center}}
    \vfil\null
\end{dedication}

\begin{abstract}This thesis aims to develop a compositional theory for the operational semantics of networks. The networks considered are described by either internal or enriched graphs. In the internal case we focus on $\Q$-nets, a generalization of Petri nets based on a Lawvere theory $\Q$. $\Q$-nets include many known variants of Petri nets including pre-nets, integer nets, elementary net systems, and bounded nets. In the enriched case we focus on graphs enriched in a quantale $R$ regarded as matrices with entries in $R$. These $R$-matrices represent distance networks, Markov processes, capacity networks, non-deterministic finite automata, simple graphs, and more. The operational semantics of $\Q$-nets is constructed as an adjunction between $\Q$-nets and categories internal to the category of models of $\Q$. The left adjoint of this adjunction sends a $\Q$-net $P$ to an internal category $F_{\Q}(P)$ whose morphisms represent all possible firing sequences in $P$. Similarly, the operational semantics of $R$-matrices is constructed as an adjunction between $R$-matrices and categories enriched in $R$. The left adjoint of this adjunction sends an $R$-matrix $M$ to the $R$-category $F_R(M)$ whose hom-objects are solutions of the algebraic path problem: a generalization of the shortest path problem to graphs weighted in $R$. For both $\Q$-nets and $R$-matrices we use the theory of structured cospans to study the compositionality of the above operational semantics. For each type of network we construct a double category whose morphisms are ``open networks", i.e.\ networks with certain vertices designated as input or output. The operational semantics gives a double functor from a double category of open networks to a double category of open enriched or internal categories. These double functors give a compositional framework for computing the operational semantics of $\Q$-nets and $R$-matrices: their functoriality and coherence give relationships between the operational semantics of a network and the operational semantics of the smaller networks from which it is composed. We introduce the black-boxing of an open network, a profunctor describing the externally observable behavior of an open network. We introduce a class of open networks called ``functional open networks" for which black-boxing preserves composition. 
\end{abstract}

\setcounter{tocdepth}{1}
\begin{small}
\ssp
\tableofcontents
\end{small}
\end{frontmatter}

\chapter{Introduction}

\begin{quote}
     [...] the whole itself may be viewed as a conceptual construction, hence the question of the ontological status of boundaries becomes of a piece with the more general issue of the conventional status of ordinary objects and events. Cfr. Goodman: “We make a star as we make a constellation, by putting its parts together and marking off its boundaries” (1980: 213) --- \cite{sep-boundary}
\end{quote}

Underlying much of scientific thought is the assumption that things can be understood by understanding the way their components join together to make the whole. We often find when applying this point of view that the old adage rings true: ``things are more than the sum of their parts". When parts are joined together to form a whole, behavior emerges between the interaction of the parts that was not present before. Therefore, the magic must be in the way that the components are glued together. 

This thesis aims to provide a general setting to study the emergence that occurs when networks are glued together from their components. Here ``network" refers to a structure with a discrete set of states and a discrete set of relationships between them. In this thesis, the set of states may be equipped with operations representing the different ways that states are allowed to join together to make new states. The relationships in this thesis may be distances, probabilities, connections, processes, or external input.

\subsubsection*{Internal and Enriched Graphs} The different types of networks considered in this thesis are unified using enriched and internal graphs. All networks considered in this thesis are either graphs internal to a category $C$ or graphs enriched in a poset $R$. 

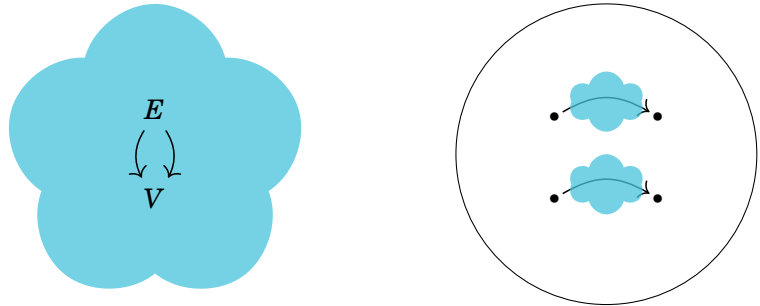
\[
\begin{tikzpicture}
  \node[cloud, cloud puffs=5, minimum width=4cm,minimum height=4cm, align=center, fill=SkyBlue,opacity=0.75] (cloud) at (-3, 0) {};
  \node at (-3,0) {$\begin{tikzcd}E \ar[d,bend left] \ar[d,bend right] \\ V \end{tikzcd}$};
  \draw (3, 0) circle (2cm);
  \node at (3,0) {$\begin{tikzcd}\bullet \ar[r,bend left,shorten <= -.5em,shorten >= -.5em] & \bullet \\
  \bullet \ar[r,bend left,shorten <=-.5em, shorten >= -.5em] & \bullet \end{tikzcd}$};
  \node[cloud, cloud puffs=6, minimum width=1cm, minimum height=.8cm, align=center, fill=SkyBlue, opacity=0.75] (cloud) at (3,.7) {};
  \node[cloud, cloud puffs=6, minimum width=1cm, minimum height=.8cm, align=center, fill=SkyBlue, opacity=0.75] (cloud) at (3,-.4) {};

\end{tikzpicture}
\]
As shown above left, a graph internal to $C$ is entirely in the clouds. In other words, its edges, vertices, source, and target maps are all objects and morphisms of $C$. An internal graph is equipped with the structure and properties of the category it lives in. If $C$ is the category of commutative monoids, then graphs internal to $C$ are presented by Petri nets. A Petri net is a diagram like this:
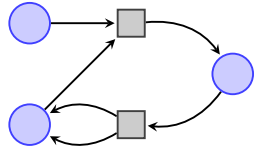
\[
\scalebox{0.9}{
\begin{tikzpicture}
  \begin{scope}
    \begin{pgfonlayer}{nodelayer}
        \node [place,tokens=0] (1a) at (-1.5,0.75) {};
        \node [place,tokens=0] (1b) at (-1.5,-0.75) {};      			
        \node [place,tokens=0] (3a) at (1.5,0) {};
        \node[transition] (2a) at (0,0.75) {};
        \node[transition] (2b) at (0,-0.75) {};
    \end{pgfonlayer}
    \begin{pgfonlayer}{edgelayer}
        \draw[style=inarrow, thick] (1a) to (2a);
        \draw[style=inarrow, thick] (1b) to (2a);
        \draw[style=inarrow, thick, bend right] (2b) to (1b);
        \draw[style=inarrow, thick, bend left] (2b) to (1b);
    
        \draw[style=inarrow, thick, bend left] (2a) to (3a);
        \draw[style=inarrow, thick, bend left] (3a) to (2b);
    \end{pgfonlayer}
  \end{scope}
\end{tikzpicture}}
\]
The circles, called ``places", represent different kinds of resources and the squares, called ``transitions", represent processes which take different resources as input and output. A discrete portion of a resource is represented by a ``token": a black dot inhabiting a place. A transition of a Petri net may ``fire" if there are enough tokens in the places with arrows going into it. When a transition fires, it removes one token from a place for each arrow going into it from that place. The transition then deposits one token into a place for each arrow coming from the transition. The following diagram represent the firing of the top transition of the above Petri net followed by a firing of the bottom transition:
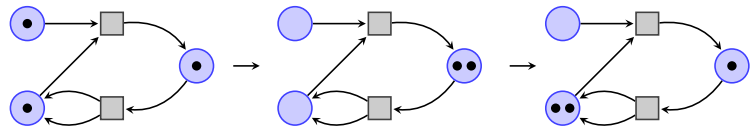
\[
\scalebox{0.75}{  
\begin{tikzpicture}
  \begin{scope}
    \begin{pgfonlayer}{nodelayer}
        \node [place,tokens=1] (1a) at (-1.5,0.75) {};
        \node [place,tokens=1] (1b) at (-1.5,-0.75) {};      			
        \node [place,tokens=1] (3a) at (1.5,0) {};
        \node[transition] (2a) at (0,0.75) {};
        \node[transition] (2b) at (0,-0.75) {};
    \end{pgfonlayer}
    \begin{pgfonlayer}{edgelayer}
        \draw[style=inarrow, thick] (1a) to (2a);
        \draw[style=inarrow, thick] (1b) to (2a);
        \draw[style=inarrow, thick, bend right] (2b) to (1b);
        \draw[style=inarrow, thick, bend left] (2b) to (1b);
    
        \draw[style=inarrow, thick, bend left] (2a) to (3a);
        \draw[style=inarrow, thick, bend left] (3a) to (2b);
    \end{pgfonlayer}
  \end{scope}

  \begin{scope}[xshift=135]
    \begin{pgfonlayer}{nodelayer}
        \node [place,tokens=0] (1a) at (-1.5,0.75) {};
        \node [place,tokens=0] (1b) at (-1.5,-0.75) {};      			
        \node [place,tokens=2] (3a) at (1.5,0) {};
        \node[transition] (2a) at (0,0.75) {};
        \node[transition] (2b) at (0,-0.75) {};
    \end{pgfonlayer}
    \begin{pgfonlayer}{edgelayer}
        \draw[style=inarrow, thick] (1a) to (2a);
        \draw[style=inarrow, thick] (1b) to (2a);
        \draw[style=inarrow, thick, bend right] (2b) to (1b);
        \draw[style=inarrow, thick, bend left] (2b) to (1b);
    
        \draw[style=inarrow, thick, bend left] (2a) to (3a);
        \draw[style=inarrow, thick, bend left] (3a) to (2b);
    \end{pgfonlayer}
  \end{scope}

    \draw[style=inarrow, thick] (2.15,0) -- (2.65,0);

  \begin{scope}[xshift=270]
    \begin{pgfonlayer}{nodelayer}
        \node [place,tokens=0] (1a) at (-1.5,0.75) {};
        \node [place,tokens=2] (1b) at (-1.5,-0.75) {};      			
        \node [place,tokens=1] (3a) at (1.5,0) {};
        \node[transition] (2a) at (0,0.75) {};
        \node[transition] (2b) at (0,-0.75) {};
    \end{pgfonlayer}
    \begin{pgfonlayer}{edgelayer}
        \draw[style=inarrow, thick] (1a) to (2a);
        \draw[style=inarrow, thick] (1b) to (2a);
        \draw[style=inarrow, thick, bend right] (2b) to (1b);
        \draw[style=inarrow, thick, bend left] (2b) to (1b);
    
        \draw[style=inarrow, thick, bend left] (2a) to (3a);
        \draw[style=inarrow, thick, bend left] (3a) to (2b);
    \end{pgfonlayer}
  \end{scope}
    \draw[style=inarrow, thick] (7.05,0) -- (7.55,0);

\end{tikzpicture}
}
\]
More generally, graphs internal to $\Mod(\Q)$, the category of models for a Lawvere theory $\Q$, are presented by $\Q$-nets. We obtain many variants of Petri nets including pre-nets, elementary net systems, $k$-safe nets, lending nets, and more by generalizing in this way. 

On the other hand, a graph enriched in a poset $R$ has only its edges in the clouds: its vertices are elements of sets, but for each pair of vertices $x$ and $y$, there is an element $R(x,y)$ of $R$ representing the connection between $X$ and $Y$. In other words, an $R$-enriched graph on a set $X$ is a square matrix \[ M \maps X \times X \to R \]
called an $R$-matrix. $R$ should be nice enough to enrich in, and it is enough for $R$ to be a commutative ``quantale": a monoidal closed poset with all joins. $R$-matrices are represented graphically as weighted graphs. When $R$ is the quantale $[0, \infty]$, of positive real numbers, the values $R(x,y)$ are regarded as distances between a set of locations. A $[0,\infty]$-matrix can be drawn as a graph where the distance $R(x,y)$ labels an edge from $x$ from $y$:
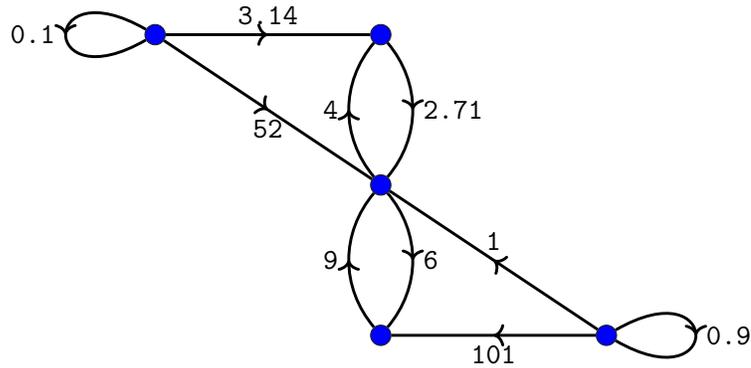
\[ 
\begin{tikzpicture}
	\begin{pgfonlayer}{nodelayer}
		\node [style=node] (0) at (-3, 2) {};
		\node [style=node] (1) at (0, -0) {};
		\node [style=node] (2) at (0, 2) {};
		\node [style=node] (3) at (0, -2) {};
		\node [style=node] (4) at (3, -2) {};
	\end{pgfonlayer}
	\begin{pgfonlayer}{edgelayer}
	    \draw [style=arrow] (0) to [out=150,in=210,looseness=30] node [left] {\tt 0.1} (0);
		\draw [style=arrow] (0) to node [above] {\tt 3.14} (2);
		\draw [style=arrow] (2) to [bend left=40] node [right] {\tt 2.71} (1);
		\draw [style=arrow] (1) to [bend left=40] node [left] {\tt 4} (2);
		\draw [style=arrow] (1) to [bend left=40] node [right] {\tt 6}(3);
		\draw [style=arrow] (3) to [bend left=40] node [left] {\tt 9} (1);
		\draw [style=arrow] (4) to node [below] {\tt 101} (3);
		\draw [style=arrow] (0) to node [below] {\tt 52} (1);
		\draw [style=arrow] (4) to node [above] {\tt 1} (1);
		\draw [style=arrow] (4) to [out=30,in=-30,looseness=30] node [right] {\tt 0.9} (4);
	\end{pgfonlayer}
\end{tikzpicture}
\]
where pairs of vertices without an edge between them are assumed to have a distance of $\infty$. Varying $R$ gives Markov processes, finite state machines, simple graphs, capacity networks, and more as instances of enriched graphs.

The following tree summarizes the networks considered in this thesis. This tree is certainly non-exhaustive as more examples may be derived from the general theory developed here.
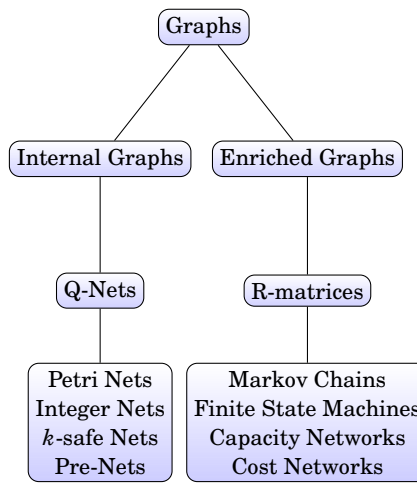
\[
\begin{tikzpicture}[scale=0.7,sibling distance=11em, level distance=7em,
  every node/.style = {shape=rectangle, rounded corners,
    draw, align=center,
    top color=white, bottom color=blue!20}]]
    \begin{footnotesize}
  \node {Graphs}
    child { node {Internal Graphs} 
        child { node {Q-Nets} 
            child {node {Petri Nets \\ Integer Nets \\ $k$-safe Nets \\ Pre-Nets} }
            }
        }
    child { node {Enriched Graphs}
      child { node {R-matrices}
        child { node {Markov Chains \\ Finite State Machines \\ Capacity Networks \\ Cost Networks} } } };
        \end{footnotesize}
\end{tikzpicture}
\] 

\subsubsection*{Operational Semantics of Networks} Networks are just formal structures until they are equipped with a semantics indicating their real world meaning. The semantics that we equip networks with is the ``operational semantics": a mathematical specification of the ways the states of the network may evolve in time. For an ordinary graph $G$, its operational semantics will be a category $F(G)$ whose objects are the vertices of $G$ and whose morphisms consist of all paths which can be formed from the edges of $G$. $F$ becomes a functor so that morphisms of graphs extend to behavior preserving functors between their operational semantics. $F$ is the left adjoint of the adjunction \[
\begin{tikzcd}
\Grph \ar[r, bend left,"F",pos=.45] \ar[r,phantom,"\bot"] & \ar[l, bend left, "U",pos=.55] \Cat
\end{tikzcd}
\]
where $U$ is the functor which sends any category to its underlying graph. This perspective allows for a uniform treatment of the operational semantics for internal and enriched graphs. In Section \ref{freecatinternalsec} we show how graphs internal to $\Mod(\Q)$ generate categories internal to $\Mod(\Q)$ via the adjunction
\[
\begin{tikzcd}
\Net{\Q} \ar[r, bend left,"\B"] \ar[r,phantom,"\bot",pos=.7] & \ar[l, bend left, "\backb"] \Grph(\Mod(\Q))
\end{tikzcd}
\]
where $\Grph(\Mod(\Q)$ is the category of graphs internal to $\Mod(\Q)$ and $\Cat(\Mod(\Q))$ is the category of small categories internal to $\Mod(\Q)$. This does not yet give the operational semantics for $\Q$-nets. In Section \ref{freeQgraph} we construct the missing piece, an adjunction 
\[ 
\begin{tikzcd}
\Net{\Q} \ar[r, bend left,"\A"] \ar[r,phantom,"\bot",pos=.7] & \ar[l, bend left, "\backa"] \Grph(\Mod(\Q))
\end{tikzcd}
\]\noindent making precise the sense in which $\Q$-nets are the generating data for graphs internal to $\Mod(\Q)$. In Theorem \ref{big} we combine these to obtain the operational semantics for $\Q$-nets
\[\begin{tikzcd}
\Net{\Q} \ar[r, bend left,"F_{\Q}"] \ar[r,phantom,"\bot"] & \ar[l, bend left, "U_{\Q}"] \cat{\Q}
\end{tikzcd} \]
In the case when $\Q$ is the Lawvere theory for commutative monoids, this operational semantics is similar to the adjunction developed by Meseguer and Montanari \cite{monoids}. Letting $\Q$ be the theory of monoids reproduces the operational semantics for pre-nets introduced in \cite{functorialsemantics}, and setting $\Q$ equal to the theory of abelian groups gives the operational semantics for lending nets developed by Genovese and Herold \cite{genovese}. Allowing $\Q$ to be other Lawvere theories gives a new categorical characterizations of the operational semantics for many other variants of Petri nets.

In Chapter \ref{enriched}, for each commutative quantale $R$, we construct an adjunction
\[
\begin{tikzcd}
\RMat \ar[r, bend left,"F_R"] \ar[r,phantom,"\bot"] & \ar[l, bend left, "U_R"] \RCat
\end{tikzcd}
\]
between matrices valued in $R$ and categories enriched in $R$. The left adjoint of this adjunction is familar: for an $R$-matrix $M$, the $R$-category $F_R(M)$ is a matrix whose entries are solutions to the algebraic path problem. The algebraic path problem is a generalization of the shortest path problem to probability, computing, matrix multiplication, and optimization \cite{tarjan1981unified,foote2015kleene}. When $R$ is the quantale of positive real numbers $([0,\infty], \mathrm{min}, +)$, a weighted graph can be regarded as an $R$-matrix, and the shortest paths of this graph are given by $F_R(M)$. The algebraic path problem allows $R$ to vary, and gets problems of a similar flavor also as the free $R$-category on an $R$-matrix. Many popular shortest path algorithms can be extended to compute solutions to the algebraic path problem in a general setting \cite{hofner2012dijkstra}. The algebraic path problem can also be implemented generically using functional programming \cite{dolan2013fun}. The above adjunction makes clear the universal property of the solutions to the algebraic path problem.



\noindent \subsubsection*{Open Networks} The main goal of this thesis is to study how the operational semantics of networks can be joined together. To compose networks we first need to equip them with boundaries. A network $G$ with vertex set $V$ is made ``open" to its surroundings by equipping it with functions $i \maps X \to V$ and $o \maps Y \to V$ designating input and output vertices respectively. This is formalized as a cospan
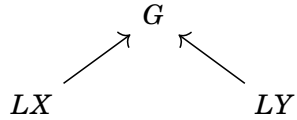
\[\begin{tikzcd}
 & G& \\
 LX \ar[ur] & & LY \ar[ul]
\end{tikzcd} \]
where $LX$ and $LY$ denote the discrete networks on the sets $X$ and $Y$. Similarly, a Petri net is made open by equipping it with functions from its input and output sets to its places. Open Petri nets will be the running example for the remainder of this introduction. Here is an open Petri net $P$ with input set $X$ and output set $Y$:
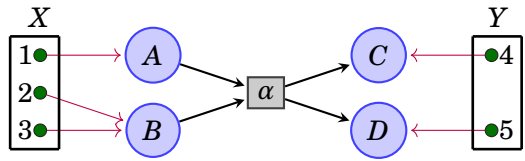
\[
\begin{tikzpicture}
	\begin{pgfonlayer}{nodelayer}
		\node [style=place] (A) at (-4, 0.5) {$A$};
		\node [style=place] (B) at (-4, -0.5) {$B$};
		\node [style=place] (C) at (-1, 0.5) {$C$};
		\node [style=place] (D) at (-1, -0.5) {$D$};
             \node [style=transition] (a) at (-2.5, 0) {$\alpha$}; 
		
		\node [style=empty] (X) at (-5.5, 1) {$X$};
		\node [style=none] (Xtr) at (-5.25, 0.75) {};
		\node [style=none] (Xbr) at (-5.25, -0.75) {};
		\node [style=none] (Xtl) at (-5.9, 0.75) {};
             \node [style=none] (Xbl) at (-5.9, -0.75) {};
	
		\node [style=inputdot] (1) at (-5.5, 0.5) {};
		\node [style=empty] at (-5.7, 0.5) {$1$};
		\node [style=inputdot] (2) at (-5.5, 0) {};
		\node [style=empty] at (-5.7, 0) {$2$};
		\node [style=inputdot] (3) at (-5.5, -0.5) {};
		\node [style=empty] at (-5.7, -0.5) {$3$};

		\node [style=empty] (Y) at (0.6, 1) {$Y$};
		\node [style=none] (Ytr) at (.9, 0.75) {};
		\node [style=none] (Ytl) at (0.25, 0.75) {};
		\node [style=none] (Ybr) at (.9, -0.75) {};
		\node [style=none] (Ybl) at (.25, -0.75) {};

		\node [style=inputdot] (4) at (0.5, 0.5) {};
		\node [style=empty] at (0.7, 0.5) {$4$};
		\node [style=inputdot] (5) at (0.5, -0.5) {};
		\node [style=empty] at (0.7, -0.5) {$5$};		
	
	\end{pgfonlayer}
	\begin{pgfonlayer}{edgelayer}
		\draw [style=inarrow] (A) to (a);
		\draw [style=inarrow] (B) to (a);
		\draw [style=inarrow] (a) to (C);
		\draw [style=inarrow] (a) to (D);
		\draw [style=inputarrow] (1) to (A);
		\draw [style=inputarrow] (2) to (B);
		\draw [style=inputarrow] (3) to (B);
		\draw [style=inputarrow] (4) to (C);
		\draw [style=inputarrow] (5) to (D);
		\draw [style=simple] (Xtl.center) to (Xtr.center);
		\draw [style=simple] (Xtr.center) to (Xbr.center);
		\draw [style=simple] (Xbr.center) to (Xbl.center);
		\draw [style=simple] (Xbl.center) to (Xtl.center);
		\draw [style=simple] (Ytl.center) to (Ytr.center);
		\draw [style=simple] (Ytr.center) to (Ybr.center);
		\draw [style=simple] (Ybr.center) to (Ybl.center);
		\draw [style=simple] (Ybl.center) to (Ytl.center);
	\end{pgfonlayer}
\end{tikzpicture}
\]
The functions from $X$ and $Y$ into the set of places indicate points at which tokens could flow in or out. We write this open Petri net as $P \maps X \to Y$ for short.
There are two fundamental operations on open networks. First, they may be composed along a shared boundary. Given another open Petri net $Q \maps Y \to Z$:
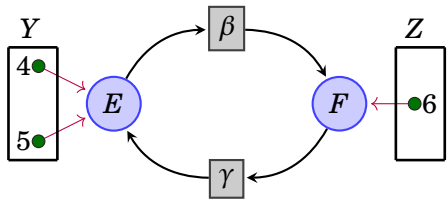
\[
\begin{tikzpicture}
	\begin{pgfonlayer}{nodelayer}

		\node [style = transition] (b) at (2.5, 1) {$\beta$};
		\node [style = transition] (c) at (2.5, -1) {$\gamma$};
		\node [style = place] (E) at (1, 0) {$E$};
		\node [style = place] (F) at (4,0) {$F$};
		
	

		\node [style=empty] (Y) at (-0.1, 1) {$Y$};
		\node [style=none] (Ytr) at (.25, 0.75) {};
		\node [style=none] (Ytl) at (-.4, 0.75) {};
		\node [style=none] (Ybr) at (.25, -0.75) {};
		\node [style=none] (Ybl) at (-.4, -0.75) {};

		\node [style=inputdot] (4) at (0, 0.5) {};
		\node [style=empty] at (-0.2, 0.5) {$4$};
		\node [style=inputdot] (5) at (0, -0.5) {};
		\node [style=empty] at (-0.2, -0.5) {$5$};		
		
		\node [style=empty] (Z) at (5, 1) {$Z$};
		\node [style=none] (Ztr) at (4.75, 0.75) {};
		\node [style=none] (Ztl) at (5.4, 0.75) {};
		\node [style=none] (Zbl) at (5.4, -0.75) {};
		\node [style=none] (Zbr) at (4.75, -0.75) {};

		\node [style=inputdot] (6) at (5, 0) {};
		\node [style=empty] at (5.2, 0) {$6$};	
		
	\end{pgfonlayer}
	\begin{pgfonlayer}{edgelayer}
		\draw [style=inarrow, bend left=30, looseness=1.00] (E) to (b);
		\draw [style=inarrow, bend left=30, looseness=1.00] (b) to (F);
		\draw [style=inarrow, bend left=30, looseness=1.00] (c) to (E);
		\draw [style=inarrow, bend left=30, looseness=1.00] (F) to (c);
		\draw [style=inputarrow] (4) to (E);
		\draw [style=inputarrow] (5) to (E);
		\draw [style=inputarrow] (6) to (F);
		\draw [style=simple] (Ytl.center) to (Ytr.center);
		\draw [style=simple] (Ytr.center) to (Ybr.center);
		\draw [style=simple] (Ybr.center) to (Ybl.center);
		\draw [style=simple] (Ybl.center) to (Ytl.center);
		\draw [style=simple] (Ztl.center) to (Ztr.center);
		\draw [style=simple] (Ztr.center) to (Zbr.center);
		\draw [style=simple] (Zbr.center) to (Zbl.center);
		\draw [style=simple] (Zbl.center) to (Ztl.center);
	\end{pgfonlayer}
\end{tikzpicture}
\]
the first step in composing $P$ and $Q$ is to put the pictures together:
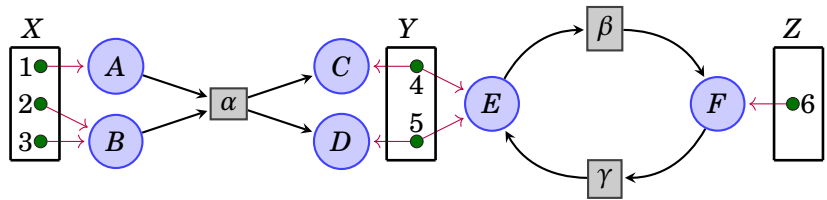
\[
\begin{tikzpicture}
	\begin{pgfonlayer}{nodelayer}
		\node [style=place] (A) at (-4, 0.5) {$A$};
		\node [style=place] (B) at (-4, -0.5) {$B$};
		\node [style=place] (C) at (-1, 0.5) {$C$};
		\node [style=place] (D) at (-1, -0.5) {$D$};
            \node [style=transition] (a) at (-2.5, 0) {$\alpha$}; 
		\node [style = transition] (b) at (2.5, 1) {$\beta$};
		\node [style = transition] (c) at (2.5, -1) {$\gamma$};
		\node [style = place] (E) at (1, 0) {$E$};
		\node [style = place] (F) at (4,0) {$F$};
		
		\node [style=empty] (X) at (-5.1, 1) {$X$};
		\node [style=none] (Xtr) at (-4.75, 0.75) {};
		\node [style=none] (Xbr) at (-4.75, -0.75) {};
		\node [style=none] (Xtl) at (-5.4, 0.75) {};
           \node [style=none] (Xbl) at (-5.4, -0.75) {};
	
		\node [style=inputdot] (1) at (-5, 0.5) {};
		\node [style=empty] at (-5.2, 0.5) {$1$};
		\node [style=inputdot] (2) at (-5, 0) {};
		\node [style=empty] at (-5.2, 0) {$2$};
		\node [style=inputdot] (3) at (-5, -0.5) {};
		\node [style=empty] at (-5.2, -0.5) {$3$};

		\node [style=empty] (Y) at (-0.1, 1) {$Y$};
		\node [style=none] (Ytr) at (.25, 0.75) {};
		\node [style=none] (Ytl) at (-.4, 0.75) {};
		\node [style=none] (Ybr) at (.25, -0.75) {};
		\node [style=none] (Ybl) at (-.4, -0.75) {};

		\node [style=inputdot] (4) at (0, 0.5) {};
		\node [style=empty] at (0, 0.25) {$4$};
		\node [style=inputdot] (5) at (0, -0.5) {};
		\node [style=empty] at (0, -0.25) {$5$};		
		
		\node [style=empty] (Z) at (5, 1) {$Z$};
		\node [style=none] (Ztr) at (4.75, 0.75) {};
		\node [style=none] (Ztl) at (5.4, 0.75) {};
		\node [style=none] (Zbl) at (5.4, -0.75) {};
		\node [style=none] (Zbr) at (4.75, -0.75) {};

		\node [style=inputdot] (6) at (5, 0) {};
		\node [style=empty] at (5.2, 0) {$6$};	
		
	\end{pgfonlayer}
	\begin{pgfonlayer}{edgelayer}
		\draw [style=inarrow] (A) to (a);
		\draw [style=inarrow] (B) to (a);
		\draw [style=inarrow] (a) to (C);
		\draw [style=inarrow] (a) to (D);
		\draw [style=inarrow, bend left=30, looseness=1.00] (E) to (b);
		\draw [style=inarrow, bend left=30, looseness=1.00] (b) to (F);
		\draw [style=inarrow, bend left=30, looseness=1.00] (c) to (E);
		\draw [style=inarrow, bend left=30, looseness=1.00] (F) to (c);
		\draw [style=inputarrow] (1) to (A);
		\draw [style=inputarrow] (2) to (B);
		\draw [style=inputarrow] (3) to (B);
		\draw [style=inputarrow] (4) to (C);
		\draw [style=inputarrow] (5) to (D);
		\draw [style=inputarrow] (4) to (E);
		\draw [style=inputarrow] (5) to (E);
		\draw [style=inputarrow] (6) to (F);
		\draw [style=simple] (Xtl.center) to (Xtr.center);
		\draw [style=simple] (Xtr.center) to (Xbr.center);
		\draw [style=simple] (Xbr.center) to (Xbl.center);
		\draw [style=simple] (Xbl.center) to (Xtl.center);
		\draw [style=simple] (Ytl.center) to (Ytr.center);
		\draw [style=simple] (Ytr.center) to (Ybr.center);
		\draw [style=simple] (Ybr.center) to (Ybl.center);
		\draw [style=simple] (Ybl.center) to (Ytl.center);
		\draw [style=simple] (Ztl.center) to (Ztr.center);
		\draw [style=simple] (Ztr.center) to (Zbr.center);
		\draw [style=simple] (Zbr.center) to (Zbl.center);
		\draw [style=simple] (Zbl.center) to (Ztl.center);
	\end{pgfonlayer}
\end{tikzpicture}
\]
At this point, if we ignore the sets $X,Y,Z$, we have a new Petri net whose set of places is the disjoint union of those for $P$ and $Q$.  The second step is to identify a place of $P$ with a place of $Q$ whenever both are images of the same point in $Y$.  We can then stop drawing everything involving $Y$, and get an open Petri net $Q \circ P \maps X \to Z$:
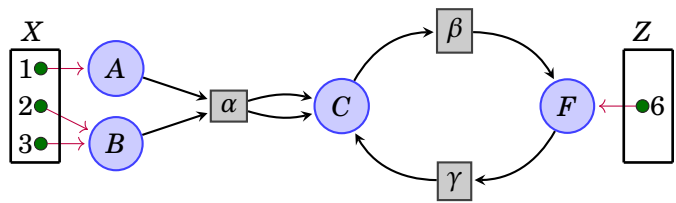
\[
\begin{tikzpicture}
	\begin{pgfonlayer}{nodelayer}
		\node [style=place] (A) at (-4, 0.5) {$A$};
		\node [style=place] (B) at (-4, -0.5) {$B$};;
             \node [style=transition] (a) at (-2.5, 0) {$\alpha$}; 
		\node [style = place] (E) at (-1, 0) {$C$};
		\node [style = place] (F) at (2,0) {$F$};

	     \node [style = transition] (b) at (.5, 1) {$\beta$};
		\node [style = transition] (c) at (.5, -1) {$\gamma$};
		
		\node [style=empty] (X) at (-5.1, 1) {$X$};
		\node [style=none] (Xtr) at (-4.75, 0.75) {};
		\node [style=none] (Xbr) at (-4.75, -0.75) {};
		\node [style=none] (Xtl) at (-5.4, 0.75) {};
             \node [style=none] (Xbl) at (-5.4, -0.75) {};
	
		\node [style=inputdot] (1) at (-5, 0.5) {};
		\node [style=empty] at (-5.2, 0.5) {$1$};
		\node [style=inputdot] (2) at (-5, 0) {};
		\node [style=empty] at (-5.2, 0) {$2$};
		\node [style=inputdot] (3) at (-5, -0.5) {};
		\node [style=empty] at (-5.2, -0.5) {$3$};	
		
		\node [style=empty] (Z) at (3, 1) {$Z$};
		\node [style=none] (Ztr) at (2.75, 0.75) {};
		\node [style=none] (Ztl) at (3.4, 0.75) {};
		\node [style=none] (Zbl) at (3.4, -0.75) {};
		\node [style=none] (Zbr) at (2.75, -0.75) {};

		\node [style=inputdot] (6) at (3, 0) {};
		\node [style=empty] at (3.2, 0) {$6$};	
		
	\end{pgfonlayer}
	\begin{pgfonlayer}{edgelayer}
		\draw [style=inarrow] (A) to (a);
		\draw [style=inarrow] (B) to (a);
	     \draw [style=inarrow, bend right=15, looseness=1.00] (a) to (E);
	     \draw [style=inarrow, bend left =15, looseness=1.00] (a) to (E);	
	     	\draw [style=inarrow, bend left=30, looseness=1.00] (E) to (b);
		\draw [style=inarrow, bend left=30, looseness=1.00] (b) to (F);
		\draw [style=inarrow, bend left=30, looseness=1.00] (c) to (E);
		\draw [style=inarrow, bend left=30, looseness=1.00] (F) to (c);	
		\draw [style=inputarrow] (1) to (A);
		\draw [style=inputarrow] (2) to (B);
		\draw [style=inputarrow] (3) to (B);
		\draw [style=inputarrow] (6) to (F);
		\draw [style=simple] (Xtl.center) to (Xtr.center);
		\draw [style=simple] (Xtr.center) to (Xbr.center);
		\draw [style=simple] (Xbr.center) to (Xbl.center);
		\draw [style=simple] (Xbl.center) to (Xtl.center);
		\draw [style=simple] (Ztl.center) to (Ztr.center);
		\draw [style=simple] (Ztr.center) to (Zbr.center);
		\draw [style=simple] (Zbr.center) to (Zbl.center);
		\draw [style=simple] (Zbl.center) to (Ztl.center);
	\end{pgfonlayer}
\end{tikzpicture}
\]The second fundamental operation of open networks comes from the morphisms between networks. These morphisms represent behavior preserving maps. When extended to open networks, these behavior preserving maps should preserve the inputs and outputs as well. For example, there is a morphism from an open Petri net $G \maps X_1 \to Y_1$:
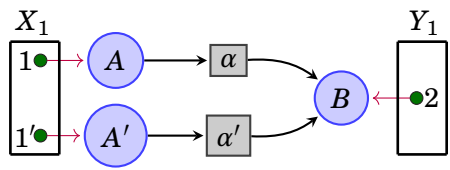
\[
\begin{tikzpicture}
	\begin{pgfonlayer}{nodelayer}

		\node [style = transition] (a) at (2.5, 0.5) {$\alpha$};
		\node [style = transition] (a') at (2.5, -0.5) {$\alpha'$};
		\node [style = place] (A) at (1, 0.5) {$A$};
		\node [style = place] (A') at (1, -0.5) {$A'$};
		\node [style = place] (B) at (4,0) {$B$};
	
		\node [style=empty] (Y) at (-0.1, 1) {$X_1$};
		\node [style=none] (Ytr) at (.25, 0.75) {};
		\node [style=none] (Ytl) at (-.4, 0.75) {};
		\node [style=none] (Ybr) at (.25, -0.75) {};
		\node [style=none] (Ybl) at (-.4, -0.75) {};

		\node [style=inputdot] (4) at (0, 0.5) {};
		\node [style=empty] at (-0.2, 0.5) {$1$};
		\node [style=inputdot] (5) at (0, -0.5) {};
		\node [style=empty] at (-0.2, -0.5) {$1'$};		
		
		\node [style=empty] (Z) at (5.1, 1) {$Y_1$};
		\node [style=none] (Ztr) at (4.75, 0.75) {};
		\node [style=none] (Ztl) at (5.4, 0.75) {};
		\node [style=none] (Zbl) at (5.4, -0.75) {};
		\node [style=none] (Zbr) at (4.75, -0.75) {};

		\node [style=inputdot] (6) at (5, 0) {};
		\node [style=empty] at (5.2, 0) {$2$};	
		
	\end{pgfonlayer}
	\begin{pgfonlayer}{edgelayer}

		\draw [style=inarrow] (A) to (a);
		\draw [style=inarrow, bend left=20, looseness=1.00] (a) to (B);
		\draw [style=inarrow] (A') to (a');
		\draw [style=inarrow, bend right=20, looseness=1.00] (a') to (B);
		\draw [style=inputarrow] (4) to (A);
		\draw [style=inputarrow] (5) to (A');
		\draw [style=inputarrow] (6) to (B);

		\draw [style=simple] (Ytl.center) to (Ytr.center);
		\draw [style=simple] (Ytr.center) to (Ybr.center);
		\draw [style=simple] (Ybr.center) to (Ybl.center);
		\draw [style=simple] (Ybl.center) to (Ytl.center);
		\draw [style=simple] (Ztl.center) to (Ztr.center);
		\draw [style=simple] (Ztr.center) to (Zbr.center);
		\draw [style=simple] (Zbr.center) to (Zbl.center);
		\draw [style=simple] (Zbl.center) to (Ztl.center);
	\end{pgfonlayer}
\end{tikzpicture}
\]
to an open Petri net $H \maps X_2 \to Y_2$
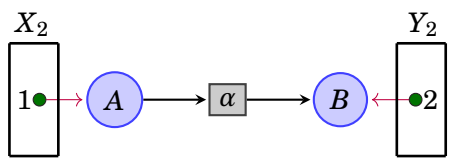
\[
\begin{tikzpicture}
	\begin{pgfonlayer}{nodelayer}

		\node [style = transition] (a) at (2.5, 0.0) {$\alpha$};
		\node [style = place] (A) at (1, 0.0) {$A$};
		\node [style = place] (B) at (4,0) {$B$};

		\node [style=empty] (Y) at (-0.1, 1) {$X_2$};
		\node [style=none] (Ytr) at (.25, 0.75) {};
		\node [style=none] (Ytl) at (-.4, 0.75) {};
		\node [style=none] (Ybr) at (.25, -0.75) {};
		\node [style=none] (Ybl) at (-.4, -0.75) {};

		\node [style=inputdot] (4) at (0, 0.0) {};
		\node [style=empty] at (-0.2, 0.0) {$1$};
		
		\node [style=empty] (Z) at (5.1, 1) {$Y_2$};
		\node [style=none] (Ztr) at (4.75, 0.75) {};
		\node [style=none] (Ztl) at (5.4, 0.75) {};
		\node [style=none] (Zbl) at (5.4, -0.75) {};
		\node [style=none] (Zbr) at (4.75, -0.75) {};

		\node [style=inputdot] (6) at (5, 0) {};
		\node [style=empty] at (5.2, 0) {$2$};	
		
	\end{pgfonlayer}
	\begin{pgfonlayer}{edgelayer}
		\draw [style=inarrow] (A) to (a);
		\draw [style=inarrow] (a) to (B);

		\draw [style=inputarrow] (4) to (A);
		\draw [style=inputarrow] (6) to (B);

		\draw [style=simple] (Ytl.center) to (Ytr.center);
		\draw [style=simple] (Ytr.center) to (Ybr.center);
		\draw [style=simple] (Ybr.center) to (Ybl.center);
		\draw [style=simple] (Ybl.center) to (Ytl.center);
		\draw [style=simple] (Ztl.center) to (Ztr.center);
		\draw [style=simple] (Ztr.center) to (Zbr.center);
		\draw [style=simple] (Zbr.center) to (Zbl.center);
		\draw [style=simple] (Zbl.center) to (Ztl.center);
	\end{pgfonlayer}
\end{tikzpicture}
\]
mapping both primed and unprimed symbols to unprimed ones. More precisely, this morphism of open Petri nets is a commutative diagram
\[\begin{tikzcd}
LX_1 \ar[r] \ar[d] & \ar[d] G & \ar[l] LY_1 \ar[d] \\
LX_2 \ar[r] & H & \ar[l] LY_2
\end{tikzcd} \]
in the category of Petri nets. We denote this morphism with the notation $\beta \maps G \Rightarrow H$. $\beta$ describes a process of ``simplifying" an open Petri net. There are also morphisms that include simple open Petri nets more complicated ones. For example, the above morphism of open Petri nets has a right inverse.

These two operations fit together into the structure of a double category. Double categories were introduced in the 1960s by Ehresmann \cite{Ehresmann63, Ehresmann65}.  More recently they have been used to study open dynamical systems \cite{Lerman,LS,N}, open electrical circuits and chemical reaction networks \cite{Bicategory}, open discrete-time Markov chains \cite{Panan}, coarse-graining for open continuous-time Markov chains \cite{BaezCourser}, and ``tile logic" for concurrency in computer science \cite{Tile}. Theorem \ref{Courser} constructs a double category where
\begin{itemize}
    \item objects are sets $X$,$Y$,$Z$,$\ldots$
    \item vertical morphisms are functions $f \maps X \to X'$,
    \item horizontal morphisms are open networks $G \maps X \to Y$,
    \item horizontal composition is the composition operation described above, and
    \item $2$-morphisms are the boundary preserving morphisms $\beta \maps G \Rightarrow H$ of open networks described above.
\end{itemize}
In Theorem \ref{openQnetdouble} we construct a double category $\Open(\Net{\Q})$ of the above form whose horizontal morphisms are open $\Q$-nets. In Theorem \ref{openRmatdouble} we construct an analogous double category $\Open(\RMat)$ for open $R$-matrices. Theorem \ref{openQnetdouble} has $\Open(\Petri)$, the double category of open Petri nets, as a special case. $\Open(\Petri)$ is a double category where objects are sets, vertical morphisms are functions, horizontal morphisms are open Petri nets, and $2$-morphisms are morphisms of open Petri nets. The axioms of a double category ensure that morphisms of open networks and composition of open networks are compatible. Besides composing open networks, we can also ``tensor" them via disjoint union: this describes networks being run in parallel rather than in series. The result is that the double category described above is upgraded to a symmetric monoidal double category. 

\noindent \subsubsection*{Composing Operational Semantics of Networks} The double categories of open networks constructed in this thesis describe a language for gluing smaller open networks into larger ones. The next step is to understand how the operational semantics of these networks can be applied to this language. We may extend the operational semantics functor for Petri nets
\[F \maps \Petri \to \CMC \]
 to a symmetric monoidal double functor
\[\Open (F) \maps \Open(\Petri) \to \Open(\CMC) \]
where $\Open(\CMC)$ is a double category whose horizontal morphisms are ``open commutative monoidal categories" $C \maps X \to Y$, i.e.\ cospans in $\CMC$ of the form
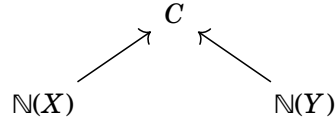
\[\begin{tikzcd}
& C & \\
\N(X) \ar[ur] & & \N(Y) \ar[ul]
\end{tikzcd}
\]
where $\N(X)$ and $\N(Y)$ are the discrete categories on the free commutative monoids on $X$ and $Y$. This double functor provides a compositional framework for composing the operational semantics of Petri nets. The key to this  double functor is that the functor $F$ preserves pushouts. Suppose a Petri net is decomposed into component open Petri nets
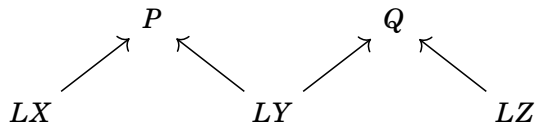
\[
\begin{tikzcd}
& P & & Q & \\
LX \ar[ur] & & LY \ar[ul] \ar[ur] & & \ar[ul] LZ
\end{tikzcd}
\]
$F$ may either be applied to their composite
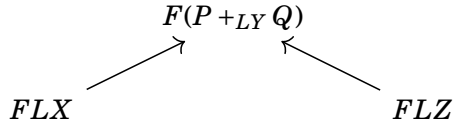
\[
\begin{tikzcd}
 & F(P+_{LY} Q) & \\
FLX \ar[ur] & & \ar[ul] FLZ
\end{tikzcd}
\]
or applied to each component and composed in $\Open(\CMC)$
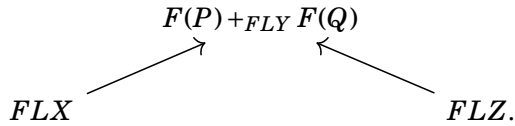
\[
\begin{tikzcd}
& F(P) +_{FLY} F(Q) & \\
FLX \ar[ur] & & FLZ. \ar[ul]
\end{tikzcd}
\]
Functoriality of $\Open(F)$ says that these two open commutative monoidal categories must be isomorphic; it therefore provides a compositionality relationship breaking down the operational semantics into smaller pieces. This is not a free lunch, the second pushout is taken in $\CMC$ which is constructed in a rather involved way. In general pushouts in $\CMC$ may be computed using Kelly's transfinite construction of free algebras \cite{kelly1980unified}. The idea behind this construction is that the pushout first takes the free commutative monoidal category 
\[F(UF(P)+_{UFLY} UF(Q)) \]
and then quotients away the redundant morphisms. The need for this second application of $F$ is clarified by the following example. Take $P$ to be this open Petri net:
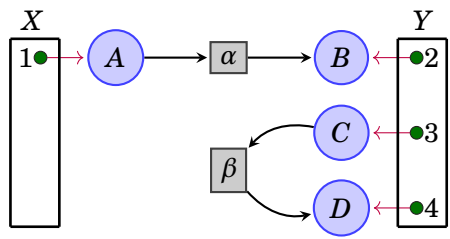
\[
\begin{tikzpicture}
	\begin{pgfonlayer}{nodelayer}
		\node [style=place] (A) at (-4, 1.5) {$A$};
		\node [style=place] (B) at (-1, 1.5) {$B$};
		\node [style=place] (C) at (-1, 0.5) {$C$};
		\node [style=place] (D) at (-1, -0.5) {$D$};
            \node [style=transition] (a) at (-2.5, 1.5) {$\alpha$}; 
             \node [style=transition] (b) at (-2.5, 0) {$\beta$}; 
		
		\node [style=empty] (X) at (-5.1, 2) {$X$};
		\node [style=none] (Xtr) at (-4.75, 1.75) {};
		\node [style=none] (Xbr) at (-4.75, -0.75) {};
		\node [style=none] (Xtl) at (-5.4, 1.75) {};
             \node [style=none] (Xbl) at (-5.4, -0.75) {};
	
		\node [style=inputdot] (1) at (-5, 1.5) {};
		\node [style=empty] at (-5.2, 1.5) {$1$};

		\node [style=empty] (Y) at (0.1, 2) {$Y$};
		\node [style=none] (Ytr) at (.4, 1.75) {};
		\node [style=none] (Ytl) at (-.25, 1.75) {};
		\node [style=none] (Ybr) at (.4, -0.75) {};
		\node [style=none] (Ybl) at (-.25, -0.75) {};

		\node [style=inputdot] (2) at (0, 1.5) {};
		\node [style=empty] at (0.2, 1.5) {$2$};
		\node [style=inputdot] (3) at (0, 0.5) {};
		\node [style=empty] at (0.2, 0.5) {$3$};
		\node [style=inputdot] (4) at (0, -0.5) {};
		\node [style=empty] at (0.2, -0.5) {$4$};		
		
		
	\end{pgfonlayer}
	\begin{pgfonlayer}{edgelayer}
		\draw [style=inarrow] (A) to (a);
		\draw [style=inarrow] (a) to (B);
		\draw [style=inarrow,bend right=30, looseness=1.00] (C) to (b);
		\draw [style=inarrow, bend right=30, looseness=1.00] (b) to (D);
		\draw [style=inputarrow] (1) to (A);
		\draw [style=inputarrow] (2) to (B);
		\draw [style=inputarrow] (3) to (C);
		\draw [style=inputarrow] (4) to (D);
	
		\draw [style=simple] (Xtl.center) to (Xtr.center);
		\draw [style=simple] (Xtr.center) to (Xbr.center);
		\draw [style=simple] (Xbr.center) to (Xbl.center);
		\draw [style=simple] (Xbl.center) to (Xtl.center);
		\draw [style=simple] (Ytl.center) to (Ytr.center);
		\draw [style=simple] (Ytr.center) to (Ybr.center);
		\draw [style=simple] (Ybr.center) to (Ybl.center);
		\draw [style=simple] (Ybl.center) to (Ytl.center);
	\end{pgfonlayer}
\end{tikzpicture}
\]
and take $Q$ to be this:
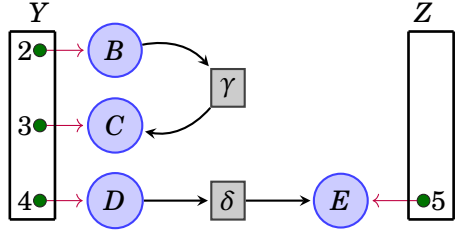
\[
\begin{tikzpicture}
	\begin{pgfonlayer}{nodelayer}
		\node [style=place] (B) at (-1, 1.5) {$B$};
		\node [style=place] (C) at (-1, 0.5) {$C$};
		\node [style=place] (D) at (-1, -0.5) {$D$};
		\node [style=place] (E) at (2, -0.5) {$E$};
            \node [style=transition] (c) at (0.5, 1) {$\gamma$}; 
             \node [style=transition] (d) at (0.5, -0.5) {$\delta$}; 
		
		\node [style=empty] (Z) at (3.1, 2) {$Z$};
		\node [style=none] (Ztr) at (2.9, 1.75) {};
		\node [style=none] (Zbr) at (2.9, -0.75) {};
		\node [style=none] (Ztl) at (3.5, 1.75) {};
             \node [style=none] (Zbl) at (3.5, -0.75) {};
	
		\node [style=inputdot] (5) at (3.1, -0.5) {};
		\node [style=empty] at (3.3, -0.5) {$5$};

		\node [style=empty] (Y) at (-2, 2) {$Y$};
		\node [style=none] (Ytr) at (-1.8, 1.75) {};
		\node [style=none] (Ytl) at (-2.4, 1.75) {};
		\node [style=none] (Ybr) at (-1.8, -0.75) {};
		\node [style=none] (Ybl) at (-2.4, -0.75) {};

		\node [style=inputdot] (2) at (-2, 1.5) {};
		\node [style=empty] at (-2.2, 1.5) {$2$};
		\node [style=inputdot] (3) at (-2, 0.5) {};
		\node [style=empty] at (-2.2, 0.5) {$3$};
		\node [style=inputdot] (4) at (-2, -0.5) {};
		\node [style=empty] at (-2.2, -0.5) {$4$};		
		
		
	\end{pgfonlayer}
	\begin{pgfonlayer}{edgelayer}
		\draw [style=inarrow,bend left=30, looseness=1.00] (B) to (c);
		\draw [style=inarrow, bend left=30, looseness=1.00] (c) to (C);
		\draw [style=inarrow] (D) to (d);
		\draw [style=inarrow] (d) to (E);
		\draw [style=inputarrow] (2) to (B);
		\draw [style=inputarrow] (3) to (C);
		\draw [style=inputarrow] (4) to (D);
		\draw [style=inputarrow] (5) to (E);
	
		\draw [style=simple] (Ytl.center) to (Ytr.center);
		\draw [style=simple] (Ytr.center) to (Ybr.center);
		\draw [style=simple] (Ybr.center) to (Ybl.center);
		\draw [style=simple] (Ybl.center) to (Ytl.center);
		\draw [style=simple] (Ztl.center) to (Ztr.center);
		\draw [style=simple] (Ztr.center) to (Zbr.center);
		\draw [style=simple] (Zbr.center) to (Zbl.center);
		\draw [style=simple] (Zbl.center) to (Ztl.center);
	\end{pgfonlayer}
\end{tikzpicture}
\]
Then their composite, $Q \circ P \maps X \to Z$, looks like this:
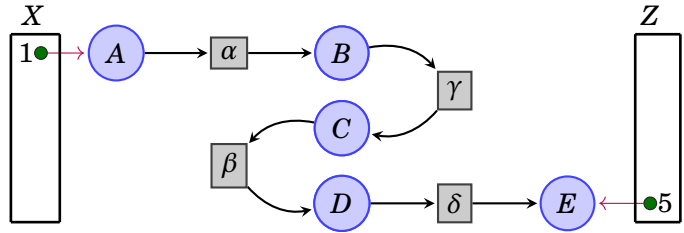
\[
\begin{tikzpicture}
	\begin{pgfonlayer}{nodelayer}
		\node [style=place] (A) at (-4, 1.5) {$A$};
		\node [style=place] (B) at (-1, 1.5) {$B$};
		\node [style=place] (C) at (-1, 0.5) {$C$};
		\node [style=place] (D) at (-1, -0.5) {$D$};
            \node [style=transition] (a) at (-2.5, 1.5) {$\alpha$}; 
             \node [style=transition] (b) at (-2.5, 0) {$\beta$}; 
		
		\node [style=empty] (X) at (-5.1, 2) {$X$};
		\node [style=none] (Xtr) at (-4.75, 1.75) {};
		\node [style=none] (Xbr) at (-4.75, -0.75) {};
		\node [style=none] (Xtl) at (-5.4, 1.75) {};
             \node [style=none] (Xbl) at (-5.4, -0.75) {};
	
		\node [style=inputdot] (1) at (-5, 1.5) {};
		\node [style=empty] at (-5.2, 1.5) {$1$};

		\node [style=place] (B) at (-1, 1.5) {$B$};
		\node [style=place] (C) at (-1, 0.5) {$C$};
		\node [style=place] (D) at (-1, -0.5) {$D$};
		\node [style=place] (E) at (2, -0.5) {$E$};
            \node [style=transition] (c) at (0.5, 1) {$\gamma$}; 
             \node [style=transition] (d) at (0.5, -0.5) {$\delta$}; 
		
		\node [style=empty] (Z) at (3.1, 2) {$Z$};
		\node [style=none] (Ztr) at (2.9, 1.75) {};
		\node [style=none] (Zbr) at (2.9, -0.75) {};
		\node [style=none] (Ztl) at (3.5, 1.75) {};
             \node [style=none] (Zbl) at (3.5, -0.75) {};
	
		\node [style=inputdot] (5) at (3.1, -0.5) {};
		\node [style=empty] at (3.3, -0.5) {$5$};

		
	\end{pgfonlayer}
	\begin{pgfonlayer}{edgelayer}
		\draw [style=inarrow] (A) to (a);
		\draw [style=inarrow] (a) to (B);
		\draw [style=inarrow,bend right=30, looseness=1.00] (C) to (b);
		\draw [style=inarrow, bend right=30, looseness=1.00] (b) to (D);
		\draw [style=inputarrow] (1) to (A);
	
		\draw [style=simple] (Xtl.center) to (Xtr.center);
		\draw [style=simple] (Xtr.center) to (Xbr.center);
		\draw [style=simple] (Xbr.center) to (Xbl.center);
		\draw [style=simple] (Xbl.center) to (Xtl.center);
		
		\draw [style=inarrow,bend left=30, looseness=1.00] (B) to (c);
		\draw [style=inarrow, bend left=30, looseness=1.00] (c) to (C);
		\draw [style=inarrow] (D) to (d);
		\draw [style=inarrow] (d) to (E);
		\draw [style=inputarrow] (5) to (E);
	
		\draw [style=simple] (Ztl.center) to (Ztr.center);
		\draw [style=simple] (Ztr.center) to (Zbr.center);
		\draw [style=simple] (Zbr.center) to (Zbl.center);
		\draw [style=simple] (Zbl.center) to (Ztl.center);
	\end{pgfonlayer}
\end{tikzpicture}
\]
This composite contains a morphism starting in $A$ and ending in $E$ which cannot be obtained from a firing sequence starting in $P$ and ending in $Q$. However, this morphism is contained in the free commutative monoidal category $F(UF(P) +_{UFLY} UF(Q))$ which accounts for all possible feedback loops and zig-zags between $P$ and $Q$. 

Often when studying open networks, we are less concerned with their internal workings than with the relationships they induce between their inputs and outputs. To represent this simplification we introduce the ``black-boxing" of an open network. For an open Petri net $P \maps X \to Y$, its black-boxing is a profunctor
\[\blacksquare(P) \maps \N(X) \times \N(Y) \to \Set \]
which sends a pair of markings $(x,y) \in \N(X) \times \N(Y)$ to the set of firing sequences in $P$ which start with $x$ and end with $Y$. Next we ask about the compositionality of black-boxing, i.e.\ how the profunctor composite $\blacksquare(Q) \circ \blacksquare(P)$ compares to $\blacksquare (Q \circ P)$. We show that black-boxing gives a lax monoidal double functor
\[\Open(\Petri) \to \Prof  \]
where $\Prof$ is a double category whose horizontal morphisms are profunctors. This double functor is only lax because the black-boxing $\blacksquare(Q \circ P)$ cannot be entirely reconstructed from the composite $\blacksquare(Q) \circ \blacksquare(P)$. As shown in the above ``zig-zag" example, $\blacksquare(Q \circ P)$ may contain morphisms which are not the composite of a morphism in $\blacksquare(P)$ with a morphism in $\blacksquare(Q)$. Next define a class of open networks that do not have this problematic behavior. We introduce ``functional open networks" based on functional Petri nets introduced by Zaitsev and Sleptsov \cite{zaitsev2005,zaitsev1997}. These are open networks for which every input is source and every output is a sink. We prove that black-boxing preserves composition on functional open networks. This gives a useful formula for composing operational semantics which can be turned into code as in \cite{compmarkov}.   

\noindent \subsubsection*{Outline of the Thesis} In Chapter \ref{graphs} we study the operational semantics of graphs and its extension to open graphs. In Section \ref{opsem} we construct the operational semantics of graphs. In Section \ref{opengraph} we define the symmetric monoidal double category of open graphs. In Section \ref{compopgraph} we extend the operational semantics of graphs to a symmetric monoidal double functor from the category of open graphs to the category of open categories. In Section \ref{blackboxgraph} we introduce the black-boxing of an open graph and show that it gives a lax double functor into the double category of profunctors. In Theorem \ref{functionalgraph} we define functional open graphs and show that black-boxing preserves their composition up to isomorphism.

In Chapter \ref{QNetchap} we define $\Q$-nets and construct their operational semantics. In Section \ref{defin} we review some definitions in Petri net theory. In Section \ref{QNet} we define $\Q$-nets and show how many existing variants of Petri nets and their relationships may be derived from this definition. In Section \ref{CMC} we construct an operational semantics adjunction for Petri nets. In Section \ref{gen} we generalize the previous section to an adjunction between $\Q$-nets and $\Q$-categories. This adjunction is factored into two parts. In Section \ref{freeQgraph} we construct the first part, turning $\Q$-nets into internal graphs, and in Section \ref{freecatinternal} we construct the second part, turning internal graphs into internal categories.

In Chapter \ref{openQNets} we study the compositionality of the operational semantics for $\Q$-nets. In Section \ref{openQnet} we define the double category of open $\Q$-nets. In Section \ref{compopseminternal} we extend the operational semantics of $\Q$-nets to a symmetric monoidal double functor from the double category of $\Q$-nets to the double category of $\Q$-categories. In Section \ref{functionalqnetsection} we define the black-boxing of an open $\Q$-net and show that black-boxing defines a lax double functor from the double category of open $\Q$-categories to a double category of profunctors. In Theorem \ref{blackboxfunctqnet} we show that black-boxing preserves composition of functional open $\Q$-nets up to isomorphism. 

In Chapter \ref{enriched} we construct an operational semantics for matrices valued in a commutative quantale $R$. In Section \ref{algpathproblem} we review the definition of $R$-matrices and the algebraic path problem. In Section \ref{freecatinternalsec} we construct an operational semantics adjunction between the category of $R$-matrices and the category of categories enriched in $R$. The left adjoint of this adjunction gives solutions to the algebraic path problem.

In Chapter \ref{compalgpathchap}, we explore how solutions of the algebraic path problem behave on open $R$-matrices. In Section \ref{openmat} we define the symmetric monoidal double category of open $R$-matrices. In Section \ref{openalg} we show how finding solutions to the algebraic path problem gives a symmetric monoidal double functor from the double category of open $R$-matrices to the double category of open $R$-categories. In Section \ref{functional} we define the black-boxing of an open $R$-matrix and show that it gives rise to a lax double functor. In Definition \ref{functmat} we define functional open $R$-matrices and in Theorem \ref{strict} we show that black-boxing preserves their composition strictly.

Lastly, in Appendix \ref{appendixdouble} we review the relevant definitions in the theory of double categories and in Appendix \ref{law} we review Lawvere theories.

\chapter{Compositionality of Graphs}\label{graphs}
This chapter serves as a blueprint for Chapters \ref{QNetchap}, \ref{openQNets}, \ref{enriched}, and \ref{openalg} by outlining the main results of this thesis in the case of ordinary graphs. In particular, in Section \ref{opsem}, we construct a well-known operational semantics of graphs in a way that lends itself to generalization to enriched and internal graphs. In Section \ref{opengraph}, we define ``open graphs", i.e.\ graphs equipped with input and output boundaries, and show that there is a symmetric monoidal double category $\Open(\Grph)$ whose horizontal morphisms are open graphs. In Section \ref{compopgraph} we show how the operational semantics of graphs can be extended a compositional setting, i.e.\ lifted to a double functor from a double category of open graphs to a double category of open categories. In Section \ref{blackboxgraph} we introduce the ``black-boxing" of an open graph. The black-boxing of an open graph is a profunctor which records the operational semantics of the open graph when restricted to the input and output boundaries. In Theorem \ref{blackbox} we prove that black-boxing lifts to a lax double functor from open graphs to a double category of profunctors. In Theorem \ref{functionalgraph} we identify a subclass of open graphs, called ``functional'', for which this double functor preserves composition up to isomorphism.
\section{Operational Semantics of Graphs}\label{opsem}
In this thesis we use the definition of graph preferred by category theorists: the edges have a direction and multiple edges are allowed between pairs of vertices.
\begin{defn}
A \define{graph} is a pair of functions
\[
\begin{tikzcd}
E \ar[r,shift left=.5ex,"s"] \ar[r,shift right=.5ex,"t",swap] & V.
\end{tikzcd}
\]
A \define{morphism of graphs} is a pair of functions $f \maps E \to E'$ and $g \maps V \to V'$ such that the following diagrams
\[
\begin{tikzcd}
E\ar[d,"f",swap] \ar[r,"s"] & V \ar[d,"g"] & E \ar[d,"f",swap] \ar[r,"t"] & V \ar[d,"g"]  \\
E' \ar[r,"s'",swap] & V' & E' \ar[r,"t'",swap] & V'
\end{tikzcd}
\]
commute. This defines a category $\Grph$ of graphs and their morphisms.
\end{defn}\noindent Let $C$ be the category generated by the graph
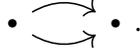
\[ \begin{tikzcd}\bullet \ar[r,bend left] \ar[r,bend right] & \bullet \end{tikzcd}. \]
Then $\Grph$ is the same as the functor category $[C,\Set]$. This fact implies the following proposition:
\begin{prop}\label{cocomplete}
$\Grph$ is complete and cocomplete with limits and colimits given pointwise in $\Set$.
\end{prop} \noindent Here a pointwise (co)limit of graphs is given by first taking the (co)limits of their underlying edges and vertices and extending the corresponding source and target maps to these new (co)limits. Paths in a graph can be constructed using pullbacks. Let $G$ be the graph\[\begin{tikzcd} E \ar[r,shift left=.5ex,"s"] \ar[r,shift right=.5ex,"t",swap] & V.\end{tikzcd}\]
If $G$ is placed next to itself:
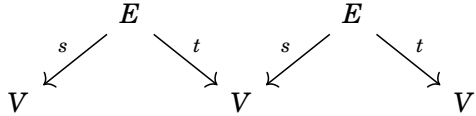
\[\begin{tikzcd}
 & \ar[dr,"t"] \ar[dl,"s",swap] E & & E \ar[dl,"s",swap] \ar[dr,"t"] & \\
 V & & V & & V
\end{tikzcd} \]
then the pullback of the center two functions is computed as
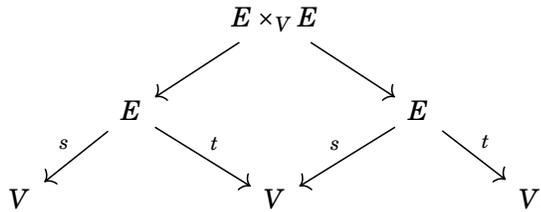
\[
\begin{tikzcd}
& & E \times_V E \ar[dr] \ar[dl] & & \\
 & \ar[dr,"t"] \ar[dl,"s",swap] E & & E \ar[dl,"s",swap] \ar[dr,"t"] & \\
 V & & V & & V.
\end{tikzcd}
\]
The outermost legs of this diagram form a graph whose edges are described explicitly as
\[\{ (e,e') \in E \times E \, |\, t(e) = s(e) \} \]
i.e.\ the paths of length $2$ in $G$. These pullbacks can be iterated $n$ times by placing $n$ copies of $G$ side by side and taking pullbacks until the outermost functions form a single span. Let $G^n$ be the graph formed by this
$n$-fold pullback. The edges of $G^n$ are given by elements of the set
\[\{(e_1,e_2,\ldots,e_n) \,| \, t(e_1)=s(e_2) \, , \, t(e_2)= s(e_3)\, , \ldots , t(e_{n-1})=s(e_n) \} \]
i.e.\ the set of paths of length $n$ in $G$. $G^0$ is defined as the span
\[\begin{tikzcd}V \ar[r,shift left=.5ex,"1_V"] \ar[r,shift right=.5ex,"1_V",swap] & V \end{tikzcd} \]
so that its edges are given by all paths of length $0$ in $G$. These powers are indeed iterated products in the following category:
\begin{defn}
Let $\Span(V,V)$ be the category where
\begin{itemize}
    \item objects are spans $V \xleftarrow{s} E \xrightarrow{t} V$ and
    \item morphisms are commutative diagrams
    
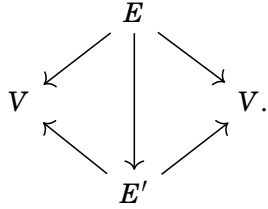
\[
    \begin{tikzcd}
    &\ar[dl] \ar[dr] E \ar[dd] & \\
    V & & V. \\
     & E' \ar[ul] \ar[ur] &
    \end{tikzcd}
    \]
\end{itemize}
\end{defn}
\noindent The product in this category is the pullback of spans defined above. For graphs $V \xleftarrow{s} E \xrightarrow{t} V$ and $V \xleftarrow{s'} E' \xrightarrow{t'} V$, their coproduct in $\Span(V,V)$ is the graph $V \xleftarrow{(s,s')} E + E' \xrightarrow{(t,t')} V$ where $(s,s')$ and $(t,t')$ represent the pairings of the source and target functions in each graph. To account for paths of any length in $G$, we must combine the graphs $G^n$ for all $n \geq 0$ using coproduct. The following proposition uses this idea to give the well-known free category construction:
\begin{prop}\label{freecat}
Let
$U \maps \Cat \to \Grph$ be the forgetful functor which sends a category to its underlying graph. Then $U$ has a left adjoint
\[F \maps \Grph \to \Cat \]
given by
\[F(G) = \sum_{n\geq 0} G^n \]
where products and sums are taken in the category of spans over the vertices of $G$.
\end{prop}
\begin{proof}A proof of this proposition can be found in many textbooks e.g.\ \cite{maclane}. Alternatively, the result can be proved in a similar way as in Sections \ref{freecatinternal} and \ref{freecatinternalsec}. In these cases as well as the above, the result follows from a general construction of free monoids over the relevant category of spans over a fixed object. Then the dependence on this object is removed using the Grothendieck construction.
\end{proof}
Because the morphisms of $F(G)$ are all paths in $G$, we borrow terminology from the theory of programming languages to call $F(G)$ the \define{operational semantics} of $G$. If $G$ represents a program where nodes are states and edges are ways of changing the state, then $F(G)$ is a category whose morphisms represent all possible runs of your program. Note that with this operational semantics, graphs are non-deterministic, i.e.\ for a given state there is in general more than one run of the program starting with that state. Non-determinism will be a feature of all the types of networks we consider in this thesis.
The operational semantics of a program encapsulates its behavior and is of use for model-checking and formal verification.

\section{Open Graphs}\label{opengraph}
To understand the compositionality of the operational semantics of graphs, we first need a paradigm where graphs are equipped with boundaries. These boundaries are represented by discrete graphs.
\begin{prop}
Let $R \maps \Grph \to \Set$ be the forgetful functor which sends a graph to its set of vertices and a function to its vertex component. Then $R$ has a left adjoint
\[L \maps \Set \to \Grph \]
which sends a set $X$ to the graph with no edges and $X$ as its set of vertices.
\end{prop}
\begin{defn}
A \define{open graph} is a cospan in $\Grph$ of the form
\[
\begin{tikzcd}
& G &\\
LX \ar[ur] & & \ar[ul] LY
\end{tikzcd}
\]
and is denoted by $G \maps X \to Y$.
\end{defn}
\noindent The idea is that the morphisms of the cospan designate certain vertices of $G$ to be either inputs or outputs. An open graph $G \maps X \to Y$ is represented by a picture like this:
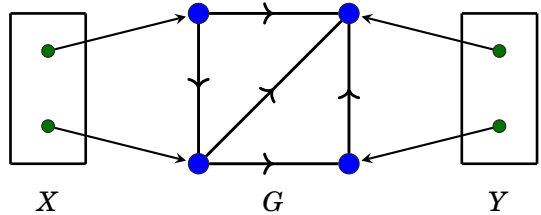
\[
\begin{tikzpicture}
	\begin{pgfonlayer}{nodelayer}
		\node [style=node] (0) at (-1, 1) {};
		\node [style=node] (1) at (-1, -1) {};
		\node [style=node] (2) at (1, 1) {};
		\node [style=node] (3) at (1, -1) {};
		\node [style=inputdot] (4) at (-3, 0.5) {};
		\node [style=inputdot] (5) at (-3, -0.5) {};
		\node [style=none] at (-3,-1.5) {$X$};
		\node [style=none] at (0,-1.5) {$G$};
		\node [style=none] at (3,-1.5) {$Y$};
		\node [style=inputdot] (6) at (3, 0.5) {};
		\node [style=inputdot] (7) at (3, -0.5) {};
		\node [style=none] (10) at (-3.5, 1) {};
		\node [style=none] (11) at (-2.5, 1) {};
		\node [style=none] (12) at (-3.5, -1) {};
		\node [style=none] (13) at (-2.5, -1) {};
		\node [style=none] (14) at (2.5, 1) {};
		\node [style=none] (15) at (3.5, 1) {};
		\node [style=none] (16) at (3.5, -1) {};
		\node [style=none] (17) at (2.5, -1) {};
	\end{pgfonlayer}
	\begin{pgfonlayer}{edgelayer}
		\draw [style=arrow] (0) to (2);
		\draw [style=arrow] (1) to (2);
		\draw [style=arrow] (0) to (1);
		\draw [style=arrow] (1) to (3);
		\draw [style=arrow] (3) to (2);
		\draw [style=inarrow] (4) to (0);
		\draw [style=inarrow] (5) to (1);
		\draw [style=inarrow] (6) to (2);
		\draw [style=inarrow] (7) to (3);
		\draw [style=simple] (10.center) to (12.center);
		\draw [style=simple] (10.center) to (11.center);
		\draw [style=simple] (11.center) to (13.center);
		\draw [style=simple] (12.center) to (13.center);
		\draw [style=simple] (14.center) to (17.center);
		\draw [style=simple] (17.center) to (16.center);
		\draw [style=simple] (16.center) to (15.center);
		\draw [style=simple] (14.center) to (15.center);

	\end{pgfonlayer}
\end{tikzpicture}
\]\noindent Given two composable open graphs
\[\begin{tikzcd}
& G & & H & \\
LX\ar[ur] & & \ar[ul] \ar[ur] LY & & \ar[ul] LZ
\end{tikzcd} \]
they are joined together using pushout
\[\begin{tikzcd}
&   & G+_{LY} H & &\\
& G \ar[ur] & & H \ar[ul]& \\
LX\ar[ur] & & \ar[ul] \ar[ur] LY & & \ar[ul] LZ
\end{tikzcd} \]
to obtain an open graph whose apex is a graph obtained by gluing $G$ and $H$ along their shared boundary. For example, if $G$ and $H$ are these open graphs:
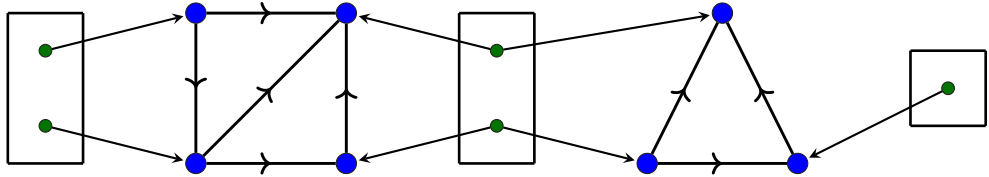
\[\begin{tikzpicture}
	\begin{pgfonlayer}{nodelayer}
		\node [style=node] (0) at (-4, 1) {};
		\node [style=node] (1) at (-4, -1) {};
		\node [style=node] (2) at (-2, 1) {};
		\node [style=node] (3) at (-2, -1) {};
		\node [style=inputdot] (4) at (-6, 0.5) {};
		\node [style=inputdot] (5) at (-6, -0.5) {};
		\node [style=inputdot] (6) at (0, 0.5) {};
		\node [style=inputdot] (7) at (0, -0.5) {};
		\node [style=none] (10) at (-6.5, 1) {};
		\node [style=none] (11) at (-5.5, 1) {};
		\node [style=none] (12) at (-6.5, -1) {};
		\node [style=none] (13) at (-5.5, -1) {};
		\node [style=none] (14) at (-0.5, 1) {};
		\node [style=none] (15) at (0.5, 1) {};
		\node [style=none] (16) at (0.5, -1) {};
		\node [style=none] (17) at (-0.5, -1) {};
		\node [style=node] (18) at (2, -1) {};
		\node [style=node] (19) at (4, -1) {};
		\node [style=node] (20) at (3, 1) {};
		\node [style=inputdot] (21) at (6, -0) {};
		\node [style=none] (22) at (5.5, 0.5) {};
		\node [style=none] (23) at (6.5, 0.5) {};
		\node [style=none] (24) at (5.5, -0.5) {};
		\node [style=none] (25) at (6.5, -0.5) {};
	\end{pgfonlayer}
	\begin{pgfonlayer}{edgelayer}
		\draw [style=arrow] (0) to (2);
		\draw [style=arrow] (1) to (2);
		\draw [style=arrow] (0) to (1);
		\draw [style=arrow] (1) to (3);
		\draw [style=arrow] (3) to (2);
		\draw [style=inarrow] (4) to (0);
		\draw [style=inarrow] (5) to (1);
		\draw [style=inarrow] (6) to (2);
		\draw [style=inarrow] (7) to (3);
		\draw [style=simple] (10.center) to (12.center);
		\draw [style=simple] (10.center) to (11.center);
		\draw [style=simple] (11.center) to (13.center);
		\draw [style=simple] (12.center) to (13.center);
		\draw [style=simple] (14.center) to (17.center);
		\draw [style=simple] (17.center) to (16.center);
		\draw [style=simple] (16.center) to (15.center);
		\draw [style=simple] (14.center) to (15.center);
		\draw [style=arrow] (18) to (20);
		\draw [style=arrow] (19) to (20);
		\draw [style=arrow] (18) to (19);
		\draw [style=simple] (22.center) to (24.center);
		\draw [style=simple] (22.center) to (23.center);
		\draw [style=simple] (23.center) to (25.center);
		\draw [style=simple] (24.center) to (25.center);
		\draw [style=inarrow] (6) to (20);
		\draw [style=inarrow] (7) to (18);
		\draw [style=inarrow] (21) to (19);
	\end{pgfonlayer}
\end{tikzpicture} \]
then their pushout is the open graph
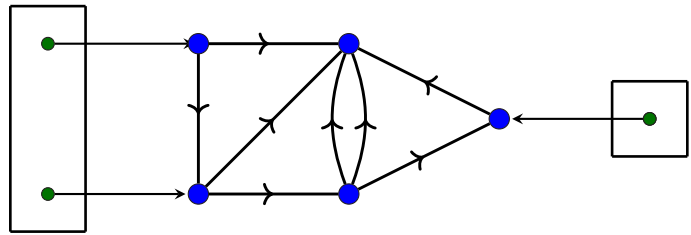
\[\begin{tikzpicture}
	\begin{pgfonlayer}{nodelayer}
		\node [style=none] (0) at (-1, 1) {};
		\node [style=node] (1) at (-1, 1) {};
		\node [style=node] (2) at (-1, -1) {};
		\node [style=node] (3) at (1, -1) {};
		\node [style=node] (4) at (1, 1) {};
		\node [style=node] (5) at (3, -0) {};
		\node [style=inputdot] (6) at (-3, 1) {};
		\node [style=inputdot] (7) at (-3, -1) {};
		\node [style=inputdot] (8) at (5, -0) {};
		\node [style=inputdot] (9) at (5, -0) {};
		\node [style=none] (10) at (-3.5, 1.5) {};
		\node [style=none] (11) at (-2.5, 1.5) {};
		\node [style=none] (12) at (-3.5, -1.5) {};
		\node [style=none] (13) at (-2.5, -1.5) {};
		\node [style=none] (14) at (4.5, 0.5) {};
		\node [style=none] (15) at (5.5, 0.5) {};
		\node [style=none] (16) at (5.5, -0.5) {};
		\node [style=none] (17) at (4.5, -0.5) {};
	\end{pgfonlayer}
	\begin{pgfonlayer}{edgelayer}
		\draw [style=inarrow] (6) to (0);
		\draw [style=inarrow] (7) to (2);
		\draw [style=inarrow] (8) to (5);
		\draw [style=arrow] (0) to (2);
		\draw [style=arrow] (2) to (3);
		\draw [style=arrow] (3) to[bend left=20] (4);
		\draw [style=arrow] (3) to[bend right=20] (4);
		\draw [style=arrow] (0) to (4);
		\draw [style=arrow] (2) to (4);
		\draw [style=arrow] (3) to (5);
		\draw [style=arrow] (5) to (4);
		\draw [style=simple] (10.center) to (12.center);
		\draw [style=simple] (12.center) to (13.center);
		\draw [style=simple] (13.center) to (11.center);
		\draw [style=simple] (10.center) to (11.center);
		\draw [style=simple] (17.center) to (14.center);
		\draw [style=simple] (14.center) to (15.center);
		\draw [style=simple] (15.center) to (16.center);
		\draw [style=simple] (16.center) to (17.center);
	\end{pgfonlayer}
\end{tikzpicture}\]
obtained by identifying all vertices which are mapped to by a common point in their shared boundary. This gluing operation forms the horizontal composition of a double category.
The formalisms developed by Fong \cite{fong2016algebra} and Courser \cite{CourserThesis} define a syntax for gluing of open systems using cospans. We will use the following result from \cite{CourserThesis} to construct a syntax for open graphs as well as the other networks considered in this thesis. Since this is a symmetric monoidal double category, it involves quite a lot of structure. The definition of symmetric monoidal double category can be found in Appendix \ref{appendixdouble}.

\begin{lem}[Courser]\label{Courser}
Let $\sA$ be a category with finite coproducts and $\X$ be a category with finite colimits. Given a left adjoint $L \maps \sA\to \X$, there exists a unique symmetric monoidal double category $_L \Cospan(\X)$,  such that:
	\begin{itemize}
		\item objects are objects of $\sA$,
		\item vertical 1-morphisms are morphisms of $\sA$,
		\item a horizontal 1-cell from $a \in \sA$ to $b\in \sA$ is a cospan in $\X$ of this form:
			\[
			\begin{tikzpicture}[scale=1.5]
			\node (A) at (0,0) {$La$};
			\node (B) at (1,0) {$x$};
			\node (C) at (2,0) {$Lb$};
			\path[->,font=\scriptsize,>=angle 90]
			(A) edge node[above,left]{} (B)
			(C)edge node[above]{}(B);
			\end{tikzpicture}
			\]
		\item{a 2-morphism is a commutative diagram in $\X$ of this form:
			\[
			\begin{tikzpicture}[scale=1.5]
			\node (E) at (3,0) {$La$};
			\node (F) at (5,0) {$Lb$};
			\node (G) at (4,0) {$x$};
			\node (E') at (3,-1) {$Lc$};
			\node (F') at (5,-1) {$Ld$.};
			\node (G') at (4,-1) {$y$};
			\path[->,font=\scriptsize,>=angle 90]
			(F) edge node[above]{} (G)
			(E) edge node[left]{$Lf$} (E')
			(F) edge node[right]{$Lg$} (F')
			(G) edge node[left]{$h$} (G')
			(E) edge node[above]{} (G)
			(E') edge node[below]{} (G')
			(F') edge node[below]{} (G');
			\end{tikzpicture}
			\]
			}
			\end{itemize}
Composition of vertical 1-morphisms is composition in $\sA$.   Composition of horizontal
1-cells is composition of cospans in $\X$ via pushout: given horizontal 1-cells
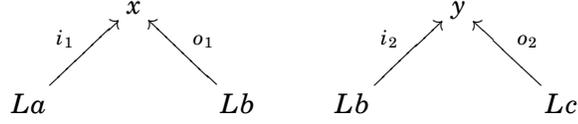
\[ \xymatrix{ & x & & & y & \\
	La \ar[ur]^{i_1} & & Lb \ar[ul]_{o_1} & Lb \ar[ur]^{i_2} & & Lc \ar[ul]_{o_2} }\]
their composite is this cospan from $La$ to $Lc$:
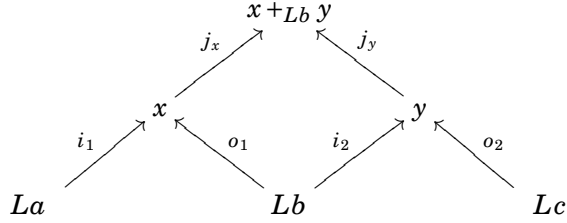
\[ \xymatrix{
	&   & x+_{Lb} y  &  & \\
	& x \ar[ur]^{j_x} &  & y \ar[ul]_{j_y} & \\
	La \quad \ar[ur]^{i_1} & & Lb \ar[ul]_{o_1}  \ar[ur]^{i_2} & & \quad Lc \ar[ul]_{o_2} }\]
where the diamond is a pushout square.  The horizontal composite of 2-morphisms
 
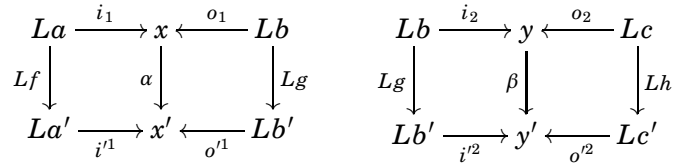
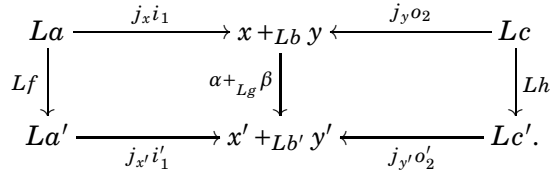
\[ \xymatrix{	La \ar[r]^{i_1} \ar[d]_{Lf} & x \ar[d]_{\alpha} & Lb \ar[l]_{o_1} \ar[d]^{Lg} \\
			  	La' \ar[r]_{i'^1}  & x'  & Lb' \ar[l]^{o'^1} }
\qquad
 \xymatrix{	Lb \ar[r]^{i_2} \ar[d]_{Lg} & y \ar[d]_{\beta} & Lc \ar[l]_{o_2} \ar[d]^{Lh} \\
			  	Lb' \ar[r]_{i'^2}  & y'  & Lc' \ar[l]^{o'^2} } \]
is given by
 \[ \xymatrix{	La \ar[rr]^-{j_x i_1} \ar[d]_{Lf} && x+_{Lb} y \ar[d]_{\alpha+_{{}_{Lg}} \beta}
 && Lc \ar[ll]_-{j_{y} o_2} \ar[d]^{Lh} \\
La' \ar[rr]_-{j_{x'} i'_1}  && x'+_{Lb'} y'  && Lc'. \ar[ll]^-{j_{y'} o'_2} } \]
The vertical composite of 2-morphisms
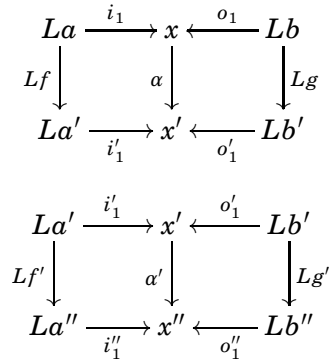
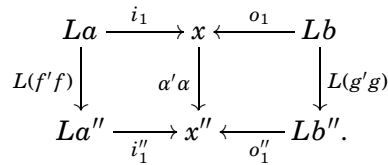
\[     \xymatrix{	La \ar[r]^{i_1} \ar[d]_{Lf} & x \ar[d]_{\alpha} & Lb \ar[l]_{o_1} \ar[d]^{Lg} \\
			  	La' \ar[r]_{i'_1}  & x'  & Lb' \ar[l]^{o'_1} } \]
\[   \xymatrix{	La' \ar[r]^{i'_1} \ar[d]_{Lf'} & x' \ar[d]_{\alpha'} &
 Lb' \ar[l]_{o'_1} \ar[d]^{Lg'} \\
			  	La'' \ar[r]_{i''_1}  & x''  & Lb'' \ar[l]^{o''_1} } \]
is given by
\[     \xymatrix{	La \ar[r]^{i_1} \ar[d]_{L(f'f)} & x \ar[d]_{\alpha'\alpha} & Lb \ar[l]_{o_1} \ar[d]^{L(g'g)} \\
			  	La'' \ar[r]_{i''_1}  & x''  & Lb''. \ar[l]^{o''_1} } \]
The tensor product is defined using
chosen coproducts in $\sA$ and $\X$.  Thus, the tensor product of two
objects $a_1$ and $a_2$ is $a_1+a_2$, the tensor product of two vertical 1-morphisms
 \[   \xymatrix{	a_1 \ar[d]_{f_1} \\  b_1 } \qquad \qquad  \xymatrix{ a_2 \ar[d]_{f_2} \\ b_2 } \]
 is
 \[ \xymatrix{	a_1+a_2 \ar[d]_{f_1+f_2} \\ b_1+b_2, } \]
the tensor product of two horizontal 1-cells
\[ \xymatrix{ 	La_1  \ar[r]^{i_1} & x_1 & Lb_1 \ar[l]_{o_1} }   \qquad \qquad
\xymatrix {	La_2  \ar[r]^{i_2} & x_2 & Lb_2 \ar[l]_{o_2} } \]
is
\[ \xymatrix{
	L(a_1 +a_2) \ar[rr]^{i_1+i_2} &&  x_1+ x_2 && L(b_1 + b_2), \ar[ll]_{o_1 + o_2} } \]
and the tensor product of two 2-morphisms
 \[ \xymatrix{	La_1 \ar[r]^{i_1} \ar[d]_{Lf_1} & x_1 \ar[d]_{\alpha_1} & Lb_1 \ar[l]_{o_1} \ar[d]^{Lg_1} \\
			  	La'_1 \ar[r]_{i'_1}  & x'_1  & Lb'_1 \ar[l]^{o'_1} } \qquad
 \xymatrix{	La_2 \ar[r]^{i_2} \ar[d]_{Lf_2} & x_2 \ar[d]_{\alpha_2} & Lb_2 \ar[l]_{o_2} \ar[d]^{Lg_2} \\
			  	La'_2 \ar[r]_{i'_2}  & x'_2  & Lb'_2 \ar[l]^{o'_2} } \]
is
\[   \xymatrix{	L(a_1+a_2) \ar[rr]^{i_1+i_2} \ar[d]_{L(f_1+f_2)} && x_1+x_2 \ar[d]_{\alpha_1+\alpha_2} && L(b_1+b_2) \ar[ll]_{o_1 + o_2} \ar[d]^{L(g_1+g_2)} \\
			  	L(a'_1+a'_2) \ar[rr]_{i'_1+i'_2} && x'_1+x'_2  && L(b'_1+b'_2). \ar[ll]^{o'_1+o'_2} } \]
The units for these tensor products are taken to be initial objects, and the symmetry is defined using the canonical isomorphisms $a + b \cong b + a$.
\end{lem}
\begin{proof}  This was proved by Baez and Courser \cite[Thm. 3.9]{structured}.  Note that we are abusing language slightly above.   We must choose a specific coproduct for each pair of
objects in $\X$ and $\sA$ to give $_L \Cospan(X)$ its tensor product.   Given
morphisms $i_1 \maps La_1 \to x_1$ and $i_2 \maps La_2 \to x_2$, their coproduct is really a morphism $i_1 + i_2 \maps La_1 + La_2 \to x_1 + x_2$ between
these chosen coproducts.  But since $L$ preserves coproducts, we can compose this
morphism with the canonical isomorphism $L(a_1 + a_2) \cong La_1 + La_2$ to obtain the morphism that we call $i_1 + i_2 \maps L(a_1 + a_2) \to x_1 + x_2$ above.
\end{proof}\noindent Now we apply this lemma to the left adjoint $L$ defined above.
\begin{thm}\label{open}
There is a symmetric monoidal double category $\Open(\Grph)$ where
\begin{itemize}
    \item objects are sets $X$,$Y$,$Z \ldots$
    \item vertical morphisms are functions $f: X \to Y$,
    \item a horizontal morphism $G \maps X \to Y$ is an open graph
    
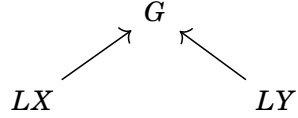
\[
    \begin{tikzcd}
     & G & \\\
    LX \ar[ur] & & LY \ar[ul]
    \end{tikzcd}
    \]
    \item vertical 2-morphisms are commutative rectangles
    
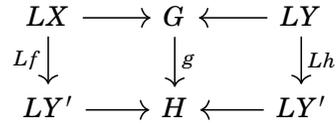
\[
    \begin{tikzcd}
    LX\ar[d,"Lf",swap] \ar[r] & G\ar[d,"g"] & \ar[l] \ar[d,"Lh"] LY \\
    LY' \ar[r] & H& \ar[l] LY'
    \end{tikzcd}
    \]
    \item vertical composition is ordinary composition of functions,
    \item and horizontal composition of an open graph $G \maps X \to Y$ and an open graph $H \maps Y \to Z$ is given by their pushout.
\end{itemize}
The symmetric monoidal structure is given by
\begin{itemize}
    \item coproducts in $\Set$ on objects and vertical morphisms,
    \item pointwise coproducts on horizontal morphisms i.e.\ for open graphs,
    
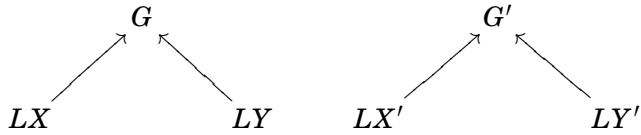
\[
    \xymatrix{ & G& & & G'& \\
	LX \ar[ur] & & LY \ar[ul] & LX' \ar[ur] & & LY' \ar[ul] }
	\]
	their coproduct is
	
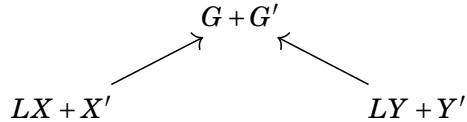
\[\begin{tikzcd}
	& G+ G'& \\
	L{X+ X'} \ar[ur] & & \ar[ul] L{Y + Y'}
	\end{tikzcd} \]
	and pointise coproduct for two vertical 2-morphisms i.e.\ for vertical 2-morphisms,
	
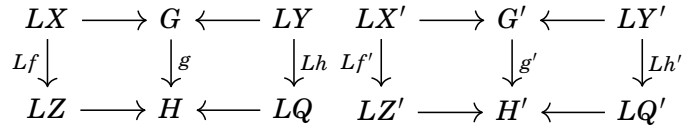
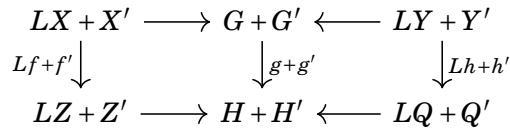
\[
	\begin{tikzcd}
    LX\ar[d,"Lf",swap] \ar[r] & G\ar[d,"g"] & \ar[l] \ar[d,"Lh"] LY \\
    LZ \ar[r] & H& \ar[l] LQ
    \end{tikzcd}
    \begin{tikzcd}
    LX'\ar[d,"Lf'",swap] \ar[r] & G'\ar[d,"g'"] & \ar[l] \ar[d,"Lh'"] LY' \\
    LZ' \ar[r] & H'& \ar[l] LQ'
    \end{tikzcd}
    \]
    their coproduct is
    \[\begin{tikzcd}
    L{X+ X'}\ar[d,"L{f+ f'}",swap] \ar[r] & G+ G'\ar[d,"g+ g'"] & \ar[l] \ar[d,"L{h+ h'}"] L{Y+ Y'} \\
    L{Z+ Z'} \ar[r] & H+ H'& \ar[l] L{Q+ Q'}
    \end{tikzcd} \]
\end{itemize}
\end{thm}

\begin{proof}

Lemma \ref{Courser} constructs this symmetric monoidal double category given that
\begin{itemize}
    \item $\Grph$ has coproducts and pushouts,
    \item and $L \maps \Set \to \Grph$ preserves pushouts and coproducts.
\end{itemize}The first point follows from Proposition \ref{cocomplete} and the second point is true because $L$ is a left adjoint.
\end{proof}

\section{Compositional Operational Semantics of Graphs}\label{compopgraph}
In this section we show how the operational semantics of graphs
\[F \maps \Grph \to \Cat \]
defined in Section \ref{opsem} is extended to open graphs. Categories can also be made open in a similar way as graphs.
\begin{defn}
An \define{open category} is a cospan in $\Cat$ of the form
\[\begin{tikzcd}
 & C & \\
FLX \ar[ur] & & \ar[ul] FLY
\end{tikzcd} \]
\end{defn}\noindent Categories are also glued together using pushout.

Pushout of open categories also gives the horizontal composition of a symmetric monoidal double category.
\begin{thm}\label{opencat}
There is a symmetric monoidal double category $\Open(\Cat)$ where
\begin{itemize}
    \item objects are sets $X$,$Y$,$Z \ldots$
    \item vertical morphisms are functions $f: X \to Y$,
    \item a horizontal morphism $C \maps X \to Y$ is an open category
    \[
    \begin{tikzcd}
     & C & \\\
    FLX \ar[ur] & & FLY \ar[ul]
    \end{tikzcd}
    \]
    \item vertical 2-morphisms are commutative rectangles
    \[
    \begin{tikzcd}
    FLX\ar[d,"FLf",swap] \ar[r] & C\ar[d,"g"] & \ar[l] \ar[d,"FLh"] FLY \\
    FLY' \ar[r] & D& \ar[l] FLY'
    \end{tikzcd}
    \]
    \item vertical composition is ordinary composition of functions,
    \item and horizontal composition of an open category $C \maps X \to Y$ and an open graph $D \maps Y \to Z$ is given by their pushout.
\end{itemize}
The symmetric monoidal structure is given by
\begin{itemize}
    \item coproducts in $\Set$ on objects and vertical morphisms,
    \item pointwise coproducts on horizontal morphisms i.e.\ for open categories,
    \[
    \xymatrix{ & C& & & C'& \\
	FLX \ar[ur] & & FLY \ar[ul] & FLX' \ar[ur] & & FLY' \ar[ul] }
	\]
	their coproduct is
	\[\begin{tikzcd}
	& C+ C'& \\
	FL(X+ X') \ar[ur] & & \ar[ul] FL(Y + Y')
	\end{tikzcd} \]
	and pointise coproduct for two vertical 2-morphisms i.e.\ for vertical 2-morphisms,
	\[
	\begin{tikzcd}
    FLX\ar[d,"FLf",swap] \ar[r] & C\ar[d,"g"] & \ar[l] \ar[d,"FLh"] FLY \\
    FLZ \ar[r] & D& \ar[l] FLQ
    \end{tikzcd}
    \begin{tikzcd}
    FLX'\ar[d,"FLf'",swap] \ar[r] & C'\ar[d,"g'"] & \ar[l] \ar[d,"FLh'"] FLY' \\
    FLZ' \ar[r] & D'& \ar[l] FLQ'
    \end{tikzcd}
    \]
    their coproduct is
    \[\begin{tikzcd}
    FL(X+ X')\ar[d,"FL{f+ f'}",swap] \ar[r] & C+ C'\ar[d,"g+ g'"] & \ar[l] \ar[d,"FL{h+ h'}"] FL{Y+ Y'} \\
    FL{Z+ Z'} \ar[r] & D+ D'& \ar[l] FL{Q+ Q'}
    \end{tikzcd} \]
\end{itemize}
\end{thm}
\begin{proof}
This double category is constructed by applying Lemma \ref{Courser} to the composite left adjoint
\[\Set \xrightarrow{L} \Grph \xrightarrow{F} \Cat. \qedhere \]
\end{proof}Theorem 4.3 of \cite{structured} shows how a commutative square of left adjoint functors lifts to a symmetric monoidal double functor. In this thesis we require a weakening of this result.
\begin{lem}\label{openfunctoriality}
Let
\[
\begin{tikzcd}
C \ar[r,"F"] & D \\
A \ar[r,"F_0",swap] \ar[u,"L"] & A' \ar[u,"L'",swap]
\end{tikzcd}
\]
be a diagram commuting up to natural isomorphism where $L$, $L'$, and $F_0$ preserve finite colimits. Then there is a lax symmetric monoidal lax double functor
\[ \Open(F) \maps \Open(C) \to \Open(D)\]
which is given by $F_0$ on objects and morphisms of $A$ and is given by pointwise application of $F$ on horizontal morphisms and $2$-cells. Explicitly, $0$-cells and vertical morphisms are mapped as follows:
\[\begin{tikzcd}
   A \ar[d,"f"{name=L}] & & F_0(A) \ar[d,"F_0(f)"{name=R}]  \\
   B & & F_0(B) \ar[mapsto, from=L,to=R,shorten <=7ex, shorten >=7ex]
\end{tikzcd}
\]\noindent horizontal morphisms are mapped as follows:
\[\begin{tikzcd}
	LX \ar[r] & G & \ar[l] LY &  \mapsto
	& L'F_0X \ar[r] & FG &\ar[l] L'F_0Y
\end{tikzcd}\]
and vertical $2$-cells are mapped as follows:
\[\begin{tikzcd}
	LX & G & LY &  & L'F_0X & FG & L'F_0Y \\
	{LX'} & H & {LY'} && {L'F_0X'} & {FG'} & {L'F_0Y'}
	\arrow[from=1-1, to=1-2]
	\arrow[from=1-3, to=1-2]
	\arrow[from=2-1, to=2-2]
	\arrow[from=2-3, to=2-2]
	\arrow["Lf"', from=1-1, to=2-1]
	\arrow["g", from=1-2, to=2-2]
	\arrow["Lh", from=1-3, to=2-3,""{name=L}]
	\arrow[from=1-5, to=1-6]
	\arrow[from=1-7, to=1-6]
	\arrow[from=2-5, to=2-6]
	\arrow[from=2-7, to=2-6]
	\arrow["Fg", from=1-6, to=2-6]
	\arrow["L'F_0f"', from=1-5, to=2-5,""{name=R}]
	\arrow["L'F_0h", from=1-7, to=2-7]
	\ar[mapsto,from=L, to=R, shorten <=7ex, shorten >=7ex]
\end{tikzcd}\]
\end{lem}
\begin{proof}
Theorem 4.3 of \cite{structured} supplies a symmetric monoidal double functor for the above square when all functors preserve finite colimits. However, a lax symmetric monoidal lax double functor (as opposed to pseudo) can nevertheless be constructed when $F$ is an arbitrary functor. The laxator of composition for this double functor
\[ \begin{tikzcd}
& F(G) +_{L'Y} F(H) \ar[dd,"\phi_{GH}"]& \\
L'X \ar[ur] \ar[dr] & & L'Z \ar[dl] \ar[ul]\\
& F(G+_{LY} H) &
\end{tikzcd}\]\noindent is induced by the universal property of pushout on the morphisms $F(i) \maps F(G) \to F(G+_{LY} H)$ and $F(j) \maps F(H) \to F(G+_{LY} H)$ where $i \maps G \hookrightarrow G+_{LY} H$ and $j \maps H \hookrightarrow G+_{LY} H$ are the canonical inclusions. Similarly, the monoidal comparison
\[ \begin{tikzcd}
 & F(G) + F(H) \ar[dd,"\psi_{GH}"]& \\
L'(X+X') \ar[ur] \ar[dr] & & L'(Y+Y') \ar[dl] \ar[ul]\\
& F(G+ H) &
\end{tikzcd}\]
is induced by the universal property of coproduct applied to the morphisms $F(i) \maps F(G) \to F(G+ H)$ and $F(j) \maps F(H) \to F(G+H)$ where $i$ and $j$ are the canonical inclusions into the coproduct. Verifying that this structure satisfies the axioms of a symmetric monoidal double functor follows the proof of Theorem 4.3 in \cite{structured} very closely.\end{proof}
The discrete graph functor and the free category functor assemble into a diagram
\[\begin{tikzcd}
\Grph\ar[r,"F"]  & \Cat \\
\Set \ar[r,equals] \ar[u,"L"]& \Set \ar[u,"FL",swap]
\end{tikzcd} \]
of commuting left adjoint functors. The following theorem is given by applying Lemma \ref{openfunctoriality} to this diagram. Because $F$ preserves colimits, horizontal composition is preserved up to isomorphism.
\begin{thm}
There is a symmetric monoidal double functor \[\Open(F) \maps \Open(\Grph) \to \Open(\Cat)\]
which is the identity on objects and functions and is given by pointwise application of $F$ on horizontal morphisms and 2-morphisms. Explicitly, the horizontal morphisms and two morphisms are mapped as follows:
\[\begin{tikzcd}
	LX \ar[r] & G & \ar[l] LY &  \mapsto
	& FLX \ar[r] & FG &\ar[l] FLY
]
\end{tikzcd}\]
\vspace{.5em}
\[\begin{tikzcd}
	LX & G & LY &  & FLX & FG & FLY \\
	{LX'} & H & {LY'} && {FLX'} & {FG'} & {FLY'}
	\arrow[from=1-1, to=1-2]
	\arrow[from=1-3, to=1-2]
	\arrow[from=2-1, to=2-2]
	\arrow[from=2-3, to=2-2]
	\arrow["Lf"', from=1-1, to=2-1]
	\arrow["g", from=1-2, to=2-2]
	\arrow["Lh", from=1-3, to=2-3,""{name=L}]
	\arrow[from=1-5, to=1-6]
	\arrow[from=1-7, to=1-6]
	\arrow[from=2-5, to=2-6]
	\arrow[from=2-7, to=2-6]
	\arrow["Fg", from=1-6, to=2-6]
	\arrow["FLf"', from=1-5, to=2-5,""{name=R}]
	\arrow["FLh", from=1-7, to=2-7]
	\ar[mapsto,from=L, to=R, shorten <=7ex, shorten >=7ex]
\end{tikzcd}\]
\end{thm}

This symmetric monoidal double functor gives a compositional framework for building the operational semantics of a graph recursively. A graph is decomposed into component open graphs, the operational semantics on each of these components is computed and then they are joined together using pushout to form the operational semantics of the total graph.
\begin{expl}\label{trainroutes}
Consider this non-exhaustive map of train routes in Southern California regarded as an open graph $G \maps \phi \to \phi$:
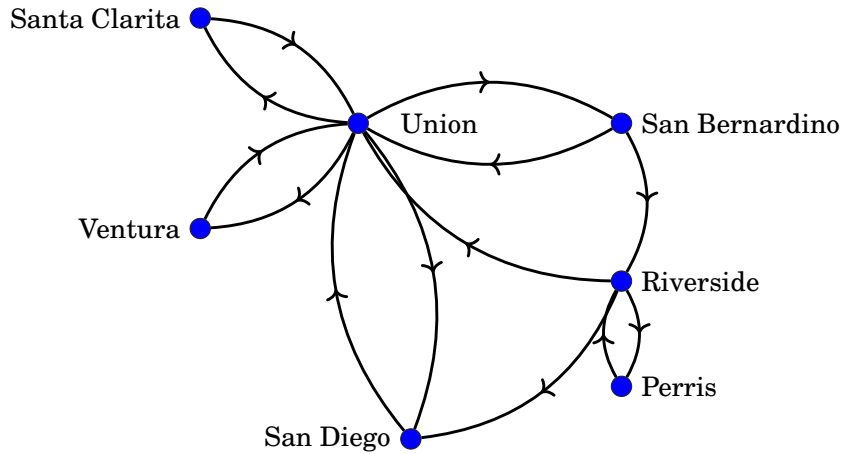
\[
\begin{tikzpicture}[scale=.7]
	\begin{pgfonlayer}{nodelayer}
		\node [style=node,label=left:{Santa Clarita}] (0) at (-3, 2) {};
		\node [style=node,label=left:{Ventura}] (1) at (-3, -2) {};
		\node [style=node,label={[label distance=.3cm]0:Union}] (3) at (0, -0) {};
		\node [style=node,label=right:{San Bernardino}] (4) at (5, -0) {};
		\node [style=node,label=right:Riverside] (5) at (5, -3) {};
		\node [style=node,label=right:Perris] (6) at (5, -5) {};
		\node [style=node,label=left:{San Diego}] (7) at (1, -6) {};

	\end{pgfonlayer}
	\begin{pgfonlayer}{edgelayer}
		\draw [style=arrow] (0) to[bend left=30] (3);
		\draw [style=arrow] (1) to[bend left=30] (3);
		\draw [style=arrow] (3) to[bend left=30] (1);
		\draw [style=arrow] (3) to[bend left=30] (0);
		\draw [style=arrow] (4) to[bend left=30] (3);
		\draw [style=arrow] (3) to[bend left=30] (4);
		\draw [style=arrow] (4) to[bend left=30] (5);
		\draw [style=arrow] (5) to[bend left=30] (3);
		\draw [style=arrow] (5) to[bend left=30] (6);
		\draw [style=arrow] (6) to[bend left=30] (5);
		\draw [style=arrow] (5) to[bend left=30] (7);
		\draw [style=arrow] (7) to[bend left=30] (3);
		\draw [style=arrow] (3) to[bend left=30] (7);
	\end{pgfonlayer}
\end{tikzpicture}
\]
This graph is decomposed into open graphs $A \maps \phi \to \{*\} $ and $B \maps \{*\} \to \phi$ as follows:
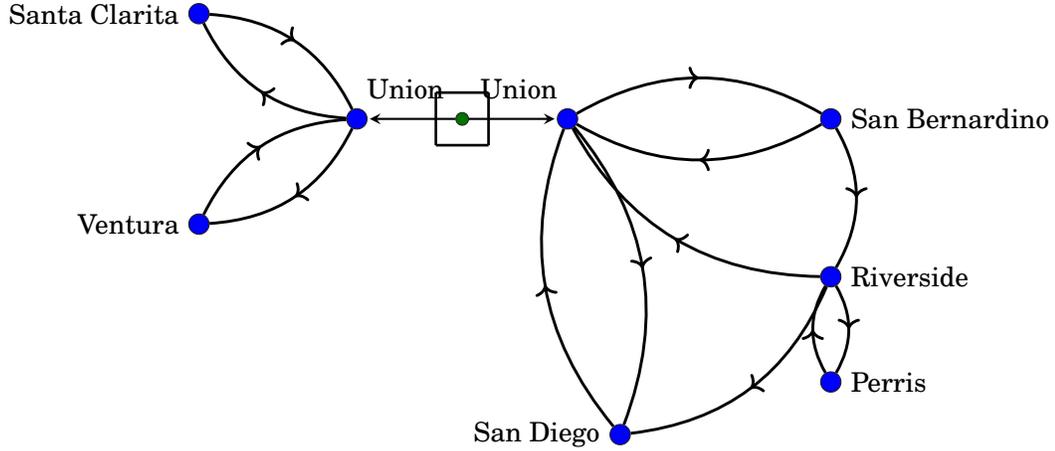
\[
\begin{tikzpicture}[scale=0.7]
	\begin{pgfonlayer}{nodelayer}
		\node [style=node,label=left:{Santa Clarita}] (0) at (-5, 2) {};
		\node [style=node,label=left:{Ventura}] (1) at (-5, -2) {};
		\node [style=node,label=87:{Union}] (2) at (-2, -0) {};
		\node [style=node,label=94:{Union}] (3) at (2, -0) {};
		\node [style=node,label=right:{San Bernardino}] (4) at (7, -0) {};
		\node [style=node,label=right:Riverside] (5) at (7, -3) {};
		\node [style=node,label=right:Perris] (6) at (7, -5) {};
		\node [style=node,label=left:{San Diego}] (7) at (3, -6) {};
		\node [style=inputdot] (8) at (0, -0) {};
		\node [style=none] (9) at (0.5, 0.5) {};
		\node [style=none] (10) at (-0.5, 0.5) {};
		\node [style=none] (11) at (-0.5, -0.5) {};
		\node [style=none] (12) at (0.5, -0.5) {};
	\end{pgfonlayer}
	\begin{pgfonlayer}{edgelayer}
		\draw [style=simple] (10.center) to (9.center);
		\draw [style=simple] (9.center) to (12.center);
		\draw [style=simple] (12.center) to (11.center);
		\draw [style=simple] (11.center) to (10.center);
		\draw [style=inarrow] (8) to (2);
		\draw [style=inarrow] (8) to (3);
		\draw [style=arrow] (0) to[bend left=30] (2);
		\draw [style=arrow] (1) to[bend left=30] (2);
		\draw [style=arrow] (2) to[bend left=30] (1);
		\draw [style=arrow] (2) to[bend left=30] (0);
		\draw [style=arrow] (4) to[bend left=30] (3);
		\draw [style=arrow] (3) to[bend left=30] (4);
		\draw [style=arrow] (4) to[bend left=30] (5);
		\draw [style=arrow] (5) to[bend left=30] (3);
		\draw [style=arrow] (5) to[bend left=30] (6);
		\draw [style=arrow] (6) to[bend left=30] (5);
		\draw [style=arrow] (5) to[bend left=30] (7);
		\draw [style=arrow] (7) to[bend left=30] (3);
		\draw [style=arrow] (3) to[bend left=30] (7);
	\end{pgfonlayer}
\end{tikzpicture}
\]
The category $F(A)$ has three objects, Santa Clarita, Union, and Ventura and the morphisms are given by all possible paths on the graph A. For example a morphism is given by the unique path
\[f=\text{Santa Clarita} \to \text{Union} \to \text{Ventura} \to \text{Union} \]
Similarly the category $F(B)$ has stations as objects and paths as morphisms. For example, $F(B)$ contains the unique morphism represented by the path
\[g= \text{Union} \to \text{San Bernardino} \to \text{Riverside} \to \text{Union} \]
Because $\Open(F)$ is a double functor, it is equipped with an isomorphism
\[ \begin{tikzcd} & F(A) +_{1} F(B) \ar[dd,"\sim"] & \\
0 \ar[ur]\ar[dr] & & \ar[ul] \ar[dl] 0 \\
 & F(G)
\end{tikzcd} \]
where 0 is the empty category and 1 is the terminal category. This isomorphism builds the operational semantics $F(G)$ using the operational semantics $F(A)$ and $F(B)$. First, paths such as $f$ and $g$ are formed then combined together using the pushout of categories. The pushout of categories requires a closure under all paths which travel back and forth between $A$ and $B$. In particular, because $f$ and $g$ both start and end at Union, $F(G)$ must contain all words in $f$ and $g$ as morphisms. 
\end{expl}

The additional axioms of a double functor ensure that building the operational semantics of an open graph from its components is independent of reassociating composition or adding identity horizontal morphisms. The second transitive closure required to join together operational semantics can be very computationally expensive and this leads to a combinatorial explosion. Although $\Open(F)$ provides a useful conceptual framework for building the operational semantics of graphs recursively, the isomorphisms provided are unlikely to provide a speed up without further assumptions or restrictions. In the next section we identify a subclass of open graphs, called functional, which can be joined together without combinatorial explosion.



\section{Black-boxing and Functional Open Graphs}\label{blackboxgraph}
In this section we define a ``black-boxing" functor, which restricts an open category to the data between its boundaries. This matches the traditional meaning of black-boxing in systems theory, i.e.\ forgetting about the internal workings of a system and concentrating only on the relationship it induces between its inputs and ouputs. In this case, the black-boxing of an open category is a profunctor whose data consists of the morphisms that travel from the input ports to the output port.

Our black-boxing is in general only laxly functorial. However, we define a subclass of open graphs for which the functoriality is strict.
\begin{defn}
For an open category
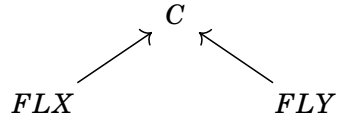
\[\begin{tikzcd}& C & \\
FLX \ar[ur] & & \ar[ul] FLY\end{tikzcd} \]
its black-boxing $\blacksquare(C)$ is a profunctor
\[\blacksquare(C) \maps FLX \times FLY \to \Set \]
given by $\blacksquare(C)(x,y) = C(i(x),j(y))$.
\end{defn}

\noindent This black-boxing operation extends to a lax double functor into the double category of profunctors \cite[Ex.\ 2.6]{framedbicat}.
\begin{defn}
Let $\mathsf{Prof}$ be the double category where
\begin{itemize}
    \item objects are categories,
    \item vertical morphisms are functors,
    \item horizontal morphisms are profunctors, and
    \item vertical 2-cells are squares
    \[
    \begin{tikzcd}
    A \ar[r,"P"] \ar[d,"f",swap]& B \ar[d,"g"] \\
    C \ar[r,"Q",swap] & D
    \end{tikzcd}
    \]
    equipped with a natural transformation
    \[
    \begin{tikzcd}
    A \times B^{op} \ar[dr,"P"{name=L},swap] \ar[rr,"f \times g^{op}"] & & C \times D^{op} \ar[dl,"Q"{name=R}]\ar[from=L,to=R,Rightarrow,"\alpha",shorten <= 5ex,shorten >= 5ex,yshift=1.5ex] \\
    & \Set &
    \end{tikzcd}
    \]
\end{itemize}
\end{defn}
\begin{thm}\label{blackbox}
Black-boxing lifts to a double functor
\[ \blacksquare \maps \Open(\Cat) \to \mathsf{Prof}\]
\end{thm}
\begin{proof}
On objects $\blacksquare$ sends a set $X$ to its discrete category $1_X$ and functions are sent to their unique extensions on these discrete categories. A $2$-cell
\[
\begin{tikzcd}
FLX \ar[d,"FLf",swap] \ar[r,"i"] & C\ar[d,"g"] & \ar[l,"j",swap] FLY \ar[d,"FLh"] \\
FLX' \ar[r,"i'",swap] & D & \ar[l,"j'"] FLY'
\end{tikzcd}
\]
is sent to the natural transformation
\[\begin{tikzcd}FLX \times FLY \ar[rr,"FLf \times FLh"] \ar[dr,"\blacksquare(C)"{name=L},swap] & & FLX' \times FLY' \ar[dl,"\blacksquare(D)"{name=R}]\ar[Rightarrow, from=L, to=R,"\blacksquare(g)",shorten <=5ex, shorten >=5ex,yshift=1.5ex] \\
& \Set & \end{tikzcd}\]
with components
\[ \blacksquare(g) (x,y) \maps \blacksquare(C)(x,y) \to  \blacksquare(D)(f(x),h(y))\]
given by the restriction of the functor $g$ to the hom-set $C(i(x),j(y))$. $\blacksquare(g)$ is well-defined because of the commutativity of the $2$-cell it comes from. Naturality of $\blacksquare(g)$ is trivial because $FLX \times FLY$ has no non-identity arrows. Functoriality of $\blacksquare$ on object and arrow categories follows immediately from the definitions. An identity $2$-cell in $\Open(\Cat)$
\[\begin{tikzcd}
& FLX & \\
FLX \ar[ur,equals] & & FLX \ar[ul,equals]
\end{tikzcd} \]
is sent to a profunctor $\blacksquare(FLX)$ given by
\[ \blacksquare(FLX)(x,x')= \delta_{xx'}\]
where $\delta_{xx'}$ is the function which returns a one element set when $x=x'$ and the empty set otherwise. $\blacksquare(FLX)$ is clearly the identity profunctor on the category $FLX$ so our double functor $\blacksquare$ preserves identities. Consider a composable pair of horizontal morphisms
\[\begin{tikzcd}
& D& & C & \\
FLX \ar[ur] & & \ar[ul] \ar[ur] FLY & & FLZ \ar[ul]
\end{tikzcd} \]
in $\Open(\Cat)$. Note that for every $y \in Y$ there is a function
\[\circ_y \maps \blacksquare(D)(x,y) \times \blacksquare(C)(y,z)  \to \blacksquare(C \circ D)(x,z)  \]
sending a pair of morphisms to their composite in $C \circ D$. The composition comparison of $\blacksquare$ is a natural transformation
\[\begin{tikzcd}FLX \times FLZ  \ar[dd,bend right=50,"\blacksquare(C) \circ \blacksquare(D)"{name=L},swap]  \ar[dd,bend left=50,"\blacksquare(C \circ D)"{name=R}]\ar[Rightarrow, from=L, to=R,"\alpha",shorten <=2ex, shorten >=2ex] \\
\\
\Set  \end{tikzcd} \]
with components
\[ \alpha_{x,z} \maps \int_{y \in Y} \blacksquare(D) (x,y) \times \blacksquare(C)(y,z) \to \blacksquare(C \circ D) (x,z)\]
given by stitching together the functions $\circ_y$ with the universal property of coends.
\end{proof}

\begin{expl}\label{loop}
Let $G\maps 1 \to 2$ be the open graph
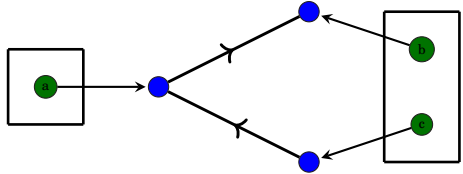
\[
\begin{tikzpicture}
	\begin{pgfonlayer}{nodelayer}
		\node [style=node] (0) at (-1, -0) {};
		\node [style=node] (1) at (1, 1) {};
		\node [style=node] (2) at (1, -1) {};
		\node [style=none] (3) at (-3, 0.5) {};
		\node [style=none] (4) at (-2, 0.5) {};
		\node [style=none] (5) at (-2, -0.5) {};
		\node [style=none] (6) at (-3, -0.5) {};
		\node [style=none] (7) at (2, 1) {};
		\node [style=none] (8) at (3, 1) {};
		\node [style=none] (9) at (2, -1) {};
		\node [style=none] (10) at (3, -1) {};
		\node [style=inputdot] (11) at (-2.5, -0) {a};
		\node [style=inputdot] (12) at (2.5, 0.5) {b};
		\node [style=inputdot] (13) at (2.5, -0.5) {c};
	\end{pgfonlayer}
	\begin{pgfonlayer}{edgelayer}
		\draw [style=arrow] (0) to (1);
		\draw [style=arrow] (2) to (0);
		\draw [style=simple] (3.center) to (6.center);
		\draw [style=simple] (3.center) to (4.center);
		\draw [style=simple] (4.center) to (5.center);
		\draw [style=simple] (6.center) to (5.center);
		\draw [style=simple] (7.center) to (9.center);
		\draw [style=simple] (9.center) to (10.center);
		\draw [style=simple] (10.center) to (8.center);
		\draw [style=simple] (8.center) to (7.center);
		\draw [style=inarrow] (11) to (0);
		\draw [style=inarrow] (12) to (1);
		\draw [style=inarrow] (13) to (2);
	\end{pgfonlayer}
\end{tikzpicture}
\]
and let $H \maps 2 \to 1$ be the open graph
  
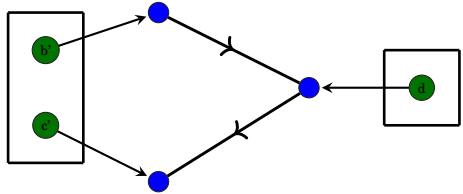
\[
  \begin{tikzpicture}

		\node [style=node] (0) at (-1, 1) {};
		\node [style=node] (1) at (-1, -1.25) {};
		\node [style=node] (2) at (1, -0) {};
		\node [style=none] (3) at (-2, 1) {};
		\node [style=none] (4) at (-2, -1) {};
		\node [style=none] (5) at (-3, -1) {};
		\node [style=none] (6) at (-3, 1) {};
		\node [style=inputdot] (7) at (-2.5, 0.5) {b'};
		\node [style=inputdot] (8) at (-2.5, -0.5) {c'};
		\node [style=none] (9) at (2, 0.5) {};
		\node [style=none] (10) at (2, -0.5) {};
		\node [style=none] (11) at (3, -0.5) {};
		\node [style=none] (12) at (3, 0.5) {};
		\node [style=inputdot] (13) at (2.5, -0) {d};

		\draw [style=arrow] (0) to (2);
		\draw [style=arrow] (2) to (1);
		\draw [style=simple] (6.center) to (5.center);
		\draw [style=simple] (5.center) to (4.center);
		\draw [style=simple] (4.center) to (3.center);
		\draw [style=simple] (3.center) to (6.center);
		\draw [style=inarrow] (7) to (0);
		\draw [style=inarrow] (8) to (1);
		\draw [style=simple] (9.center) to (10.center);
		\draw [style=simple] (9.center) to (12.center);
		\draw [style=simple] (12.center) to (11.center);
		\draw [style=simple] (10.center) to (11.center);
		\draw [style=inarrow] (13) to (2);
\end{tikzpicture}
\]
so that their composite $H \circ G \maps 1 \to 1$ is
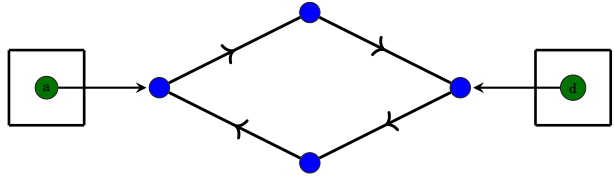
\[
\begin{tikzpicture}
	\begin{pgfonlayer}{nodelayer}
		\node [style=node] (0) at (-1, -0) {};
		\node [style=node] (1) at (1, 1) {};
		\node [style=node] (2) at (1, -1) {};
		\node [style=none] (3) at (-3, 0.5) {};
		\node [style=none] (4) at (-2, 0.5) {};
		\node [style=none] (5) at (-2, -0.5) {};
		\node [style=none] (6) at (-3, -0.5) {};
		\node [style=inputdot] (7) at (-2.5, -0) {a};
		\node [style=node] (8) at (3, -0) {};
		\node [style=none] (9) at (4, 0.5) {};
		\node [style=none] (10) at (4, -0.5) {};
		\node [style=none] (11) at (5, 0.5) {};
		\node [style=none] (12) at (5, -0.5) {};
		\node [style=inputdot] (13) at (4.5, -0) {d};
	\end{pgfonlayer}
	\begin{pgfonlayer}{edgelayer}
		\draw [style=arrow] (0) to (1);
		\draw [style=arrow] (2) to (0);
		\draw [style=simple] (3.center) to (6.center);
		\draw [style=simple] (3.center) to (4.center);
		\draw [style=simple] (4.center) to (5.center);
		\draw [style=simple] (6.center) to (5.center);
		\draw [style=inarrow] (7) to (0);
		\draw [style=arrow] (8) to (2);
		\draw [style=arrow] (1) to (8);
		\draw [style=simple] (9.center) to (10.center);
		\draw [style=simple] (9.center) to (11.center);
		\draw [style=simple] (11.center) to (12.center);
		\draw [style=simple] (12.center) to (10.center);
		\draw [style=inarrow] (13) to (8);
	\end{pgfonlayer}
\end{tikzpicture}
\]
The black-boxing $\blacksquare(G) \maps FL(1) \times FL(2) \to \Set$ has the singleton as $\blacksquare (G) (a,b)$ representing the only path from $a$ to $b$ and $\blacksquare (G)(a,c)$ is the empty set. Similarly, $\blacksquare(H) \maps FL(2) \times FL(1) \to \Set$ has the singleton for $\blacksquare(H)(b',d)$ and the empty set for $\blacksquare(c',d)$. Their composite $\blacksquare (H) \circ \blacksquare(G)\maps FL(1) \times FL(1) \to \Set$ has only one value and it's given by the coend formula
\[\blacksquare(H) \circ \blacksquare(G) (a,d) = \int_{x\in 2} \blacksquare(G)(a,x) \times \blacksquare(H)(x,d) \]
In this case the coend is again the singleton, whose unique element represents the path in $H \circ G$ traversing the top two edges. Note that there are many more elements in $\blacksquare(H \circ G)(a,d)$. Paths in $H \circ G$ may circle around any natural number of times before arriving at their destination and this gives an element of $\blacksquare(H \circ G)$ for every finite natural number $n$. The laxator of composition for $\blacksquare$ sends the unique element of $\blacksquare(H) \circ \blacksquare(G)(a,d)$ to element $0$ in $\blacksquare(H \circ D)(a,d)$ representing the path which goes directly from $a$ to $d$ without appending any loops.
\end{expl}
\noindent The previous example may be somewhat discouraging, in general $\blacksquare(H \circ G)$ will have much larger sets than $\blacksquare(H) \circ \blacksquare(G)$ so any technique for constructing the former from the latter will be rife with difficulty. However, we can define a sort of open graph for which the looping behavior of Example \ref{loop} is disallowed.
\begin{defn}\label{functgraph}
Let $G=\begin{tikzcd} E \ar[r,shift left=.5ex,"s"] \ar[r,shift right=.5ex,"t",swap] & V \end{tikzcd}$ be a graph. Then a vertex $x \in V$ is a \define{source} if it has no incoming edges i.e.\ $t^{-1}(x)=\phi$. Similarly, $x$ is a \define{sink} if it has no outgoing edge i.e.\ $s^{-1}(x) =\phi$.
An open graph
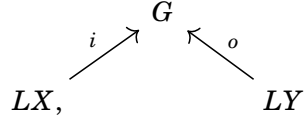
\[ \begin{tikzcd}
 & G & \\
 LX \ar[ur,"i"], & & \ar[ul,"o",swap] LY
\end{tikzcd} \]
is \define{functional} if for every $x \in X$, $i(x)$ is a source and for $y \in Y$, $j(y)$ is a sink.
\end{defn}\noindent The following theorem relies on a lemma.
\begin{lem}\label{binomialroot}
For functional open graphs $G \maps X \to Y$ and $H \maps Y \to Z$,
\[\blacksquare((H \circ G)^n) \cong \sum_{i+j=n} \blacksquare(H^i) \circ \blacksquare(G^j) \]
where the powers refer to iterated pullbacks of a graph with itself.
\end{lem}
\begin{proof} For elements $x \in X$ and $z \in Z$, let $f$ be an element of $\blacksquare((H \circ G)^n)(x,z)$. Then $f$ is a sequence of composable edges $(e_1,e_2,\ldots,e_n)$ such that the source of $e_1$ is $x$ and the target of $e_n$ is $z$. Because $G$ and $H$ are functional, if $e_{i}$ is an edge of $H$, then $e_{i+1}$ cannot be an element of $G$. Therefore there is a last occurence $1 \leq i_* \leq n$ such that $e_{i_*}$ is an edrnge of $G$ and $e_i$ is an edge of $H$ for all $i>i_*$. The
$n$-tuple $(e_1,e_2,\ldots,e_n)$ can be split into tuples $(e_1,e_2,\ldots, e_{i_*}) \in \blacksquare (G^{i_*})(x,y)$ and $(e_{i_*+1},e_{i_* +2},\ldots,e_n) \in \blacksquare(H^{n-i_*})(y,z)$. The composite profunctor
\[\blacksquare(H^{n-i_*}) \circ \blacksquare(G^{i_*}) (x,z) = \int_{y \in Y} \blacksquare(G^{i_*})(x,y) \times \blacksquare(H^{n-i_*})(y,z)\]
accounts for every possible $y$ value that the path can stop in. Because $i_*$ can occur at any value, we need to to take the coproduct of the above composite for every power of $H$ and $G$ summing to $n$ in order to account for every element of $\blacksquare((H \circ G)^n)(x,z)$.
\end{proof}
\begin{thm}\label{functionalgraph}
The composite
\[\blacksquare \circ \Open(F) \maps \Open(\Grph) \to \mathsf{Prof}\]
preserves horizontal composition on functional open graphs up to isomorphism.
\end{thm}
\begin{proof}
We show that for functional open graphs $G \maps X \to Y$ and $H \maps Y \to Z$, there is an isomorphism
\[\blacksquare(\sum_{n \geq 0} (H \circ G)^n) \cong \blacksquare(\sum_{i \geq 0} H^i) \circ \blacksquare(\sum_{j \geq 0} G^j)  \]
Starting with the right hand side, we pull the $\blacksquare$'s inside the sums and use the distributive law
\[ \blacksquare(\sum_{i \geq 0} H^i) \circ \blacksquare(\sum_{j \geq 0} G^j )\cong \sum_{i \geq 0} \sum_{j \geq 0} \blacksquare(H^{i}) \circ \blacksquare(G^{j}) \]
The values of $i$ and $j$ in this sum are represented by the grid points on the cartesian axes below
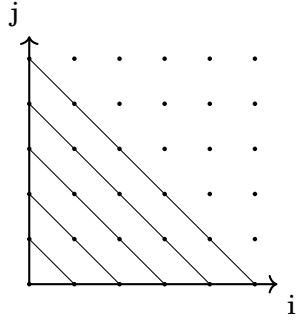
\[
\begin{tikzpicture}[scale=.6]
    \foreach \x in {0,1,...,5} {
        \foreach \y in {0,1,...,5} {
            \fill[color=black] (\x,\y) circle (0.05);
        }
    }
\draw (0,1) -- (1,0);
\draw (0,2) -- (2,0);
\draw (0,3) -- (3,0);
\draw (0,4) -- (4,0);
\draw (0,5) -- (5,0);
\draw[thick,->] (0,0) -- (5.5,0) node[anchor=north west] {i};
\draw[thick,->] (0,0) -- (0,5.5) node[anchor=south east] {j};
\end{tikzpicture}
\]
\noindent The above double coproduct is arranged so that entries are summed horizontally first and then going up vertically. A rearrangement this sum is represented by the diagonal lines, starting from $(0,0)$ and summing each diagonal before moving on to the next. The $n$-th diagonal line crosses all coordinates $(i,j)$ with $i+j=n$ so the above sum can be rearranged as the sum going over the diagonals i.e.
\[\sum_{n \geq 0} \sum_{i+j=n} \blacksquare(H^{i}) \circ \blacksquare(G^{j}).   \]
By Lemma \ref{binomialroot} this composite is isomorphic to
\[\sum_{n \geq 0} \blacksquare((H \circ G)^n) \]
Pulling $\blacksquare$ out of the sum gives the desired result.
\end{proof}
This theorem says the operational semantics of a composite of open graphs can be computed using profunctor composition when paths can only go from the first open graph to the second. Lemma \ref{binomialroot}, says that the paths which occur in exactly $n$ steps on a composite of functional open graphs can be computed compositionally using the given formula. In what follows we will generalize these results to a larger range of networks which can be described by enriched and internal graphs.

\chapter{Operational Semantics of Q-Nets}\label{QNetchap}
In this chapter we construct the operational semantics for Petri nets and many of their related variants by generalizing the free category construction of Proposition \ref{freecat} to the case of graphs internal to $\Mod(\Q)$, the category of models of a Lawvere theory $\Q$. Petri nets are not internal graphs. However, they are the generating data for graphs internal to $\CMon$, the category of commutative monoids. In Section \ref{defin} we review the basic definitions of Petri net theory. If Petri nets are the generating data for graphs internal to $\CMon$, more generally, what is the generating data for graphs internal to $\Mod(\Q)$? This chapter answers this question by introducing $\Q$-nets, a generalization of Petri net based on the operations and axioms of a Lawvere theory $\Q$. In Section \ref{QNet}, we define $\Q$-nets and show how many well-known variants of Petri nets and the relationships between them can be understood using this definition. The main goal of this chapter is to construct an adjunction representing the operational semantics of $\Q$-nets. We start with a motivating example: when $\Q$ is the Lawvere theory for commutative monoids we obtain Petri nets. In Section \ref{CMC}, we construct an adjunction
\[\begin{tikzcd}\Petri \ar[r,bend left,"F"] & \ar[l,bend left,"U"] \CMC \end{tikzcd}\]
where $\CMC$ is the category of commutative monoidal categories. In Section \ref{gen} we state Theorem \ref{big}, the main result of this chapter: for every Lawvere theory $\Q$, there is an adjunction
\[\begin{tikzcd}
\Net{\Q} \ar[r,bend left,"F_\Q"] & \ar[l,"U_\Q",bend left] \cat{\Q}
\end{tikzcd}
\]
between the category of $\Q$-nets and the category of ``$\Q$-categories", i.e., models of $\Q$ internal to the category of categories. For a $\Q$-net $P$, the $\Q$-category $F_\Q(P)$ is a category whose morphisms represent the execution sequences of $P$. To construct this adjunction we factor it as
\[
\begin{tikzcd}
\Net{\Q} \ar[r,bend left,"\A_\Q"] \ar[r,phantom,"\bot",pos=.6] & \ar[l,bend left, "{\backa}_\Q"] \Grph(\Mod(\Q)) \ar[r,bend left, "\B_\Q"]\ar[r,phantom,"\bot"] & \ar[l,bend left,"\backb_\Q"] \Mod(\Q,\Cat)
\end{tikzcd}
\]
where $\Grph(\Mod(\Q))$ is the category of graphs internal to $\Mod(\Q)$. In Section \ref{freeQgraph} we construct the first part $\A \dashv \backa$ and in Section \ref{freecatinternal} we construct the second part $\B \dashv \backb$.

\section{Petri Nets and Their Executions}\label{defin}
Petri nets are the motivating example for the generalizations considered here. Therefore we review their properties here in depth to provide intuition for the more abstract treatment later on. Petri nets are regarded as a graph, whose source and target land in a free commutative monoid. To describe this graph we need the following adjunction. 
\begin{defn}\label{N}
Let $L \maps \Set \to \CMon$ be the free commutative monoid functor, that is, the left adjoint of the functor $R \maps \CMon \to \Set$ that sends commutative monoids to their underlying sets and monoid homomorphisms to their underlying functions. Let 
\[ \N \maps \Set \to \Set \] be the \define{free commutative monoid monad} given by the composite $R\circ L$.
\end{defn}
\noindent For any set $X$, $\N[X]$ is the set of formal finite linear combinations of elements of $X$ with natural number coefficients.   The unit of $\N$ is given by the natural inclusion of $X$ into $\N[X]$, and for any function $f \maps X \to Y$, $\N[f] \maps \N[X] \to \N[Y]$ is the unique monoid homomorphism that extends $f$.   

\begin{defn}\label{PetriNet}
We define a \define{Petri net} to be a pair of functions of the following form:
\[\xymatrix{ T \ar@<-.5ex>[r]_-t \ar@<.5ex>[r]^-s & \N[S]. } \]
We call $T$ the set of \define{transitions}, $S$ the set of \define{places}, $s$ the \define{source} function, and $t$ the \define{target} function. 
	
\end{defn}

\begin{defn}\label{PetriMorphism}
A \define{Petri net morphism} from the Petri net 
$\xymatrix{ T \ar@<-.5ex>[r]_-t \ar@<.5ex>[r]^-s & \N[S] }$ to
the Petri net $\xymatrix{ T' \ar@<-.5ex>[r]_-{t'} \ar@<.5ex>[r]^-{s'} & \N[S']}$ 
is a pair of functions $(f \maps T \to T', g \maps S \to S')$ such that the diagrams 

	\[
	\xymatrix{ 
		T \ar[d]_f  \ar[r]^-{s} & \N[S] \ar[d]^-{\N[g]} \\	
		T' \ar[r]_-{s'} & \N[S'] 
	}
	\qquad
	\xymatrix{ 
		T \ar[d]_f  \ar[r]^-{t} & \N[S] \ar[d]^-{\N[g]} \\	
		T' \ar[r]_-{t'} & \N[S'] . 
	}
	\]
	commute.
\end{defn}

\begin{defn}
Let $\Petri$ be the category of Petri nets and Petri net morphisms, with composition 
defined by 
\[  (f, g) \circ (f',g') = (f \circ f' , g \circ g')  .\]
\end{defn}
\noindent Our definition of Petri net morphism differs from the earlier definition used by Degano--Meseguer--Montanari \cite{DMM} and Sassone \cite{SassoneStrong,SassoneAxiom,functorialsemantics}. The difference is that our definition requires that the homomorphism between free commutative monoids come from a function between the sets of places whereas the above references allow arbitrary commutative monoid homomorphisms. This difference of definition is present in our definition of the category of $\Q$-nets as well. With this change, the categories $\Petri$ and $\Net{\Q}$ become complete and cocomplete as shown in Proposition \ref{qnetcocomplete}. This is important for the compositionality results in Chapter \ref{openQNets}.

Petri nets have a natural semantics which is described by "the token game''. This is a game where each place of a Petri net is equipped with a natural number of tokens. Players are then allowed to shuffle the tokens from place to place using the transitions. The token game is formalized by the notions of marking and firing. 
\begin{defn}
A \define{marking} of a Petri net $P = \begin{tikzcd}T \ar[r,shift left=.5ex,"s"] \ar[r,shift right=.5ex,"t",swap] & \N[S] \end{tikzcd}$ is an element $m \in \N[S]$, or equivalently, a function $m \maps S \to \N$ which is zero on all but a finite number of elements. A \define{firing} of $P$ is a tuple $(\tau,x,y)$, where $\tau$ is a transition and $x$ and $y$ are markings of $P$ with $x - s(\tau) \geq 0$ and $x-s(\tau) + t(\tau) = y$.
\end{defn}
Firings can be chained together in sequence: for a firing $(\tau,x,y)$ and a firing $(\sigma,y,z)$ we can define their composite as a tuple $(\sigma \circ \tau,x,z)$ where $\circ$ is a formal symbol. Firings can also be performed in parallel: for two firings $(\tau,x,y)$ and $(\tau',x',y')$ there is a parallelization $(\tau +\tau', x+x',y+y')$. This suggests that firings of a Petri net have the structure of a monoidal category. Meseguer and Montanari were the first to notice this and show how Petri nets can be turned into commutative monoidal categories \cite{monoids}.
\begin{defn}\label{CMCdef}
	A \define{commutative monoidal category} is a commutative monoid object internal to $\Cat$. Explicitly, a commutative monoidal category is a strict monoidal category $(C,\otimes,I)$, such that, for all objects $a$, $b$ and morphisms $f$, $g$ in $C$ 
	\[a \otimes b = b \otimes a \text{ and } f \otimes g = g \otimes f.\]
\end{defn}
 
\noindent Note that a commutative monoidal category is the same as a strict symmetric monoidal category where the symmetry isomorphisms $\sigma_{a,b} \maps a \otimes b \stackrel{\sim}{\longrightarrow} b \otimes a$ are all identity morphisms. In fact, a commutative monoidal category is precisely a category where the objects and morphisms form commutative monoids and the structure maps are commutative monoid homomorphisms. A commutative monoidal category where the morphisms represent sequences of firings of a Petri net $P$ will be referred to as the "operational semantics" of $P$. In this chapter we characterize this semantics construction as an adjunction between the category of Petri nets and the following category.
\begin{defn}
	Let $\CMC$ be the category whose objects are commutative monoidal categories and whose morphisms are strict monoidal functors.
\end{defn}
Note that every monoidal functor between commutative monoidal categories is automatically a strict symmetric monoidal functor, so the adjective symmetric is not included in the above definition.

\section{Q-Nets}\label{QNet}
Petri nets need not have a free commutative monoid of places, and this aspect can be generalized using Lawvere theories. A review of the basic definitions and properties of Lawvere theories can be found in in Appendix \ref{law}. As in Definition \ref{models}, let $\Mod(\law{Q})$ be the category of models of $\law{Q}$ in $\Set$,
\[ \begin{tikzcd} \Set \ar[r,bend left,"L_\Q",pos=.6] \ar[r,phantom,"\bot",pos=.6] & \ar[l,bend left,"R_\Q",pos=.4] \Mod(\law{Q}) \end{tikzcd}\]
be the adjunction generating free models of $\law{Q}$ and let $M_\law{Q}$ be the composite $R_\law{Q} \circ L_\law{Q}$.
\begin{defn}\label{QNetd}	
	Let $\Net{\law{Q}}$ be the category where
\begin{itemize}
	    \item objects are $\law{Q}$\define{-nets}, i.e.\ pairs of functions of the form
	    \[ \begin{tikzcd}
	    T \ar[r, shift left=.5ex, "s"] \ar[r, shift right=.5ex,"t",swap] & M_{\law{Q}} S
	    \end{tikzcd} \]
		\item a morphism from the $\law{Q}$-net 
$\xymatrix{ T \ar@<-.5ex>[r]_-t \ar@<.5ex>[r]^-s & M_\law{Q} S }$ to
the $\law{Q}$-net $\xymatrix{ T' \ar@<-.5ex>[r]_-{t'} \ar@<.5ex>[r]^-{s'} & M_\law{Q} S'}$ 
is a pair of functions $(f \maps T \to T', g \maps S \to S')$ such that the following diagrams commute:
	\[
	\xymatrix{ 
		T \ar[d]_f  \ar[r]^-{s} & M_\law{Q} S \ar[d]^-{M_\law{Q} g} \\	
		T' \ar[r]_-{s'} & M_\law{Q} S' 
	}
	\qquad
	\xymatrix{ 
		T \ar[d]_f  \ar[r]^-{t} & M_\law{Q}S \ar[d]^-{M_\law{Q} g} \\	
		T' \ar[r]_-{t'} & M_\law{Q} S' . 
	}
	\]
	\end{itemize}
\end{defn}
	This definition lifts to a functor. Let $M_\law{Q}$ be as before and let $M_{\law{R}} \maps \Set \to \Set$ be corresponding monad induced by a Lawvere theory $\law{R}$. Every morphism of Lawvere theories $f \maps \law{Q} \to \law{R}$ induces a functor
	\[f_* \maps \Mod(\law{R}) \to  \Mod(\law{Q})\]
	which composes every model of $\law{R}$ with $f$. A left adjoint
	\[f^* \maps \Mod(\law{Q}) \to \Mod(\law{R}) \] is given by the left Kan extension of each model along $f$ \cite{ttt,buckley}. Now, we have the following commutative diagram of functors
	
	\[
	\xymatrix{\Mod(\law{Q}) \ar[dr]_{R_\law{Q}}  & &\ar[ll]_{f_*} \Mod(\law{R}) \ar[dl]^{R_{\law{R}}} \\
	& 
	 \Set  & }
	\]
	all of which have left adjoints. Given this set of assumptions, there is a morphism of monads $M^f$ given by
	\[M^f = R_\law{Q} \eta L_\law{Q} \maps M_\law{Q} \Rightarrow  M_{\law{R}} \]
	where $\eta$ is the unit of the adjunction $f^* \dashv f_*$. This can either be verified directly, or by using the adjoint triangle theorem \cite{triangles}. In what follows we will use this morphism of monads to translate between different types of $\law{Q}$-nets.

\begin{defn}\label{netf}
	    Let 
	    \[ \Net{(-)}  \maps \Law \to \Cat \]
	    be the functor which sends a Lawvere theory $\law{Q}$ to the category $\Net{\law{Q}}$ and sends a morphism $f\maps \law{Q} \to \law{R}$ of Lawvere theories to the functor  $\Net{f} \maps \Net{\law{Q}} \to \Net{\law{R}}$ which sends a $\law{Q}$-net 
            \[
            \begin{tikzcd}
            T \ar[r, shift left=.5ex, "s"] \ar[r,shift right=.5ex,"t",swap] & M_{\law{Q}} S
            \end{tikzcd}
            \]
to the $\law{R}$-net
            \[ 
            \begin{tikzcd}
            T \ar[r, shift left=.5ex, "s"] \ar[r,shift right=.5ex,"t",swap] & M_{\law{Q}} S \ar[r,"M_S^{f}"] & M_{\law{R}} S
            \end{tikzcd} 
            \] 

\noindent For a morphism of $\law{Q}$-nets $(g\maps T \to T',h \maps S \to S')$, $\Net{f} (g,h)$ is $(g,h)$. This is well-defined because of the naturality of $M^f$.
\end{defn}
\noindent Varying the Lawvere theory $\Q$ gives many known variants of Petri nets.
\begin{expl}
Setting $\Q$ equal to $\law{CMON}$, the Lawvere theory for commutative monoids, we obtain the category of Petri nets.
\end{expl}

\begin{defn}\label{prenetdef}
Let $(-)^* \maps \Set \to \Set$ denote the monad that the Lawvere theory $\law{MON}$ induces via the correspondence in \cite{linton}. For a set $X$, $X^*$ is given by the underlying set of the free monoid on $X$.
 A \define{pre-net} is a pair of functions of the form
 \[
 \begin{tikzcd}
 T \ar[r, shift left=.5ex, "s"] \ar[r, shift right=.5ex, "t",swap] & S^*
 \end{tikzcd}
 \]
 A morphism of pre-nets from a pre-net $(s,t \maps T \to S^*)$ to a pre-net $(s',t'\maps T' \to S'^*)$ is a pair of functions $(f \maps T \to T', g\maps S \to S')$ which preserves the source and target as in Definition \ref{PetriMorphism}. This defines a category $\PreNet$.
\end{defn}
\begin{expl} \label{prenet}
If we take $\law{Q}= \law{MON}$, the Lawvere theory of monoids, we get the category $\PreNet$.
\end{expl} \noindent A description of $\law{MON}$ can be found in the Appendix. $\PreNet$ has the same objects as the category introduced in \cite{functorialsemantics} but the morphisms are restricted as in Definition \ref{PetriMorphism}. Pre-nets are the same as tensor schemes introduced by Joyal and Street in \cite{scheme}. The authors define a notion of free category on a tensor scheme and Bruni, Meseguer, Montanari, and Sassone construct an adjunction between pre-nets and a subcategory of the category of strict symmetric monoidal categories $\SSMC$ \cite{functorialsemantics}.

In \cite{master_2020}, a closely related adjunction 
	
	\[\begin{tikzcd}
	\PreNet \ar[r,bend left,"Z"] \ar[r,phantom,"\bot"] & \ar[l,bend left, "K"] \SSMC
	\end{tikzcd} \]is constructed which does not require the restriction to a subcategory of $\SSMC$. Pre-nets are useful because after forming an appropriate quotient, the category $Z(P)$ for a pre-net $P$ is equivalent the category of \textit{strongly concatenable processes} which can be performed on the net. This equivalence is important for realizing the individual token philosophy \cite{functorialsemantics}. The individual token philosophy, as opposed to the collective token philosophy, gives identities to the individual tokens and keeps track of the causality in the executions of a Petri net.
\begin{expl}
In 2013 Bartoletti, Cimoli, and Pinna introduced lending Petri nets \cite{lending}. These are Petri nets where arcs can have a negative multiplicity and tokens can be borrowed in order to fire a transition. Lending nets are also equipped with a partial labeling of the places and transitions so they can be composed and are required to have no transitions which can be fired spontaneously. In 2018 Genovese and Herold introduced integer nets \cite{genovese}. Let $\law{ABGRP}$ be the Lawvere theory of abelian groups. This Lawvere theory contains three generating operations
\[e \maps 0 \to 1 \text{,  } i \maps 1 \to 1 \text{, and } m \maps 2 \to 1\]
representing the identies, inverses, and multiplication of an abelian group. These generating morphisms are required to satisfy the axioms of an abelian group; associativity, commutativity, and the existence of inverses and an identity.
The category of integer nets, modulo a change in the definition of morphisms, can be obtained by taking $\law{Q} = \law{ABGRP}$ in the definition of $\Net{\law{Q}}$.
\begin{defn}\label{integernets}
    Let $\Z \maps \Set \to \Set$ be the free abelian group monad which for a set $X$ generates the free abelian group $\Z [X]$ on the set $X$. Note that $\Z$ is the monad induced by the Lawvere theory $\law{ABGRP}$ via the correspondence in \cite{linton}. An \define{integer net} is a pair of functions of the form
    \[\xymatrix{ T \ar@<-.5ex>[r]_-t \ar@<.5ex>[r]^-s & \Z[S]. } \]
    A \define{morphism of integer nets} is a pair $(f \maps T \to T', g \maps S\to S')$ which makes the diagrams analogous to the definition of Petri net morphism (Definition \ref{PetriMorphism}) commute. 
    Let $\Net{\Z}$ be the category where objects are integer nets and morphism are morphisms of integer nets.
\end{defn}
\noindent Integer nets are useful for modeling the concepts of credit and borrowing. There is a correspondence between lending Petri nets and propositional contract logic; a form of logic useful for ensuring that complex networks of contracts are honored \cite{lending}. Genovese and Herold constructed a categorical semantics for integer nets \cite{genovese}. In \cite{master_2020}, the author constructed a variation of this semantics which uses the general framework developed in this chapter.
\end{expl}

\begin{expl}
Elementary net systems, introduced by Rozenberg and Thiagarajan in 1986, are are Petri nets with a maximum of one edge between a given place and transition \cite{elementary}.

\begin{defn}
 An \define{elementary net system} is a pair of functions 
 \[
 \begin{tikzcd}
 T  \ar[r, shift left=.5ex, "s"] \ar[r, shift right=.5ex, "t", swap] & 2^S
 \end{tikzcd}
 \] 
 where $2^S$ denotes the power set of $S$.
\end{defn}

\noindent Elementary net systems can be obtained from our general formalism. Let $\law{SLAT}$ be the Lawvere theory for semi-lattices, i.e.\ commutative idempotent monoids. This Lawvere theory contains morphisms
\[ m \maps 2 \to 1 \text{ and } e \maps 0 \to 1\]
as in $\law{MON}$ the theory of monoids. Also similar to $\law{MON}$, $\law{SLAT}$ is quotiented by the associativity and unitality axioms given in Example \ref{mon}. In addition, $\law{SLAT}$ has the following axioms representing commutativity and idempotence
\[
\begin{tikzcd}
2 \ar[rr, "\sigma"] \ar[dr, "m",swap] & & 2 \ar[dl, "m"] \\
& 1 &
\end{tikzcd}
\quad \quad 
\begin{tikzcd}
1\ar[rr, "\Delta"] \ar[dr, "id",swap] & & 2 \ar[dl, "m"] \\
& 1 &
\end{tikzcd}
\]
where $\sigma \maps 2 \to 2$ is the braiding of the cartesian product and $\Delta \maps 1 \to 2$ is the diagonal. For models in $\Set$, The first diagram says that you can multiply two elements in either order and the get the same thing. The second diagram says that if you multiply an element by itself you get itself. As in Definition \ref{models}, $\law{SLAT}$ corresponds to a monad on $\Set$. It is well known that this monad is the covariant power set monad
\[
2^{(-)} \maps \Set \to \Set
\]\noindent which sends a set $X$ to its set of finite subsets and a function to the mapping which sends subsets of $X$ to their image. This motivates the following:
\begin{defn}
 Let $\Net{\law{SLAT}}$ be the category of elementary net systems obtained as in Definition \ref{QNet} for $\law{Q}= \law{SLAT}$.
\end{defn} 
\end{expl}
\begin{expl}\label{ksafe}
Generalizing the previous example, for any natural number $k$ we form the Lawvere theory $\law{k}$ for k-idempotent monoids. $\law{k}$ has the same operations and axioms as $\law{SLAT}$ except the idempotency axiom is replaced with the axiom
\[
\begin{tikzcd}
1\ar[rr, "\Delta^k"] \ar[dr, "id",swap] & & k \ar[dl, "m^k"] \\
& 1 &
\end{tikzcd}
\]
where $\Delta^k$ is the $k$-fold diagonal map and $m^k$ is the $k$-fold multiplication map. Via the correspondence of Linton, $\law{k}$ gives the $k$-powerset monad
\[k^{(-)} \maps \Set \to \Set \]
sending a set $X$ to the set of finitely supported functions $\{X \to \{0,1,2,\ldots,k-1\}=k\}$ \cite{linton}. For a function $f \maps X \to Y$, $k^f \maps k^X \to k^Y$ is defined by
\[k^f(X \xrightarrow{a} k)(y) = \sum_{x \in f^{-1}(y)} a(x) \]
i.e.\ the analog of direct image for $k$-multisets.
\begin{defn}
A \define{$k$-safe} net is a pair of functions
\[\begin{tikzcd}T \ar[r,shift left=.5ex, "s"] \ar[r,shift right=.5ex,"t",swap] & k^S\end{tikzcd} \]
\noindent A morphism of $k$-safe nets is a pair of functions
\[
\begin{tikzcd}T\ar[d,"f",swap] \ar[r,shift left=.5ex, "s"] \ar[r,shift right=.5ex,"t",swap] & k^S \ar[d,"g"] \\
T' \ar[r,shift left=.5ex, "s'"] \ar[r,shift right=.5ex,"t'",swap] & k^{S'}\end{tikzcd} 
\]
commuting as in the previous definitions. This defines a category $\Net{\law{k}}$ of $k$-safe nets and their morphisms. 
\end{defn}
\end{expl}
The functoriality of Definition \ref{QNet} can be exploited to generate functors between different categories of $\law{Q}$-nets. 
There is the diagram in $\mathsf{Law}$
\begin{center}
\begin{tikzcd}
\law{SLAT}  &            &\\
\law{CMON} \ar[r,"b"] \ar[u,"a"] & \law{ABGRP} &\\
\law{MON} \ar[u,"c"] \ar[r,"d",swap] & \law{GRP} \ar[u,"e",swap]
\end{tikzcd}
\end{center}
where all the morphisms send their generating operations to their counterparts in the Lawvere theory of their codomain. These target Lawvere theories either have extra axioms or operations making the above functors not necessarily full or faithful:
\begin{itemize}
    \item $a$ sends the morphism $m \circ \Delta$ in $\law{CMON}$ to $id_1 \maps 1 \to 1$ in $\law{SLAT}$ to impose the idempotent law. All other other generating components are sent to their natural counterparts.
    \item $b$ and $d$ send every object and morphism to its natural analog. However, $\law{ABGRP}$ and $\law{GRP}$ have an extra operation $i \maps 1 \to 1$ representing inverses. This makes the functors $b$ and $d$ faithful but not full.
    \item $c$ and $e$ add the commutativity law; they send both the multiplication $m \maps 2 \to 1$ and the composite $\sigma \circ m\maps 2 \to 2$ of the braiding $\sigma \maps 2 \to 2$ to the multiplication map $m \maps 2 \to 2$ in the target Lawvere theory. This makes $c$ and $e$ not faithful.
\end{itemize}

\noindent Definition \ref{netf} is used to give a network between different flavors of $\law{Q}$-nets. By applying $\Net{(-)}$ to the above diagram we get the diagram of categories

\begin{center}
\begin{tikzcd}
\Net{\law{SLAT} } &            &\\
\Petri \ar[r,"\Net{b}"] \ar[u,"\Net{a}"] & \Net{\Z} &\\
\PreNet \ar[u,"\Net{c}"] \ar[r,"\Net{d}",swap] & \Net{\law{GRP}} \ar[u,"\Net{e}",swap]
\end{tikzcd}
\end{center}
The functors in this diagram are described as follows:
\[\Net{c} \maps \PreNet \to \Petri \]
 is often called ``abelianization" because it sends a pre-net to the Petri net which forgets about the ordering on the input and output of each transition. The authors of \cite{functorialsemantics} use $\Net{c}$ to explore the relationship between pre-nets and Petri nets. The functor $\Net{e}\maps \Net{\law{GRP}} \to \Net{\Z}$ gives the analogous relationship for integer nets.
\[\Net{b} \maps \Petri \to \Net{\Z}\]
is the functor which does not change the source and target of a given place. The only difference is that the markings of a $\Z$-net coming from a Petri net are thought of as elements of a free abelian group rather than a free abelian monoid. $\Net{d}$ is the analogous functor for pre-nets. 
\[ \Net{a} \maps \Petri \to \Net{\law{SLAT}} \]
is the functor which sends a Petri net to the $\law{SLAT}$-net which forgets about the multiplicity of the edges between a given source and transition.
Before moving on to the semantics of $\law{Q}$-nets, we discuss a property of the category $\Net{\law{Q}}$.

\begin{prop}\label{qnetcocomplete}
$\Net{\law{Q}}$ is cocomplete.
\end{prop}
\begin{proof}
We can construct $\Net{\law{Q}}$ as the comma category. 
\[  \Set \downarrow (\times \circ \Delta \circ M_\law{Q} ) \]
where $M_\law{Q} \maps \Set \to \Set$ is the monad corresponding to the Lawvere theory $\law{Q}$, $\Delta \maps \Set \to \Set \times \Set$ is the diagonal, and $\times \maps \Set \times \Set \to \Set$ is the cartesian product in $\Set$. An object in this category is a map
\[T \to M_\law{Q} S \times M_\law{Q} S\]
\noindent which corresponds to a pair of maps $s,t \maps T \to M_\law{Q} S$ which become the source and target maps of a $\law{Q}$-net. Morphisms in this comma category are commutative squares

\[ 
\xymatrix{ T\ar[r] \ar[d]_{f} & M_\law{Q} S \times M_\law{Q} S \ar[d]^{M_{\law{Q}} g \times M_{\law{Q}} g } \\
            T' \ar[r] & M_\law{Q} S' \times M_\law{Q} S'}
\]
giving a map of $\law{Q}$-nets $(f \maps T \to T',g \maps S \to S')$. The commutativity of the above square ensures that this map of $\law{Q}$-nets is well-defined. 

Theorem 3, Section 5.2 of \textit{Computational Category Theory} \cite{comp} says that given $S \maps A \to C$ and $T \maps B \to C$ then the comma category $(S \downarrow T)$ is cocomplete if
\begin{itemize}
\item $S$ is cocontinuous, and
\item $A$ and $B$ are cocomplete,
\end{itemize}
Because $\Set$ is cocomplete and the identity functor $1_{\Set} \maps \Set \to \Set$ preserves all colimits, we have that $\Set \downarrow (\times \circ \Delta \circ M_\law{Q} )$ is cocomplete. Because $\Net{\law{Q}}$ is equivalent to this category, it is cocomplete as well.
\end{proof}
\section{Generating Free Commutative Monoidal Categories}\label{CMC}
In this section we examine in detail the motivating example for the main result of this chapter, an adjunction generating the semantics of $\law{Q}$-nets for every Lawvere theory $\law{Q}$. This result can feel abstract on its own and the example of Petri nets provides invaluable intuition. A confident reader may skip this section, as it is not strictly necessary for the rest of the chapter.

The operational semantics for Petri nets will take the form of an adjunction
\[
\begin{tikzcd}
\Petri \ar[r,bend left,"F"] \ar[r,phantom,"\bot"] & \CMC \ar[l,bend left,"U"] 
\end{tikzcd}
\]
For a given Petri net $P$, this adjunction will be constructed in two steps: first the transitions of $P$ will be freely closed under a commutative monoidal sum and then freely closed under composition. This will take the form of factoring the adjunction into the composite
\[
\begin{tikzcd}
\Petri \ar[r,bend left,"\A"] \ar[r,phantom,"\bot",pos=.6] & \Grph(\CMon) \ar[l,bend left,"\backa",pos=.45] \ar[r,phantom,"\bot",pos=.4] \ar[r,bend left,"\B"] & \CMC \ar[l,bend left,"\backb",pos=.55].
\end{tikzcd}
\]
Here a left adjoint is indicated by a bullet on the left and a right adjoint is indicated by a bullet on the right. $\Grph(\CMon)$ is the category of graphs internal to $\CMon$.
\begin{defn}
A \define{commutative monoidal graph} is a graph
\[
\begin{tikzcd}
E \ar[r,shift left=.5ex,"s"] \ar[r, shift right=.5ex,"t",swap] & V
\end{tikzcd}
\]
where $E$ and $V$ are commutative monoids and $s$ and $t$ are commutative monoid homomorphisms. A morphism of commutative monoidal graphs is a tuple of commutative monoid homomorphisms $(f \maps E \to E', g \maps V \to V')$ making the diagrams
	\[
	\xymatrix{ 
		E \ar[d]_f  \ar[r]^-{s} & V \ar[d]^-{g} \\	
		E' \ar[r]_-{s'} & V' 
	}
	\qquad
	\xymatrix{ 
		E \ar[d]_f  \ar[r]^-{t} & V \ar[d]^-{g} \\	
		E' \ar[r]_-{t'} & V' 
	}
	\]commute. This defines a category $\Grph(\CMon)$ where objects are commutative monoidal graphs and morphisms are as above. In short, $\Grph(\CMon)$ is the category of graphs internal to $\CMon$.
\end{defn}
We will now define these adjunctions but omit the proofs that they are indeed well-defined adjunctions, as this follows from the more general results of Section \ref{gen}.
The left adjoint $\A \maps \Petri \to \Grph(\CMon)$ is defined as follows:
\begin{defn}\label{A}
Let 
\[ \A \maps \Petri \to \Grph(\CMon)\]
be the functor which sends a  Petri net 
\[
\begin{tikzcd}
P = T \ar[r,shift left=.5ex,"s"] \ar[r, shift right=.5ex,"t",swap] & \N[S]
\end{tikzcd}
\]
to the commutative monoidal graph
\[
\begin{tikzcd}
\A P = LT \ar[r,shift left=.5ex,"\phi^{-1}(s)"] \ar[r,shift right=.5ex,"\phi^{-1}(t)",swap] & LS
\end{tikzcd}
\]
where $L$ is the left adjoint of the adjunction in Definition \ref{N} and $\phi \maps \Hom (LT, LS) \xrightarrow{\sim} \Hom(T,RLS)$ is the natural isomorphism of that adjunction. $\A$ sends a morphism of Petri nets
\[(f \maps T \to T', g \maps S \to S')\]
to the morphism of commutative monoidal graphs given by 
\[(Lf \maps LT \to LT',Lg \maps LS \to LS') \]
\end{defn}
In words, $\A$ freely generates a commutative monoidal structure on the transitions of a Petri net and $\A$ uniquely extends each component of a Petri net morphism to a commutative monoid homomorphism. 
The right adjoint of this functor is non-trivial:
\begin{defn}
Let 
\[\backa \maps \Grph(\CMon) \to  \Petri \]
be the functor which sends a commutative monoidal graph
\[\begin{tikzcd}
Q=E \ar[r,shift left=.5ex,"s"] \ar[r, shift right=.5ex,"t",swap] & V
\end{tikzcd} \]
to the Petri net 
\[ 
\begin{tikzcd}
\backa Q = \bar{E} \ar[r,shift left=.5ex,"\bar{s}"] \ar[r, shift right=.5ex,"\bar{t}",swap] & \N[RV]
\end{tikzcd}
\]
$\bar{E}$ is defined as 
\[\bar{E}= \{(e,x,y) \in RE \times \N[RV] \times \N[RV] | R\epsilon_V(x)=s(e) \text{ and } R \epsilon_V(y)=t(e) \} \]
where $\epsilon_V$ is the counit of the adjunction $L \dashv R$. $\bar{s}$ and $\bar{t}$ are given by the projection of $\bar{E}$ onto its first and second coordinates respectively. $\backa$ sends a morphism of commutative monoidal graphs
\[(f\maps E \to E',g\maps V \to V')\]
to the morphism of Petri nets 
\[(h \maps \bar{E} \to \bar{E'}, Rg \maps RV \to RV') \]
where $h$ is the function which makes the assignment
\[(e,x,y) \mapsto (\phi(f)(e),\N[Rg](x),\N[Rg](y) ) \]
\end{defn}
\begin{rmk}\label{failure}
Petri nets must have a free commutative monoid of places, so it is necessary to regard $RV$ as the set of places for $\backa Q$ rather than having $V$ be the commutative monoid of places itself. The reader at this point may guess a simpler formula for the right adjoint $\backa$ which keeps the $RE$ as the set of transitions and uses the unit of $\N$ to construct the source and target maps. Unfortunately this construction is doomed to fail. For a commutative monoidal graph $\begin{tikzcd}E \ar[r,shift left=.5ex,"s"] \ar[r,shift right=.5ex,"t",swap] & V \end{tikzcd}$, suppose that the right adjoint $\A$ sends this graph to the Petri net
\[
\begin{tikzcd}
RE \ar[r,shift left=.5ex,"\eta \circ Rs"] \ar[r,shift right=.5ex,"\eta\circ Rt",swap]& \N[RV].
\end{tikzcd}
\]
A problem arises because this process unnaturally chunks the source and target of each transition. To see this consider the commutative monoidal graph 
\[
\begin{tikzcd}
Q= \N[\tau] \ar[r,shift left=.5ex] \ar[r,shift right=.5ex] & \N[\{a,b,c\}] \cong \N^3
\end{tikzcd}
\]The edge $\tau$ in $Q$ can be depicted as
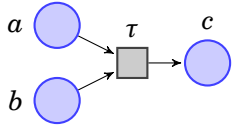
\[
\begin{tikzpicture}[node distance=1.3cm,>=stealth',bend angle=45,auto]
\node [place] (a) [label=left:$a$] at (0,0) {};
\node [place] (b) [label=left:$b$] at (0,-1) {};
\node [place] (c) [label=above:$c$] at (2,-.5) {};
\node [transition] (t) [label=$\tau$] at (1,-.5) {}
    edge [pre]          (a)
    edge [pre]          (b)
    edge [post]         (c);
\end{tikzpicture}
\]
With the above (faulty) description, $\backa Q$ is given by
\[
\begin{tikzcd}
\N[\tau] \ar[r,shift left=.5ex] \ar[r, shift right=.5ex] & \N[\N^3]
\end{tikzcd}
\]
To avoid confusion, we denote the outer sum in $\N[\N^3]$ by $\times$ and the sum in $\N^3$ by $+$. Then, the faulty $\backa$ would turn $\tau$ into the transition
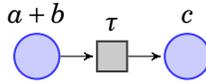
\[
\begin{tikzpicture}[node distance=1.3cm,>=stealth',bend angle=45,auto]
\node [place] (ab) [label=above:$a+b$] at (-1,0) {};
\node [place] (c) [label=above:$c$] at (1,0) {};
\node [transition] (t) [label=$\tau$] at (0,0) {}
    edge [pre]          (ab)
    edge [post]         (c);
\end{tikzpicture}
\]
To find a counit for this adjunction we seek a morphism
\[ \A \backa Q \to Q  \]
A morphism of this sort is defined by its assignment on generators. A natural choice of morphism sends the places $a+b$ to the sum of the places $a$ and $b$ using the counit of $L \dashv R$. However, then the assignment of $\tau \mapsto \tau$ does not respect the source of $\tau$ and is therefore not a morphism of commutative monoidal graphs. The problem is that we want the source of $\tau$ in $\A \backa Q$ to be $a \times b$ and not $a+b$. To fix this we force the source of $\tau$ to be $a\times b$ by upgrading $\tau$ to the tuple $(\tau, a \times b, c)$ in $\N[\tau] \times \N[\N^3] \times \N[\N^3]$. Now, the natural choice for the counit which sends $a \times b$ to $a +b$ respects the source of $\tau$.
\end{rmk}

The next part of the semantics adjunction for Petri nets freely generates the structure of a category on a given commutative monoidal graph. In Section \ref{gen} this is accomplished by rephrasing this construction in terms of free monoids. Here we provide an explicit description in the case of Petri nets. 

\begin{defn}
Let 
\[ \backb \maps \CMC \to \Grph(\CMon)\]
be the forgetful functor which sends a commutative monoidal category to its underlying commutative monoidal graph and a strict monoidal functor to its underlying morphism of commutative monoidal graphs. Then $\backb$ has a left adjoint 
\[ \B \maps \Grph(\CMon) \to \CMC\]
which sends a commutative monoidal graph
\[ 
\begin{tikzcd}
Q= E \ar[r,shift left=.5ex,"s"] \ar[r,shift right=.5ex,"t",swap] & V
\end{tikzcd}
\]
to the commutative monoidal category $\B Q$ with objects given by $V$ and morphisms generated inductively by the rules:
\begin{itemize}
    \item for every edge $e \in E$ a morphism $e \maps s(e) \to t(e)$,
    \item for every pair of morphisms $e \maps x \to y$ and $d \maps y \to z$, a morphism $d \circ e \maps x \to z$,
    \item for every object $v \in V$ a morphism $1_v \maps v \to v$
\end{itemize}
This defines an evident composition operation $\circ$ on $\B Q$.
There is also a sum on the $\B Q$ defined using the sum of $V$ on objects. If $e$ and $e'$ are edges of $Q$ then the morphisms $e \maps x \to y$ and $e' \maps x' \to y'$ already have a sum given by
\[e+e' \maps x +x' \to y+y' \]
 The morphisms of $\B Q$ are quotiented by the relations:
\begin{itemize}
    \item for all tuples of morphisms $(f \maps a \to b, g \maps b \to c, h \maps c \to d)$ 
    \[(f \circ g) \circ h = f \circ (g \circ h)\]
    \item for all morphisms $f \maps x \to y$ 
    \[1_y \circ f = f = f \circ 1_x \] 
    \item We require that composition is a commutative monoid homomorphism. For tuples of morphisms $(e \maps x \to y, d \maps y \to z, e' \maps x' \to y', d' \maps y' \to z')$ we can form their sum and composite in two different ways. We quotient the morphisms of $\B Q$ so that these are equal, i.e.
    \[(d \circ e) + (d' \circ e') = (d +d') \circ (e+e')\]
    \item We require that the assignment of identities is a commutative monoid homomorphism. For objects $x$ and $x'$ in $V$ we set
    \[1_{x+x'}=1_x +1_x'\]
\end{itemize}
\end{defn}
\section{Semantics Functors for Generalized Nets}\label{gen}
In this section we state our construction of semantics categories for $\Q$-nets; categories whose morphisms represent possible sequences of firings which can be performed using a given $\Q$-net.
Let $\Q$ be a Lawvere theory and let
\[\begin{tikzcd}
\Set \ar[r,bend left=45,"L"]\ar[r,phantom,"\bot",pos=.7] & \ar[l,bend left=45,"R"] \Mod(\Q, \Set)
\end{tikzcd}
\] be the adjunction it induces on $\Set$. In this section we will use this adjunction to construct an adjunction
\[
\begin{tikzcd}
  \Net{\Q} \ar[r,bend left,"F_{\Q}"] \ar[r,phantom,"\bot",pos=.6] & \ar[l,bend left,"U_{\Q}"] \Mod(\Q,\Cat)
\end{tikzcd}
\] which is analogous to the adjunction in Section \ref{CMC} and where $\Mod(\Q,\Cat)$ is the category of models of $\Q$ in $\Cat$. This adjunction factors as

\[
\begin{tikzcd}
\Net{\Q} \ar[r,bend left,"\A_\Q"] \ar[r,phantom,"\bot",pos=.6] & \ar[l,bend left, "{\backa}_\Q"] \Grph(\Mod(\Q)) \ar[r,bend left, "\B_\Q"]\ar[r,phantom,"\bot"] & \ar[l,bend left,"\backb_\Q"] \Mod(\Q,\Cat)
\end{tikzcd}
\]
where $\A_\Q$ freely generates a model of $\Q$ on the transitions of a given $\Q$-Net and $\B_\Q$ freely generates the structure of a category on a given $\Q$-graph.

These adjunctions are heavily motivated by the case when $\Q= \law{CMON}$ as this gives Petri nets. The main result of this chapter is as follows:

\begin{thm}\label{big}
There is an adjunction
\[
\begin{tikzcd}
  \Net{\Q} \ar[r,bend left,"F_{\Q}"] \ar[r,phantom,"\bot",pos=.6]& \ar[l,bend left,"U_{\Q}"] \Mod(\Q,\Cat).
\end{tikzcd}
\]
\end{thm}\noindent The left adjoint can be described using inference rules.
Let $P$ be the $\Q$-net
\[
\begin{tikzcd}
T \ar[r, shift left=1, "s"] \ar[r, shift right=1, "t",swap] & M_\Q S
\end{tikzcd}
\]
The objects of $F_{\Q}$ are given by $L_{\Q} S$. That is, for every morphism $o \maps n \to m$ in $\Q$ and every tuple of places $x_1, x_2, \ldots, x_n$ there is an object $\mathbf{o}(x_1, x_2, x_3, \ldots, x_n)$. For an equation of morphisms in $\Q$
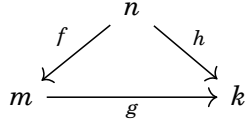
\begin{center}
\begin{tikzcd}
& n \ar[dl, "f",swap] \ar[dr,"h"] \\
m \ar[rr,"g",swap] & &k
\end{tikzcd}
\end{center}
the objects generated by each path must be equal. This means that there are $k$ equations of objects
\[ \mathbf{g_j} (\mathbf{f_1} (x_{1}, x_{2}, \ldots, x_{n}), \mathbf{f_2} (x_{1},x_{2}, \ldots, x_{n} ) , \ldots, \mathbf{f_m} (x_{1},x_{2},\ldots, x_{n}) ) = \mathbf{h_j} (x_{1}, x_2, \ldots, x_n)  \]
where the unlabeled index runs over the components of $f$ and the index $j$ runs over the components of $g$ and $h$. The morphisms of $F_Q P$ are generated inductively by the rules
\[
\infer{\tau \maps s(\tau) \to t(\tau) \in \Mor\, F_\Q P}{\tau \in T} \quad \quad \infer{1_x \maps x \to x \in \Mor\, F_\Q P}{x \in \Ob\, F_\Q P} \quad \quad \infer{g \circ f \maps x \to z \in \Mor \, F_\Q P}{f \maps x \to y \text{ and } g \maps y \to z \in \Mor\, F_\Q P }
\]
\[
\infer{\mathbf{o} (f_1,f_2, \ldots, f_n) \maps \mathbf{o} (x_1,x_2, \ldots,x_n) \to \mathbf{o} (y_1,y_2, \ldots, y_n)}{o \maps n \to 1 \in \Mor\, \Q \text{ and } f_1\maps x_1 \to y_1,f_2 \maps x_2 \to y_2, \ldots f_n \maps x_n \to y_n \in \Mor \, F_\Q P}
\]
and is quotiented to satisfy the following:
\begin{itemize}
\item The morphisms must satisfy the same equations that the objects satisfy. That is, for an equation of morphisms in $\Q$, the objects generated by each path must again be equal.

\item $\Mor\, F_\Q P$ is quotiented to satisfy the axioms of a category including the associative and unital laws
\[ (f \circ g) \circ h = f \circ (g \circ h) \text{  and  } 1_y \circ f = f = f \circ 1_x \]
for all morphisms $f$, $g$ and $h$ in $\Mor \, F_{\Q} P$.
\item $\Mor\, F_\Q P$ is quotiented so that the structure maps of a category (source, target, identity and composition) are $\Q$-model homomorphisms.
\end{itemize}

For a morphism of $\Q$-nets, $(f,g) \maps P \to P'$, the $\Q$-functor
\[ F_{\Q} (f,g) \maps F_{\Q} P \to F_{\Q} P' \]
is the unique extension of $f$ and $g$ which respects composition, unitality, and the operations of $\Q$.
The proof will require several lemmas. The first step is to show how $\Q$-nets freely generate graphs internal to the category of models of $\Q$.
\section{$\Q$-nets freely generate $\Q$-graphs}\label{freeQgraph}
In this section we show how $\Q$-nets freely generate graphs internal to $\Mod(\Q)$. In the following proofs we will write the monad $M_\Q$ as $R L$ and make use of the natural isomorphism
\[\phi \maps \hom(L X, Y) \xrightarrow{\sim} \hom(X, R Y)\]
for all sets $X$ and objects $Y$ in $\Mod(\Q)$.
\begin{defn}
Let\[ \A_\Q \maps \Net{\Q} \to \Grph(\Mod(\Q))\] be the functor which makes the assignment
\[
\begin{tikzcd}
T \ar[d,"f",swap] \ar[r, shift left=.5ex,"s"] \ar[r,shift right=.5ex,"t",swap] & RLS \ar[d,"RLg"{name=L}] & & LT \ar[d,"Lf"{name=R},swap] \ar[r, shift left=.5ex,"\phi^{-1}(s)"] \ar[r,shift right=.5ex,"\phi^{-1}(t)",swap] & LS \ar[d,"Lg"] \\
T' \ar[r, shift left=.5ex,"s'"] \ar[r,shift right=.5ex,"t'",swap] & RLS'& & LT' \ar[r, shift left=.5ex,"\phi^{-1}(s')"] \ar[r,shift right=.5ex,"\phi^{-1}(t')",swap] & LS' \ar[mapsto,from=L,to=R,shorten <=2ex,shorten >=2ex]
\end{tikzcd}
\]
on objects and morphisms.
\end{defn}

\begin{lem}
    $\A_\Q$ is well-defined.
\end{lem}
The next few proofs will make heavy use of the naturality equations for $\phi$ and its inverse:
\[  \phi(a \circ b \circ Lc) = R a \circ \phi(b) \circ c\] and
\[\phi^{-1} (Ra \circ b \circ c) = a \circ \phi^{-1}(b) \circ Lc.\]
\begin{proof}
First we show that $Le$ commutes with the source of $\A_\Q P$. This follows from the chain of equalities:
\begin{align*}
 & \phi^{-1} (s') \circ Lf \\
 &= \phi^{-1} (s' \circ f) \\
 &= \phi^{-1} (RLg \circ s)\\
 &= Lg \circ \phi^{-1} (s).
 \end{align*}
 A similar equation holds for the target maps.
\end{proof}
Let $G$ be the $\Q$-graph
\[\begin{tikzcd}
E \ar[r,shift left=.5ex,"s"] \ar[r,shift right=.5ex,"t",swap] & V.
\end{tikzcd} \]Because $V$ is not a free model of $\Q$, there is no obvious forgetful way to turn this into a $\Q$-net.  A first guess for the $\Q$-net $\backa_\Q (G)$ might be the $\Q$-net
\[\begin{tikzcd}RE \ar[r,shift left=.5ex,"\eta_{RV} \circ s"] \ar[r,shift right=.5ex,"\eta_{RV} \circ t",swap] & M_\Q (RV)\end{tikzcd} \]
where $\eta_{RV}$ is the unit of the monad $M_\Q$ applied to the set $RV$. However, as explained in Remark \ref{failure}, this fails to be a right adjoint. An alternative approach was suggested by Mike Shulman in the comments of an $n$Caf\'e blog post \cite{nCafe}. This solution was inspired by the construction of the free category on a tensor scheme introduced in \emph{The Geometry of Tensor Calculus I} \cite{scheme}. Instead of using $RE$ as the set of transitions, we use
    \[ \bar{E} = \{ (e,x,y) \in RE \times M_\Q RV \times M_\Q RV \, | \, s(e) = R\epsilon_{V} (x)\text{ and } t(e) = R \epsilon_{V}(y) \}\]
where $\epsilon_{V}$ is the $V$ component of the counit for $M_\Q$. Here and in what follows we are using $s$ to denote $Rs$ and $t$ to denote $Rt$ for notational simplicity. The source and target maps of the resulting $\Q$-net are given by the projections of $\bar{E}$ onto its second and third coordinates. The set $\bar{E}$ can be described formally using pullbacks.
\begin{defn}\label{backa}
Let
\[ \backa_{\Q} \maps \Grph(\Q)\to \Net{\Q} \]
be the functor which makes the assignment on objects and morphisms
\[
\begin{tikzcd}
E \ar[d,"f",swap] \ar[r, shift left=.5ex,"s"] \ar[r,shift right=.5ex, "t",swap] & V  \ar[d,"g"{name=L}] & \bar{E}\ar[d,"\bar{f}"{name=R},swap]\ar[mapsto,shorten <=1ex,shorten >=1ex,from=L,to=R] \ar[r, shift left=.5ex,"\bar{s}"] \ar[r, shift right=.5ex,"\bar{t}",swap] & M_\Q RV \ar[d,"M_\Q R g"]\\
E  \ar[r, shift left=.5ex,"s"] \ar[r,shift right=.5ex, "t",swap] & V   & \bar{E'}  \ar[r, shift left=.5ex,"\bar{s'}"] \ar[r, shift right=.5ex,"\bar{t'}",swap] & M_\Q RV' &
\end{tikzcd}
\]
where
\begin{itemize}
    \item $\bar{E}$ is the pullback of sets
    \[
    \begin{tikzcd}[column sep={4em,between origins}]
 & \ar[dl,"i",swap] \bar{E} \ar[dr,"j"] & \\
    RE \ar[dr,"{(s,t)}",swap] & & M_\Q RV \times M_\Q RV \ar[dl,"R\epsilon_V \times R\epsilon_V"] \\
    & RV\times RV &
    \end{tikzcd}
    \]
   where $(s,t)$ denotes the pairing of $s$ and $t$, and $\epsilon_{RV}\times \epsilon_{RV}$ denotes the cartesian product of the counits.
    \item $\bar{s} \maps \bar{E} \to M_\Q RV$ is the composite
    \[
    \begin{tikzcd}
    \bar{E} \ar[r,"j"] & M_\Q RV \times M_\Q RV \ar[r,"\pi_1"] & M_Q RV
    \end{tikzcd}
    \]
    and $\bar{t} \maps \bar{E} \to M_\Q RV$ is the composite
      \[
    \begin{tikzcd}
    \bar{E} \ar[r,"j"] & M_\Q RV \times M_\Q RV \ar[r,"\pi_2"] & M_Q RV
    \end{tikzcd}
    \]
    that is the maps which send an element $(e,x,y)$ of $\bar{E}$ to its second and third coordinates.
   \item $\bar{f} \maps \bar{E} \to \bar{E'}$ is induced by the universal property of $\bar{E}$ as shown below
 
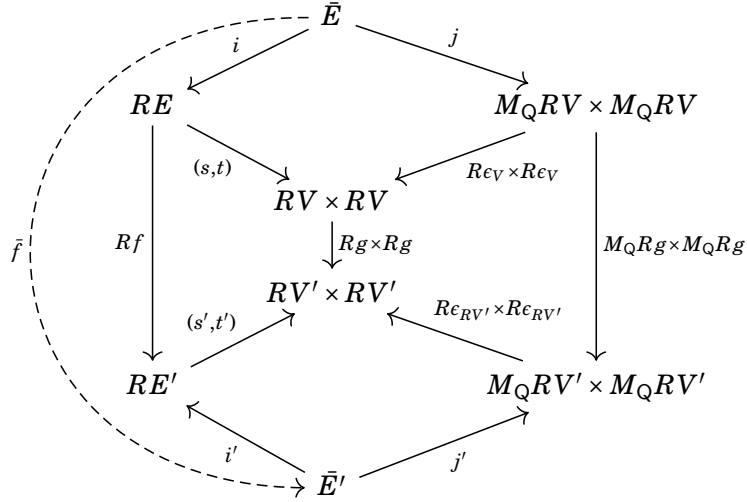
\[
    \begin{tikzcd}
 & \ar[dl,"i",swap] \bar{E} \ar[ddddd,dashrightarrow,bend right=90,looseness=2,"\bar{f}",swap] \ar[dr,"j"] & \\
    RE \ar[ddd,"Rf",swap]\ar[dr,"{(s,t)}",swap] & & M_\Q RV \times M_\Q RV \ar[dl,"R\epsilon_V \times R\epsilon_V"] \ar[ddd,"M_\Q Rg \times M_\Q Rg"] \\
    & RV\times RV \ar[d,"Rg \times Rg"] &\\
    & RV' \times RV' & \\
    RE' \ar[ur, "{(s',t')}"] & & M_\Q RV' \times M_\Q RV' \ar[ul, "R\epsilon_{RV'} \times R\epsilon_{RV'}",swap,pos=.7] \\
     & \bar{E'} \ar[ur,"j'",swap] \ar[ul,"i'"]&
    \end{tikzcd}
    \]
    More simply, $\bar{f}$ makes the assignment
    \[ (e,x,y) \mapsto (Rf(e),M_\Q Rg (x), M_\Q Rg (y) )\]
\end{itemize}
\end{defn}

\begin{lem}
    $\backa_\Q$ is well-defined.
\end{lem}

\begin{proof}
We must show that $(\bar{f}, Rg)$ is a well-defined morphism of $\Q$-nets.
$\bar{f}$ and $Rg$ commute with the source and target maps. Indeed, using the elementary descriptions we get that
\begin{align*}
    \bar{s'} \circ \bar{f} (\tau,x,y) &= \bar{s} ( Rf(\tau), M_\Q Rg(x), M_\Q Rg (y) ) \\
    &= M_\Q Rg (x) \\
    &= M_\Q Rg ( \bar{s} (\tau,x,y) )
\end{align*}
A similar equation holds for the target maps. $(\bar{f}, Rg)$ commutes with the identity maps:
\begin{align*}
    \bar{f} \circ \bar{e} (x) &= \bar{f} (Re (x), \eta_{RV} (x), \eta_{RV} (x) ) \\
    &=(Rf \circ Re(x), M_\Q Rg \circ \eta_{RV} (x), M_\Q Rg \circ \eta_{RV} (x) ) \\
    &= (Re' \circ Rg (x), M_\Q Rg \circ \eta_{RV} (x), M_\Q Rg \circ \eta_{RV} (x) )\\
    &= (Re' \circ Rg (x),\eta_{RV'} \circ Rg (x),\eta_{RV'} \circ Rg (x) )\\
    &= \bar{e'} \circ Rg (x)
\end{align*}
where the last two steps follow from naturality of $\eta$ and $(f,g)$ being a morphism of $\Q$-graphs.
\end{proof}

\begin{lem}
$\backa_\Q$ is a right adjoint to $\A_\Q$.
\end{lem}

\begin{proof}
Let $P$ be the $\Q$-net
\[
\begin{tikzcd}
T \ar[r,shift left=.5ex,"s"] \ar[r,shift right=.5ex,"t",swap] & RLS
\end{tikzcd}
\]
and $Q$ be the $\Q$-graph
\[
\begin{tikzcd}
E \ar[r,shift left=.5ex,"s'"] \ar[r, shift right=.5ex,"t'",swap] & V
\end{tikzcd}
\]
We define a natural isomorphism
\[ \Phi \maps  \Hom(\B P, Q) \xrightarrow{\sim} \Hom(P, \backb_\Q Q)\]
by the rule
\[
\begin{tikzcd}
LT \ar[d,"f",swap] \ar[r, shift left=.5ex,"\phi^{-1} (s)"] \ar[r, shift right=.5ex,"\phi^{-1}(t)",swap] & LS \ar[d,"g"{name=L}]  & & T \ar[d,"h"{name=R},swap] \ar[r, shift left=.5ex,"s"] \ar[r,shift right=.5ex,"t",swap] & RLS \ar[d,"RL \phi(g)"] \ar[mapsto,from=L,to=R,shorten <=3ex, shorten >=3ex]\\
E \ar[r,shift left=.5ex,"s'"] \ar[r, shift right=.5ex,"t'",swap]  & V & & \bar{E} \ar[r,shift left=.5ex,"\bar{s'}"] \ar[r,shift right=.5ex,"\bar{t'}",swap] & RLRV
\end{tikzcd}
\]
$h$ is defined by the universal property induced by $\bar{E}$ and the diagram
\[
\begin{tikzcd}[column sep={6em,between origins}]
& \ar[ddl] \bar{E} \ar[ddrr] & & \\
&                          & & \\
RE \ar[dr,"{(Rs',Rt')}",swap] &\ar[l,"\phi(f)",swap] T \ar[uu,dashed,"h"] \ar[rr,"{(RL \phi(g) \circ s, RL \phi(g) \circ t)}"] & & RLRV \times RLRV \ar[dll,"R\epsilon_V \times R\epsilon_V"] \\
& RV \times RV & &
\end{tikzcd}
\]
This diagram is well defined because $T$ is a competitor to the pullback $\bar{E}$ i.e.\ it makes the lowest triangle commute. Checking this amounts to showing that the bottom square commutes and this can be verified componentwise:
\begin{align*}
R\epsilon_V \circ RL \phi(g) \circ s &=R(\epsilon_V \circ L \phi(g)) \circ s \\
&= R( \phi^{-1} (1_RV) \circ L \phi(g) ) \circ s \\
&= R( \phi^{-1} ( 1_{RV} \circ \phi(g) ) ) \circ s \\
& = R( \phi^{-1} (\phi(g) ) ) \circ s \\
&= Rg \circ s \\
& = Rg \circ \phi( \phi^{-1} (s) ) \\
&= \phi (g \circ \phi^{-1} (s) ) \\
&= \phi(s' \circ f) \\
&= Rs' \circ \phi(f)
\end{align*}
and similar equations hold for the target maps. Therefore, $h$ is well defined.  Explicitly $h$ is the map which makes the assignment on transitions in $T$
\[ \tau \mapsto ( \phi(f)(\tau), RL \phi(g) \circ s(\tau), RL \phi(g) \circ t(\tau) ). \]
 $(h, \phi(g) )$ is a well-defined morphism of $\Q$-graphs by construction. The source and target functions map elements to their second and third coordinates so the equation
\[ \bar{s'} \circ h = RL \phi(g) \circ s\]
is true.

An inverse to $\Phi$,
\[\Phi^{-1} \maps \Hom ( P, \backb_{\Q} Q) \to \Hom( \B_{\Q} P , Q), \]
is defined as follows
\[
\begin{tikzcd}
T \ar[d,"h",swap] \ar[r, shift left=.5ex,"s"] \ar[r, shift right=.5ex,"t",swap] & RLS \ar[d,"RLg"{name=L}]  & & LT \ar[d,"\phi^{-1}(a)"{name=R},swap] \ar[r, shift left=.5ex,"\phi^{-1}(s)"] \ar[r,shift right=.5ex,"\phi^{-1}(t)",swap] & LS \ar[d," \phi^{-1}(g)"] \ar[mapsto,from=L,to=R,shorten <=2ex, shorten >=2ex]\\
\bar{E} \ar[r,shift left=.5ex,"\bar{s'}"] \ar[r, shift right=.5ex,"\bar{t'}",swap]  & RLRV & & E \ar[r,shift left=.5ex,"s'"] \ar[r,shift right=.5ex,"t'",swap] & V
\end{tikzcd}
\]$a$ is defined by the universal property of $\bar{E}$ and the diagram
\[
\begin{tikzcd}
& \ar[dl,"i",swap] \bar{E} \ar[dr,"j"] &  \\
RE \ar[dr,"{(Rs',Rt')}",swap] &\ar[l,dashed,"a",swap] T \ar[u,"h"] \ar[r,"b",dashed] &  RLRV \times RLRV \ar[dl,"R\epsilon_V \times R\epsilon_V"] \\
& RV \times RV & &
\end{tikzcd}
\]
To show that $(\phi^{-1}(a), \phi^{-1}(g))$ is a well defined morphism of $\Q$-graphs we perform the computation:
\begin{align*}
    s' \circ \phi^{-1}(a) &= \phi^{-1}(Rs' \circ a)\\
    &= \phi^{-1} (R \epsilon_V \circ \pi_1 \circ b) \\
    &= \phi^{-1} (R \epsilon_V \circ \bar{s'} \circ h)
\end{align*}
where $\pi_1 \maps RLRV \times RLRV \to RLRV$ is the projection and the last two steps follow from the definition of $\bar{s'}$ and commutativity of the above diagram. This can be reduced using the fact that $h$ commutes with the source and target of $P$ and $\backb_{\Q} Q$ and naturality of $\phi^{-1}$. Indeed,
\begin{align*}
    \phi^{-1} ( R \epsilon_V \circ \bar{s'} \circ h ) & = \phi^{-1} (R \epsilon_V \circ RLg \circ s) \\
    &=\phi^{-1} ( R(\epsilon_V \circ Lg) \circ s ) \\
    &= \epsilon_V \circ Lg \circ \phi^{-1} (s) \\
    &=\phi^{-1} (1_{RV}) \circ Lg \circ \phi^{-1} (s)\\
    &= \phi^{-1} (1_{RV} \circ g)\circ \phi^{-1} (s)\\
    &= \phi^{-1} (g) \circ \phi^{-1} (s)
\end{align*}
A similar equation holds for target so this is a well-defined morphism of $\Q$-graphs. $\Phi$ is a natural isomorphism if it is a natural and a bijection in the places component and the transitions component. The places component is only an application of $\phi$ so it is both natural and a bijection. For the transition component let $D \maps C \to \Set$ be the diagram
\[
\begin{tikzcd}
RE \ar[dr,"{(Rs,Rt)}",swap]& & M_{\Q} RV \times M_{\Q}RV \ar[dl,"R\epsilon_V \times R\epsilon_V"]\\
& RV \times RV & \\
\end{tikzcd}
\]
where $C$ is the walking cospan. Let $\Delta_T \maps C \to \Set$ be the constant diagram which sends every object to $T$ and every morphism to $1_T$. Then the universal property of $\bar{E}$ can be expressed as the natural isomorphism
\[\Psi \maps \mathrm{Nat} (\Delta_T, D) \xrightarrow{\sim}  \Hom (T, \bar{E})\]
where $\mathrm{Nat}(\Delta_T, D)$ denotes the set of natural transformations from $\Delta_T$ to $D$. With this description, the transition component of $\Phi$ can be described as follows
\[
\Phi \maps \, f \,\mapsto \,\Psi(\langle \phi(f), (RL \phi(g) \circ s, RL \phi(g) \circ t)  \rangle)
\]
where the angle brackets encase the components of a natural transformation. Similarly, the transition component of $\Phi^{-1}$ can be described as
\[
\Phi^{-1} \maps \, h\, \mapsto \,\phi^{-1}(\Psi^{-1} (h)_{RE})
\]
where the subscript $RE$ indicates that we take the $RE$ component of the natural transformation. With this description, we can verify that they are inverses on the transition component:
\begin{align*}
    f & \mapsto \Psi(\langle \phi(f), (RL \phi(g) \circ s, RL \phi(g) \circ t)  \rangle) \\
    & \mapsto \phi^{-1}(\Psi^{-1} (\Psi(\langle \phi(f), (RL \phi(g) \circ s, RL \phi(g) \circ t)  \rangle))_{RE} ) \\
    &= \phi^{-1}(\langle \phi(f), (RL \phi(g) \circ s, RL \phi(g) \circ t)  \rangle_{RE} )\\
    &= \phi^{-1}(\phi(f)) \\
    &= f
\end{align*}
and the other direction:
\begin{align*}
    h & \mapsto \phi^{-1} ( \Psi^{-1}(h)_{RE} ) \\
      & \mapsto \Psi( \langle \phi ( \phi^{-1} (\Psi^{-1} (h)_{RE})), (RL \phi(g) \circ s, RL \phi(g) \circ t) \rangle ) \\
      &= \Psi \Psi^{-1} \langle h, (RL \phi(g) \circ s, RL \phi(g) \circ t \rangle_{RE} \\
      &= \langle h, (RL \phi(g) \circ s, RL \phi(g) \circ t \rangle_{RE}\\
      &=h.
\end{align*}
The transition component of $\Phi$ and $\Phi^{-1}$ are natural because they are made up of components which are individually natural transformations.
\end{proof}
The next step in the proof of Theorem \ref{big} is to construct an adjunction between $\Grph(\Mod(\Q))$ and $\Mod(\Q,\Cat)$, i.e.\ the free graph construction internal to $\Mod(\Q)$.

\section{Free Categories Internal to $\Mod(\Q)$}\label{freecatinternal} In this section we will construct an adjunction
\[\begin{tikzcd}
\Mod(\Q,\Grph) \ar[r,bend left,"F_\Q"] & \Mod(\Q,\Cat) \ar[l,bend left,"U_\Q"]
\end{tikzcd} \]
to complete the proof of Theorem \ref{big}. A general property of algebraic theories $\law{P}$ and $\Q$ is that models of $\law{P}$ in the category of models of $\Q$ are the same as models of  $\Q$ in the category of models of $\law{P}$. In particular for a Lawvere theory $\Q$, a model of $\Q$ in $\Cat$ is the same as a category internal to $\Mod(\Q)$ and a model of $\Q$ in $\Grph$ is the same as a graph internal to $\Mod(\Q)$. This extends to an equivalences of categories
\[\Mod(\Q,\Cat) \cong \Cat(\Mod(\Q))\text{ and } \Mod(\Q,\Grph) \cong \Grph(\Mod(\Q)).  \]
Therefore, in this section we instead construct an adjunction
\[\begin{tikzcd}\Grph(\Mod(\Q) \ar[r,bend left,"F_\Q"] & \ar[l,bend left,"U_\Q"] \Cat(\Mod(\Q)) \end{tikzcd} \]
i.e.\ we construct free categories internal to $\Mod(\Q)$. This adjunction is not new, and was first given in \cite{baues1997}. In this section we obtain it by applying a construction of Lack \cite{lack} to the following monoidal category:

\begin{defn}
Let $\Grph(\Mod(\Q))(V)$ be the monoidal category where
\begin{itemize}
    \item objects are given by graphs $s,t \maps E \to V$ in $\Mod(\Q)$,
    \item morphisms are given by maps $f \maps E \to E'$ making the diagram
    \[
    \begin{tikzcd}
    & \ar[dl,"s",swap] E \ar[dd,"f"] \ar[dr,"t"] & \\
    V & & V\\
    &\ar[ul,"s'"] E'\ar[ur,"t'",swap]&
    \end{tikzcd}
    \]commute.
    \item monoidal product is given by chosen pullbacks. That is, for spans
    \[ \begin{tikzcd}
 & \ar[dl,"a",swap] E \ar[dr,"b"] & & & \ar[dl,"c",swap] F \ar[dr,"d"] & \\
 V & & V & V & & V
\end{tikzcd}
\]
their monoidal product is the chosen pullback
\[
\begin{tikzcd}
& &\ar[dl] E \times_{V} F \ar[dr] & & \\
  &\ar[dl,"a",swap]E\ar[dr,"b"]&  &\ar[dl,"c",swap]F\ar[dr,"d"]& \\
V & &V & & V
\end{tikzcd}
\]
On morphisms $f \maps E \to E'$ and $g \maps F \to F'$ is the unique map
\[(f,g) \maps E \times_V F \to E' \times_V F'\]
induced by the universal property of $E' \times_V F'$.
\end{itemize}
\end{defn}A monoid in this monoidal category is a span $s,t \maps E \to V$ along with multiplication and unit maps
\[ \circ \maps E \times_V E \to E \text{ and } e \maps V \to E\]
satisfying associativity and unitality. Interpreting $\circ$ as composition and $e$ as the map assigning identity morphisms, the monoid becomes a category internal to $\Mod(\Q)$. Indeed, a category with object model $V$ is exactly a monoid in the category $\Grph(\Mod(\Q)) (V)$ \cite{betti1996formal}.
Therefore it suffices to show that $\Grph(\Mod(\Q)(V)$ admits a free monoid construction. In Proposition \ref{freecat} we handled the case when $\Mod(\Q)$ is $\Set$ using the geometric series formula
\[F(G) = 1 + G + G^2 + \ldots = \sum_{n \geq 0} G^n \]
in the category of graphs $\Grph(V)$ over a fixed vertex set. The multiplication of $F(G)$ has the type
\[ \sum_{n \geq 0} G^n \times \sum_{m \geq 0} G^m \to \sum_{n \geq 0} G^n.\]
Because products distribute over coproducts in $\Grph(V)$, we can factor the multiplication as
\[ \sum_{n \geq 0} G^n \times \sum_{m \geq 0} G^m \xrightarrow{\sim} \sum_{m,n \geq 0} G^m \times G^n \xrightarrow{m} \sum_{n \geq 0} G^n
\]
where $m$ is the unique map induced by the natural concatenation morphisms $G^m \times G^n \to G^{m+n}$. In the general case, because products may not distribute over coproducts in the category $\Grph(\Mod(\Q))(V)$, this multiplication won't in general exist and a different construction of free monoids is necessary. Luckily Lack offers an alternative construction which replaces the coproducts of the previous approach with filtered colimits and reflexive coequalizers \cite{lack}. These colimits are sifted so they commute with finite products and a multiplication map based on concatenation can be naturally defined.
\begin{thm}\label{lack}{\normalfont [Lack]}
Let $(C, \otimes)$ be a monoidal category with
\begin{itemize}
    \item finite limits,
    \item countable colimits, and
    \item the functors $- \otimes A$ and $A \otimes -$ preserve reflexive coequalizers and colimits of countable chains.
\end{itemize}
Then $C$ admits a free monoid construction, that is, a left adjoint to the forgetful functor
\[ \mathrm{Mon} (C) \to C\]
that sends every monoid to its underlying object of $C$.
\end{thm}\noindent We now apply this Theorem to $\Grph(\Mod(\Q))(V)$.
\begin{prop}\label{freeinternalv}For each object $V$ in $\Mod(\Q)$, there is an adjunction
\[
\begin{tikzcd}
\Grph(\Mod(\Q))(V) \ar[r,bend left, "\B_V"] & \Cat(\Mod(\Q))(V) \ar[l,bend left,"\backb_V"]
\end{tikzcd}
\]
where $\Cat(\Mod(\Q))(V)$ is the category of categories internal to $\Mod(\Q)$ whose model of objects is $V$.
\end{prop}

\begin{proof}
The hypotheses of Theorem \ref{lack} require that the following conditions hold: \begin{itemize}
    \item $\Grph(\Mod(\Q)) (V)$ has finite limits and countable colimits. $\Mod(\Q)$ has these limits and colimits as shown in Theorem 3.4.5 of \cite{Borceux2}. The corresponding limits and colimits in $\Grph(\Mod(\Q))(V)$ are computed on the edges of each graph.
    \item The product of $\Grph(\Mod(\Q))(V)$ preserves colimits of countable chains and reflexive coequalizers. This is true because colimits of countable chains and reflexive coequalizers are sifted colimits so they commute with finite products.
\end{itemize}
Applying Theorem \ref{lack} to the category $\Grph(\Mod(\Q))(V)$ gives the desired result.
\end{proof}
To complete the adjunction exhibiting the operational semantics of $\Q$-nets, we need to remove the dependence on the model of vertices $V$. To accomplish this, we use the Grothendieck construction \cite{Borceux2}.

\begin{defn}
Let
\[\Grph(\Mod(\Q))(-) \maps \Mod(\Q) \to \CAT\]
be the functor which sends an object $V$ to the category $\Grph(\Mod(\Q))(V)$ of graphs over $V$. For a morphism $f \maps V \to W$ in $\Mod(\Q)$, let
\[\Grph(\Mod(\Q))(f) \maps \Grph(\Mod(\Q))(V) \to \Grph(\Mod(\Q))(W) \]
be the functor which makes the assignment
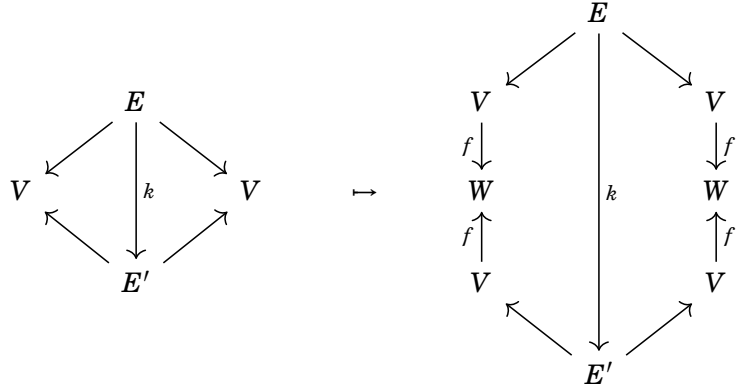
\[
\begin{tikzcd}
 & & & & &\ar[dl] E \ar[dddd,"k"] \ar[dr] &\\
 & \ar[dl] E \ar[dd,"k"] \ar[dr] & & & V \ar[d,"f",swap] & & V \ar[d,"f"] \\
 V & & V & \mapsto & W & & W \\
 & \ar[ul] E' \ar[ur] & & & V \ar[u,"f"] & & V \ar[u,"f",swap] \\
 & & & & &\ar[ul] E' \ar[ur] &\\
\end{tikzcd}
\]on objects and morphisms. Let
\[ \Cat(\Mod(\Q))(-) \maps \Mod(\Q) \to \CAT \]
be the functor which sends an object $V$ to the category of small categories internal to $\Mod(\Q)$ with object model of $\Q$ given by $V$. For a morphism $f \maps V \to W$, let
\[\Cat(\Mod(\Q))(f) \maps \Cat(\Mod(\Q))(V) \to \Cat(\Mod(\Q))(W) \]
be the functor which which makes the assignment
\[
\begin{tikzcd}
\Mor\, C \ar[d,"k",swap] \ar[r, shift left=.5ex,"s"] \ar[r, shift right=.5ex,"t",swap] & V \ar[d,"\mathrm{id}"{name=L}] &  & \Mor\, C  \ar[d,"k"{name=R},swap]\ar[r, shift left=.5ex,"s"] \ar[r, shift right=.5ex,"t",swap] & V\ar[d,"\mathrm{id}"] \ar[r,"f"] & W \ar[d,"\mathrm{id}"] \ar[mapsto,from=L,to=R,shorten <=3ex, shorten >=3ex] \\
\Mor\, C' \ar[r, shift left=.5ex,"s'"] \ar[r,shift right=.5ex,"t'",swap] & V & & \Mor \, C'  \ar[r, shift left=.5ex,"s'"] \ar[r, shift right=.5ex,"t'",swap] & V \ar[r,"f",swap] & W
\end{tikzcd}
\]on the underlying graphs of objects and morphisms.
\end{defn}\noindent Proposition \ref{freeinternalv} is reframed in this context.
\begin{prop}\label{Grothendieckinternal}
The family of adjunctions
\[
\begin{tikzcd}
\Grph(\Mod(\Q)) (V) \ar[r,bend left, "\B_V"] & \Cat(\Mod(\Q))(V) \ar[l,bend left,"{\backb}_V"]
\end{tikzcd}
\]
form components of natural transformations
\[ C \maps \Grph(\Mod(\Q))(-) \Rightarrow \Cat(\Mod(\Q))(-) \text{  and  } \backb \maps \Cat(\Mod(\Q))(-) \Rightarrow \Grph(\Mod(\Q))(-) \]
Furthermore, $C$ and $\backb$ form an adjoint pair in the 2-category $[\Mod(\Q), \Cat]$ where
\begin{itemize}
    \item objects are functors $F \maps \Mod(\Q) \to \Cat$,
    \item morphisms are natural transformations $\alpha \maps F \Rightarrow G$ whose components $\alpha_c \maps F(c) \to G(c)$ are functors and,
    \item 2-morphisms are modifications $\gamma \maps \alpha \to \beta$. That is, for every object $c$ in $\Mod(\Q)$ a natural transformation of the type
    \[
    \begin{tikzcd}
    F(c) \ar[d,bend right=50,"\alpha_c"{name=L},swap] \ar[d, bend left=50,"\beta_c"{name=R}] \\
    G(c) \ar[Rightarrow, from=L, to=R,shorten <= 1.7ex, shorten >= 1.7ex,"\gamma_c"].
    \end{tikzcd}
    \]
\end{itemize}
\end{prop}

\begin{proof}
For naturality, it suffices to show that the squares
\[
\begin{tikzcd}
\Grph(\Mod(\Q))(V) \ar[d,"\Grph(\Mod(\Q))(f)",swap] \ar[r,"\B_V"] & \Cat(\Mod(\Q))(V) \ar[d,"\Cat(\Mod(\Q))(f)"] & \\
\Grph(\Mod(\Q))(W) \ar[r,"\B_W",swap] & \Cat(\Mod(\Q))(W)
\end{tikzcd}
\]
and
\[
\begin{tikzcd}
\Grph(\Mod(\Q))(V) \ar[d,"\Grph(\Mod(\Q))(f)",swap]  &\ar[l,"{\backb}_V",swap] \Cat(\Mod(\Q))(V) \ar[d,"\Cat(\Mod(\Q))(f)"]\\
\Grph(\Mod(\Q))(W) & \ar[l,"{\backb}_V"]  \Cat(\Mod(\Q))(W)
\end{tikzcd}
\] commute. This is verified by direct computation. To show that $\B$ and $\backb$ are an adjoint pair we need the following fact: $\B$ is a left adjoint to $\backb$ in $[\Mod(\Q),\Cat]$ if and only if the components
\[
\begin{tikzcd}
\Grph(\Mod(\Q)) (V) \ar[r,bend left, "\B_V"] & \Cat(\Mod(\Q))(V) \ar[l,bend left,"{\backb}_V"]
\end{tikzcd}
\]
form an adjoint pair in $\Cat$. The counit-unit definition of adjunction requires that we have modifications $\epsilon \maps \B \circ \backb \to 1_{\Cat(\Mod(\Q))(-)}$ and $\eta \maps 1_{\Grph(\Mod(\Q))(-)} \to \backb  \circ \B$ satisfying the snake equations. Unpacking this gives components $\epsilon_V \maps \B_V \circ \backb_V \Rightarrow 1_{\Cat(\Mod(\Q))(V)}$ and $\eta_V \maps \backb_V \circ \B_V \Rightarrow 1_{\Grph(\Mod(\Q))(V)}$ satisfying the snake equations. This is equivalent to each component being an adjunction. However, Theorem \ref{lack} says that each component is an adjunction so the claim is shown.
\end{proof} \noindent So far we have the diagram
\[
\begin{tikzcd}
\Mod(\Q)\, \ar[r,bend left=70,"\Grph(\Mod(\Q))(-)"{name=U}] \ar[r,bend right=70,"\Cat(\Mod(\Q))(-)"{name=D},swap] & \ar[Rightarrow,from=U, to=D,shorten <= 1.7ex, shorten >= 1.7ex,bend right,"\B"description] \ar[Rightarrow,from=D, to=U,shorten <= 1.7ex, shorten >= 1.7ex,bend right,"\backb"description]\Cat
\end{tikzcd}
\]of adjoint 1-cells in
$[\Mod(\Q),\Cat]$. We apply the Grothendieck construction to this diagram to get
\[
\begin{tikzcd}
\int \Grph(\Mod(\Q))(-) \ar[r,bend left,"\int \B"] & \int \Cat(\Mod(\Q))(-) \ar[l,bend left,"\int \backb"]
\end{tikzcd}
\]The Grothendieck construction is a 2-functor $\int \maps [\Mod(\Q), \CAT]\to \CAT/\Mod(\Q)$ where $\CAT$ denotes the 2-category of large categories, functors, and natural transformations. When composed with the forgetful 2-functor $\CAT/\Mod(\Q) \to \CAT$ which remembers only the domain of each functor, we obtain the composite
\[\int \maps [\Mod(\Q),\CAT] \to  \CAT\]
which we denote as $\int$ in an abuse of notation.

A fundamental fact is that every 2-functor preserves adjunctions. Therefore the above diagram is an adjunction. Moreover, the following proposition shows that it is the adjunction we are looking for.

\begin{prop}\label{equiv}
The category $\int \Grph(\Mod(\Q))(-)$ is equivalent to $\Grph(\Mod(\Q))$ and the category $\int \Cat(\Mod(\Q))(-)$ is equivalent to $\Mod(\Q,\Cat)$.
\end{prop}

\begin{proof}
\noindent $\int \Grph(\Mod(\Q))(-)$ has
\begin{itemize}
    \item pairs $(V,V \leftarrow E \to V)$ as objects and,
    \item pairs $(f \maps V \to V', g \maps E \to E')$ such that the diagram
    \[
    \begin{tikzcd}
    & \ar[dl] E \ar[ddd,"g"] \ar[dr] & \\
    V \ar[d,"f",swap]& & V\ar[d,"f"]\\
    V' & & V'\\
    & \ar[ur] E' \ar[ul] &
    \end{tikzcd}
    \]in $\Mod(\Q)$ commutes as morphisms.

\end{itemize}
An equivalence $\int \Grph(\Mod(\Q))(-) \xrightarrow{\sim} \Grph(\Mod(\Q))$ sends $(V,V \leftarrow E \to V)$ to the graph $\begin{tikzcd}
E \ar[r, shift left=.5ex] \ar[r, shift right=.5ex,swap] & V
\end{tikzcd}$
and a morphism $(f,g)$ to the evident morphism of graphs $(f \maps E \to E', g \maps V \to V')$. \noindent $\int \Cat(\Mod(\Q))(-)$ has
\begin{itemize}
    \item pairs $(V,C)$ where $C$ is a category over $V$ as objects and,
    \item pairs $(f\maps V \to V', g \maps C \to C')$ where $g$ is an object fixing functor from $\Cat(\Mod(\Q))(f) (C)$ to $C'$ as morphisms.
\end{itemize}
An equivalence $\int \Cat(\Mod(\Q))(-) \xrightarrow{\sim} \Cat(\Mod(\Q)$ is given by sending objects $(V,C)$ to their second component and morphisms $(f,g)$ to the functor whose object component is $f$ and whose morphism component is the morphism component of $g$.
\end{proof}\noindent We denote the compositions of $\int \B$ and $\int \backb$ with the above equivalences by $\B_{\Q}$ and $\backb_{\Q}$ respectively.

\noindent \textbf{Proof of Theorem \ref{big}.}
The composite adjunction $F_\Q \dashv U_\Q$ is constructed by setting $F_\Q = \B_{\Q} \circ \A_{\Q}$ and $U_\Q = \backb_\Q \circ \backa_\Q$. \hfill \qedsymbol\\
\smallskip

\chapter{Compositionality of Q-Nets}\label{openQNets}In this chapter we use the operational semantics developed in the previous chapter to develop a compositional theory of the behavior of $\Q$-nets. In Section \ref{openQnet} we define ``open" $\Q$-nets, i.e.\ $\Q$-nets equipped with input and output ports. Open $\Q$-nets are glued together via pushout. There is a symmetric monoidal double category called $\Open(\Net{\Q})$ where the horizontal morphisms are open $\Q$-nets and horizontal composition is pushout. In Section \ref{compopseminternal} we show how the operational semantics functor
\[F_\Q \maps \Net{\Q} \to \cat{\Q} \]
of Theorem \ref{big} lifts to a symmetric monoidal double functor
\[\Open(F) \maps \Open(\Net{\Q}) \to \Open(\cat{\Q}). \] Because composition in $\Open(\Net{\Q})$ is gluing, functoriality of this double functor gives relationships between the behavior of a $\Q$-net and the behaviors of its components. In Section \ref{functionalqnetsection} we define the black-boxing of an open $\Q$-category: a profunctor that encapsulates the morphisms from the input ports to the output ports. In Theorem \ref{blackboxqnet} we show that black-boxing lifts to a lax double functor
\[\blacksquare \maps \Open(\cat{\Q}) \to \Prof. \]
We define functional open $\Q$-nets, based off of the functional Petri nets of Zaitsev and Sleptsov \cite{zaitsev1997}, as open $\Q$-nets for which every input port is a source and every output port is a sink. Functional open $\Q$-nets generalize the functional open graphs of Definition \ref{functgraph}. In Theorem \ref{blackboxfunctqnet}, we show that $\blacksquare \circ \Open(F)$ is strictly functorial on functional open $\Q$-nets. This gives a straightforward expression of the compositionality of functional open $\Q$-nets which does not suffer from combinatorial explosion.
\section{Open $\Q$-nets}\label{openQnet} In this section we define open $\Q$-nets and construct a double category
$\Open(\Net{\Q})$ with open $\Q$-nets as horizontal 1-morphisms. As in Chapter \ref{graphs}, to define open $\Q$-nets, we require a functor $L \maps \Set \to \Net{\Q}$ that maps any set $S$ to a $\Q$-net with $S$ as its set of places, and we need $L$ to be a left adjoint.

\begin{defn}
Let $L \maps \Set \to \Net{\Q}$ be the functor defined on sets and functions as follows:
\[
	\xymatrix{
		X \ar[d]_f^{ \quad \mapsto} & \emptyset \ar[d] \ar@<-.5ex>[r] \ar@<.5ex>[r] & M_{\Q}[X] 		\ar[d]^{M_{\Q}[f]} \\
		Y & \emptyset  \ar@<-.5ex>[r] \ar@<.5ex>[r] & M_{\Q}[Y]
	}
\]
where the unlabeled maps are the unique maps of their type.
\end{defn}

\begin{lem}
\label{L}
The functor $L$ has a right adjoint $R \maps \Net{\Q} \to \Set$ that acts as follows on
$\Q$-nets and their morphisms:
\[
	\xymatrix{
		T \ar[d]_f  \ar@<-.5ex>[r]_-{t} \ar@<.5ex>[r]^-{s} & M_{\Q}[S] \ar[d]^{M_{\Q}[g] \quad \mapsto} &  S \ar[d]^{g} \\
		T' \ar@<-.5ex>[r]_-{t'} \ar@<.5ex>[r]^-{s'} & M_{\Q}[S] & S'.
	}
\]
\end{lem}

\begin{proof}
For any set $X$ and $\Q$-net  $P = (s,t \maps T \to M_{\Q}[S])$ we have natural isomorphisms
\[
\begin{array}{ccl}
\hom_{\Net{\Q}}\big(L(X), \xymatrix{ T \ar@<-.5ex>[r]_-t \ar@<.5ex>[r]^-s & M_{\Q}[S] }\big)
& \cong & \hom_{\Net{\Q}}\big(\xymatrix{ \emptyset \ar@<-.5ex>[r] \ar@<.5ex>[r] & M_{\Q}[X] },
\xymatrix{ T \ar@<-.5ex>[r]_-t \ar@<.5ex>[r]^-s & M_{\Q}[S] }\big) \\
&\cong& \hom_{\Set}(X,S) \\
&\cong& \hom_{\Set}\big(X, R(\!\xymatrix{ T \ar@<-.5ex>[r]_-t \ar@<.5ex>[r]^-s & M_{\Q}[S]}\!) \big).
\quad \qedhere \end{array}
\]
\end{proof}

An ``open" $\Q$-net is a $\Q$-net $P$ equipped with maps from two sets $X$ and $Y$ into its set of places, $RP$.  We can write this as a cospan in $\Set$ of the form
\[ \xymatrix{ & RP & \\
					X \ar[ur] & & Y. \ar[ul] }
\]
Using the left adjoint $L$ we can reexpress this as a cospan in $\Net{\Q}$, and this gives our official definition:

\begin{defn}
\label{defn:openpetri}
An \define{open $\Q$-net} is a diagram in $\Net{\Q}$ of the form
\[ \xymatrix{ & P & \\
					LX \ar[ur]^i & & LY \ar[ul]_o }
\]
for some sets $X$ and $Y$.  We sometimes write this as $P \maps X \nrightarrow Y$ for short.
\end{defn}

We now introduce the main object of study for this section: the double category $\Open(\Net{\Q})$, which has open $\Q$-nets as its horizontal 1-cells.

	\begin{thm}
	\label{openQnetdouble}
		There is a symmetric monoidal double category $\Open(\Net{\Q})$ for which:
		\begin{itemize}
			\item objects are sets
			\item vertical 1-morphisms are functions
			\item horizontal 1-cells from a set $X$ to a set $Y$ are open $\Q$-nets
				
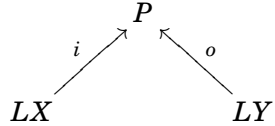
\[ \xymatrix{ & P & \\
					LX \ar[ur]^i & & LY \ar[ul]_o } \]
			\item 2-morphisms $\alpha \maps P \Rightarrow P'$ are commutative diagrams
			  
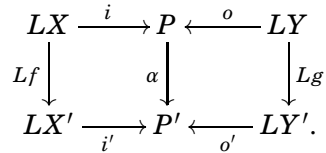
\[ \xymatrix{	LX \ar[r]^i \ar[d]_{Lf} & P \ar[d]_{\alpha} & LY \ar[l]_o \ar[d]^{Lg} \\
			  	LX' \ar[r]_{i'}  & P'  & LY'. \ar[l]^{o'} } \]
			  	in $\Net{\Q}$.
	\end{itemize}
Composition of vertical 1-morphisms is the usual composition of functions.   Composition of horizontal
1-cells is composition of cospans via pushout: given two horizontal 1-cells
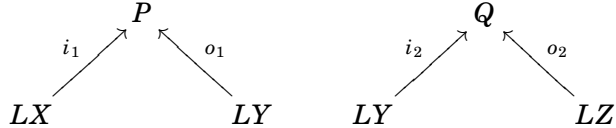
\[ \xymatrix{ & P & & & Q & \\
	LX \ar[ur]^{i_1} & & LY \ar[ul]_{o_1} & LY \ar[ur]^{i_2} & & LZ \ar[ul]_{o_2} }\]
their composite is given by this cospan from $LX$ to $LZ$:
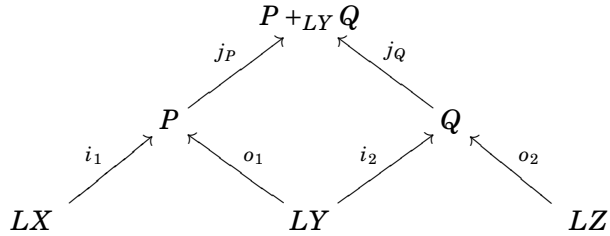
\[ \xymatrix{
	&   & P+_{LY} Q  &  & \\
	& P \ar[ur]^{j_P} &  & Q \ar[ul]_{j_Q} & \\
	LX \quad \ar[ur]^{i_1} & & LY \ar[ul]_{o_1}  \ar[ur]^{i_2} & & \quad LZ \ar[ul]_{o_2} }\]
where the diamond is a pushout square.  The horizontal composite of 2-morphisms
 
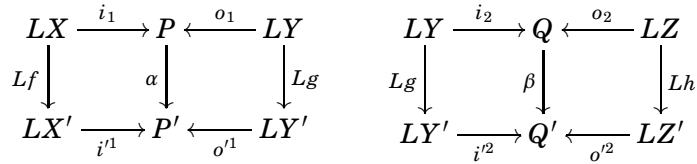
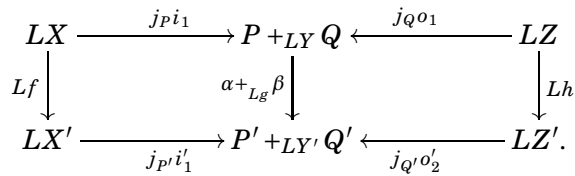
\[ \xymatrix{	LX \ar[r]^{i_1} \ar[d]_{Lf} & P \ar[d]_{\alpha} & LY \ar[l]_{o_1} \ar[d]^{Lg} \\
			  	LX' \ar[r]_{i'^1}  & P'  & LY' \ar[l]^{o'^1} }
\qquad
 \xymatrix{	LY \ar[r]^{i_2} \ar[d]_{Lg} & Q \ar[d]_{\beta} & LZ \ar[l]_{o_2} \ar[d]^{Lh} \\
			  	LY' \ar[r]_{i'^2}  & Q'  & LZ' \ar[l]^{o'^2} } \]
is given by
 \[ \xymatrix{	LX \ar[rr]^{j_P i_1} \ar[d]_{Lf} && P+_{LY} Q \ar[d]_{\alpha+_{{}_{Lg}} \beta}
 && LZ \ar[ll]_{j_Q o_1} \ar[d]^{Lh} \\
LX' \ar[rr]_{j_{P'} i'_1}  && P'+_{LY'} Q'  && LZ'. \ar[ll]^{j_{Q'} o'_2} } \]
Vertical composition of 2-morphisms is done using composition of functions.
The symmetric monoidal structure comes from coproducts in $\Set$ and $\Net{\Q}$.
\end{thm}

\begin{proof}
We construct this symmetric monoidal double category using Lemma \ref{Courser}. This lemma requires that $\Net{\Q}$ has all colimits and this was proved in Proposition \ref{qnetcocomplete}.
\end{proof}
\begin{expl}Setting $\Q$ equal to the theory of commutative monoids gives open Petri nets. This double category and its properties are explored in detail in the paper \cite{open} by the author and Baez.
\end{expl}

In the remarks following Example \ref{ksafe}, we exploited the functoriality of Definition \ref{netf} to describe functorial relationships between different categories of $\Q$-nets. This functoriality can be extended to produce symmetric monoidal double functors between double categories of open $\Q$-nets. 
\begin{prop}\label{openqnetfunct}
Every morphism of Lawvere theories $f\maps Q \to R$ induces a lax symmetric monoidal lax double functor
\[\Open(\Net{f}) \maps \Open(\Net{\Q}) \to \Open(\Net{\law{R}}) \]
\end{prop}

\begin{proof}
From Definition \ref{netf} and Proposition \ref{L} we have the following diagram
\[\begin{tikzcd}
\Net{\Q} \ar[r,"\Net{f}"]&  \Net{\law{R}} \\
\Set \ar[r,equals] \ar[u,"L_{\law{R}}"] & \ar[u,"L_\Q",swap] \Set
\end{tikzcd} \]and it is straightforward to verify that this diagram commutes up to natural isomorphism. The desired double functor is obtained by applying Lemma \ref{openfunctoriality} to this diagram.\end{proof}
In Section \ref{QNet} we constructed the following diagram of functors
\[
\begin{tikzcd}
\Net{\law{SLAT} } &            &\\
\Petri \ar[r,"\Net{b}"] \ar[u,"\Net{a}"] & \Net{\Z} &\\
\PreNet \ar[u,"\Net{c}"] \ar[r,"\Net{d}",swap] & \Net{\law{GRP}} \ar[u,"\Net{e}",swap]
\end{tikzcd}\]
This diagram of functors can be continued one step further to categories of open $\Q$-nets via Proposition \ref{openqnetfunct}:
\[
\begin{tikzcd}
\Open(\Net{\law{SLAT} }) &            &\\
\Open(\Petri) \ar[r,"\Open(\Net{b})"] \ar[u,"\Open(\Net{a})"] & \Open(\Net{\Z}) &\\
\Open(\PreNet) \ar[u,"\Open(\Net{c})"] \ar[r,"\Open(\Net{d})",swap] & \Open(\Net{\law{GRP}}) \ar[u,"\Open(\Net{e})",swap]
\end{tikzcd}\]
This diagram says that the above functors between $\Q$-nets can be extended in a coherent way to open $\Q$-nets.

\section{Compositionality of the Operational Semantics for $\Q$-nets}\label{compopseminternal}

In Theorem \ref{big} we saw how a $\Q$-net $P$ gives a $\Q$-category $F_\Q(P)$, and in Theorem \ref{openQnetdouble} we constructed a double category $\Open(\Net{\Q})$ of open $\Q$-nets.  Now we construct a double category $\Open(\cat{\Q})$ of ``open $\Q$-categories" and a map
\[             \Open(F_\Q) \maps \Open(\Net{\Q}) \to \Open(\cat{\Q}) .\]
This can be seen as providing an operational semantics for open $\Q$-nets in which any open $\Q$-net is mapped to the $\Q$-category it presents. 
The key is this commutative diagram of left adjoint functors:
\[ \xymatrix{ \Set \ar[r]^L \ar[dr]_{FL} & \Net{\Q} \ar[d]^F\\
	& \cat{\Q}}
\]
where $FL$ sends any set to the free $\Q$-category on this set: $FLX$ has $M_\Q[X]$, the free $\Q$-model on $X$, as its set of objects and only identity morphisms.  Using Lemma \ref{Courser},  we can produce two symmetric monoidal double categories from this diagram.   We have already seen one: $\Open(\Net{\Q})$ obtained from the left adjoint $L$. We now obtain $\Open(\cat{\Q})$ from the left adjoint $FL$.

\begin{thm}
\label{thm:openCMC}
There is a symmetric monoidal double category $\Open(\cat{\Q})$ for which:
\begin{itemize}
		\item objects are sets
		\item vertical 1-morphisms are functions
			\item horizontal 1-cells from a set $X$ to a set $Y$ are \define{open $\Q$-categories} $C \maps X \nrightarrow Y$, that is, cospans in $\cat{\Q}$ of the form
				\[ \xymatrix{ & C & \\
					L'X \ar[ur]^i & & L'Y \ar[ul]_o } \]
		where $C$ is a $\Q$-category and $i,o$ are strict $\Q$-functors,
			\item 2-morphisms $\alpha \maps C \Rightarrow C'$ are commutative diagrams in $\cat{\Q}$
			of the form
			  \[ \xymatrix{	L'X \ar[r]^i \ar[d]_{L'f} & C \ar[d]_{\alpha} & L'Y \ar[l]_o \ar[d]^{L'g} \\
			  	L'X' \ar[r]_{i'}  & C'  & L'Y'. \ar[l]^{o'} } \]
	\end{itemize}
and the rest of the structure is given as in Lemma \ref{Courser}.
\end{thm}

\begin{proof}
To apply Lemma \ref{Courser} to the functor $FL \maps \Set \to \cat{\Q}$ we just need to check that $\cat{\Q}$ has finite colimits.  First note that
\[\cat{\Q} \simeq \mathsf{Mod}(\mathsf{Q},\Cat).\]
The cocompleteness of this category then follows from various classical results, some listed in the introduction of a paper by Freyd and Kelly \cite{FK}.  More recently, Trimble \cite[Prop.\ 3.1]{Trimble} showed that for any Lawvere theory $\Q$ and any cocomplete cartesian category $\X$ with finite products distributing over colimits, the category of finite-product-preserving functors $\mathsf{Mod}(\Q,\X)$ is cocomplete.
\end{proof}

The functor $F \maps \Net{\Q} \to \cat{\Q}$ induces a map sending open $\Q$-nets to open $\Q$-categories. This map is actually part of a symmetric monoidal double functor. 
\begin{thm}
\label{thm:functoriality}
There is a symmetric monoidal double functor
\[   \Open(F) \maps \Open(\Net{\Q}) \to \Open(\cat{\Q}) \]
that is the identity on objects and vertical 1-morphisms, and makes the following assignments on horizontal 1-cells and 2-morphisms:
  \[ \xymatrix{	LX \ar[r]^i \ar[d]_{Lf} & P \ar[d]_{\alpha} & LY \ar[l]_o \ar[d]^{Lg \qquad {\Huge{\mapsto}} \quad } & &
  L'X \ar[r]^{Fi} \ar[d]_{L'f} & FP \ar[d]_{F\alpha} & L'Y \ar[l]_{Fo} \ar[d]^{L'g }
   \\
			  	LX' \ar[r]_{i'} & P'  & LY' \ar[l]^{o'} & &
			  	L'X' \ar[r]_{Fi'}  & FP'  & L'Y'. \ar[l]^{Fo'} } \qedhere
\]
\end{thm}

\begin{proof}
We apply Lemma \ref{openfunctoriality} to the commutative square
\[
\begin{tikzcd}
\Net{\Q} \ar[r,"F"] & \cat{\Q} \\
\Set \ar[u,"L"] \ar[r,equals] & \Set \ar[u,"F \circ L",swap]
\end{tikzcd}
\]
of left adjoint functors.
\end{proof}

We can think of the $\Q$-category $FP$ as providing an operational semantics for the $\Q$-net $P$: morphisms in this category are processes allowed by the $\Q$-net. The above theorem says that this semantics is compositional.  That is, if we write $P$ as a composite (or tensor product) of smaller open $\Q$-nets, $FP$ will be the composite (or tensor product) of the corresponding open $\Q$-categories.

\section{Black-boxing and Functional Open $\Q$-nets}\label{functionalqnetsection}
In this section we show how the black-boxing functor introduced in Theorem \ref{blackbox} can be applied to the categorical operational semantics of open $\Q$-nets. We also introduce functional open $\Q$-nets and prove that black-boxing preserves composition of functional open $\Q$-nets up to isomorphism. To define black-boxing for open $\Q$-nets we make use of the following functor.
\begin{defn}
Let \[U \maps \cat{\Q} \to \Cat\]
be the forgetful functor that regards every $\Q$-category as an ordinary category and every $\Q$-functor as an ordinary functor.
\end{defn}

\begin{thm}\label{blackboxqnet}
There is a lax double functor
\[\blacksquare \maps \Open(\cat{\Q}) \to \Prof  \]
that
\begin{itemize}
\item sends sets $X$,$Y$ and functions $f \maps X \to Y$ to the discrete categories and functors between them.
\item An open $\Q$-category
\[
\begin{tikzcd}
& C & \\
FLX \ar[ur,"i"] & & \ar[ul,"j",swap] FLY
\end{tikzcd}
\]
to the profunctor
\[\blacksquare(C) \maps UFLX \times UFLY \to \Set \]
given by $\blacksquare(C)(x,y)= UFP(i(x),j(y))$.
\item A 2-cell of open $\Q$-nets
\[
\begin{tikzcd}
FLX \ar[r,"i"]\ar[d,"Lf",swap] & C\ar[d,"g"] & FLY\ar[d,"Lh"] \ar[l,"o",swap]\\
FLX' \ar[r,"i'",swap] & D & FLY' \ar[l,"o'"]
\end{tikzcd}
\]
is sent to the 2-cell of profunctors
\[
\begin{tikzcd}
UFLX \times UFLY \ar[rr,"UFLf \times UFLh"] \ar[dr,"\blacksquare(C)"{name=L},swap] & & \ar[dl,"\blacksquare(D)"{name=R}]  UFLX' \times UFLY'\ar[from=L,to=R,Rightarrow,"\alpha",shorten <= 5ex,shorten >= 5ex,yshift=1.5ex] \\
& \Set &
\end{tikzcd}\]
where the components of $\alpha$ 
\[ \alpha_{x,y} \maps UC(ix,oy) \to UD(i'fx,o'hy)\]
are given by pointwise application of the functor $Ug$.
\end{itemize}
\end{thm}

\begin{proof}
We factor $\blacksquare$ into three parts
\[\Open(\cat{\Q}) \hookrightarrow \Csp(\cat{\Q}) \xrightarrow{\Open(U)} \Csp(\Cat) \xrightarrow{\blacksquare_{g}} \Prof\]
where
\begin{itemize}
    \item the inclusion $\Open(\cat{\Q}) \hookrightarrow \Csp(\cat{\Q})$ is given by pointwise application of $FL \maps \Set \to \cat{\Q}$ on sets and functions and is given by the identity on horizontal morphisms and $2$-cells.
    \item The functor $\Open(U)$ is given by applying Lemma \ref{openfunctoriality} to the square
    \[\begin{tikzcd}
    \cat{\Q} \ar[r,"U"] & \Cat \\
    \cat{\Q} \ar[u,equals] \ar[r,"U",swap] & \Cat \ar[u,equals]
    \end{tikzcd} \]
    Explicitly it is given by pointwise application of $U$ everywhere. Note that because $U$ does not preserve finite colimits, this double functor will preserve horizontal composition laxly.
    \item $\blacksquare_g$ is an extension of the black-boxing functor 
    \[\blacksquare \maps \Open(\Cat) \to \Prof \]
    of Theorem \ref{blackbox} to the domain $\Csp(\Cat)$. Explicitly, $\blacksquare_g$ is the identity on categories and functors and the black-boxing of horizontal morphisms and $2$-cells is exactly as in Theorem \ref{blackbox}.
\end{itemize} 
The desired double functor is obtained by composing the three double functors above. \end{proof}
The phenomenon of Example \ref{loop} persists for the black-boxing of open $\Q$-nets: the double functor $\blacksquare$ is lax rather than strict because profunctor composition only accounts for firing sequences which go from the first component to the second and do not come back. Functional $\Q$-nets are a class of $\Q$-nets for which firing sequences on a composite open Petri net only flow from the the first component to the second. Functional Petri nets were first introduced by Zaitsev and Sleptsov \cite{zaitsev1997}. Their definition generalizes straightforwardly to $\Q$-nets.
\begin{defn}
Let $P\maps X \to Y$ be the open $\Q$-net
\[\begin{tikzcd}
& P & \\
LX \ar[ur,"i"] & & LY \ar[ul,"j",swap]
\end{tikzcd}
\]
An element $x \in X$ is a \define{source} if $i(x)$ is not the target of any transition in $P$ and an element $y \in Y$ is a \define{sink} if $j(y)$ is not the source of any transition in $P$. $P \maps X \to Y$ is \define{functional} if every $x \in X$ is a source and every $y \in Y$ is a sink.
\end{defn}
For functional open $\Q$-nets profunctor composition accounts for all firing sequences on a composite.  Therefore black-boxing preserves horizontal composition up to isomorphism when the component $\Q$-nets are functional.
\begin{thm}\label{blackboxfunctqnet}
The composite double functor
\[\blacksquare \circ \Open(F) \maps \Open(\Net{\Q}) \to \Prof \]
preserves horizontal composition of functional $\Q$-nets up to isomorphism.
\end{thm}

 \begin{proof}
Let $P \maps X \to Y$ and $\Q \maps Y \to Z$ be functional open $\Q$-nets. Then their black-boxings are equipped with a composition comparison 
\[\alpha \maps \blacksquare \circ \Open(F) (P) \circ \blacksquare \circ \Open(F) (Q) \to \blacksquare \circ \Open(F) (P \circ Q) \]
with components 
\[ \alpha_{x,z} \maps \int_{y \in M_\Q(Y)} \blacksquare(\Open(P))(x,y) \times \blacksquare(\Open(Q))(y,z) \to \blacksquare(\Open(P) \circ \Open(Q)) (x,z) \]
given by sending a pair of morphisms $(g,f)$ to their composite $g \circ f$ in $\Open(P) \circ \Open(Q)$. Because $P$ and $Q$ are functional, every morphism from $x$ to $z$ in $\Open(P) \circ \Open(Q)$ is of this form. Therefore the composition comparison is an isomorphism.
\end{proof}

\begin{rmk}
Let $\dRel$ be the double category where objects are sets, vertical morphisms are functions, horizontal morphisms are relations, and 2-cells are rectangles
\[\begin{tikzcd}
A \ar[d,"f",swap] \ar[r,squiggly,"R"] & B \ar[d,"g"] \\
A' \ar[r,squiggly,"R'",swap] & B
\end{tikzcd} \]
such that there is an inclusion
\[(f \times g) \circ R \subseteq R'. \]
In \cite{open}, the author and Baez construct a different black-boxing functor
\[ \blacksquare \maps \Open(\Petri) \to \dRel\]
that sends an open Petri net $LX \xrightarrow{i} P \xleftarrow{o} LY$ to its reachability relation i.e.\ the relation 
$\blacksquare(P) \subseteq \N[X] \times \N[Y]$ given by
\[ \blacksquare(P)= \{(x,y) \in \N[X] \times \N[Y] \, | \, \text{ there exists a morphism }f\maps x \to y\text{ in }FP(ix,oy) \}.\] The reachability semantics functor is the decategorification of the black-boxing functor of Theorem \ref{blackboxqnet}. Any profunctor $P \maps C \times D^{\op} \to \Set$ can be turned into a relation $R \subseteq \Ob\, C \times \Ob\, D$ given by
\[
R = \{ (x,y) \, | \, P(x,y) \text{ is nonempty}. \} \]
When the profunctor $\blacksquare_{\law{CMON}}(FLX \xrightarrow{Fi} FP \xleftarrow{o} FLY)$ is turned into a relation in this way it becomes the reachability relation for $P$. This process may be called decategorification as it arises from the change of enrichment $\Set \to \{0,1\}$ where $\{0,1\}$ is boolean monoid regarded as a monoidal category. The change of enrichment is a functor which sends a set to $1$ if it is non-empty and extends to a double functor
\[\Prof \to \dRel. \]
To obtain the black-boxing of Open Petri nets, the black-boxing of this thesis is composed with the above decategorification to obtain the reachability semantics double functor of \cite{open}.
\end{rmk}

In \cite{zaitsev2005}, Zaitsev gives a polynomial time algorithm for decomposing a Petri net into functional open Petri nets. This decomposition provides a speedup for computing invariants of Petri nets. This chapter offers a formal language and general language to understand this strategy of decomposition in a larger context.

 \chapter{Operational Semantics of Enriched Graphs}\label{enriched}
  The algebraic path problem is a generalization of the shortest path problem to probability, computing, matrix multiplication, and optimization \cite{tarjan1981unified,foote2015kleene}. Let $R$ be the quantale of positive real numbers $([0,\infty], \mathrm{min}, +)$. A weighted graph is regarded as an $R$-matrix, and the shortest paths of this graph are computed as the operational semantics studied in this chapter. The algebraic path problem allows $R$ to vary, and gets solutions to other problems of a similar flavor within the same framework. In Section \ref{algpathproblem} we review the relevant definitions for quantales and the algebraic path problem. In Section \ref{freecatinternalsec} we generalize the free category construction of Proposition \ref{freecat} to graphs enriched in $R$. In Theorem \ref{freecatinternal} we obtain for every quantale $R$ an adjunction
\[
\begin{tikzcd}
\Mat_R \ar[r,bend left,"F"] \ar[r,phantom,"\bot"] & \RCat \ar[l,bend left,"U",pos=.55] 
\end{tikzcd} 
\]
between $R$-matrices and categories enriched in $R$. Remarkably, the free $R$-category $F(M)$ on an $R$-matrix $M$ is the solution to its algebraic path problem. The adjunction above encapsulates the universal property of this solution. Note that the functors in the above adjunction are referred to as $F_R$ and $U_R$ in the introduction to disambiguate from the other functors in this thesis. However, because this chapter is self contained we drop the subscript to reduce notational clutter.

\section{The Algebraic Path Problem}\label{algpathproblem}
The networks considered in this chapter will be parameterized by commutative quantales.
\begin{defn}
A \define{quantale} is a monoidal closed poset with all joins. Explicitly, a quantale is a a poset $R$ with an associative, unital, and monotone multiplication $\cdot \maps R \times R\to R$ such that
\begin{itemize}
    \item all joins, $\sum_{i \in I} x_i$, exist for arbitrary index set $I$ and
    \item  $\cdot$ preserves all joins, i.e.\ 
    \[a \cdot \sum_{i \in I}  x_i = \sum_{i \in I} a \cdot x_i\]
    for all joins over an arbitrary index set $I$.
\end{itemize}
A quantale is commutative if its multiplication operation, $\cdot$, is commutative.
\end{defn}\noindent A motivating example of such a quantale is the poset $[0,\infty]$ with $+$ as its monoidal product and with join given by infimum. Note that this poset is equipped with the reverse of the usual ordering on $[0,\infty]$. Fong and Spivak show how the shortest path problem on this quantale computes the shortest paths between all pairs of vertices in a given $[0,\infty]$-weighted graph \cite[\S 2.5.3]{fong2019invitation}. Other motivating examples include the rig $([0,1],\mathrm{sup},\times)$ (whose algebraic path problem corresponds to most likely path in a Markov chain) and the powerset of the language generated by an alphabet (whose algebraic path problem corresponds to the language decided by a nondeterministic finite automata (NFA))\cite{foote2015kleene}.
\begin{defn}
For a commutative quantale $R$ and sets $X$ and $Y$, an \define{$R$-matrix} $M \maps X \to Y$ is a function $M \maps X \times Y \to R$. For $R$-matrices $M \maps X \to Y$ and $N \maps Y \to Z$, their matrix product $M N$ is defined by the rule
\[MN(i,k) = \sum_{j\in Y} M(i,j) N(j,k)\]
\end{defn}
\noindent If $R$ is a commutative quantale, $R$-matrices form a quantale as well. 
\begin{defn}
Let $\RMat(X)$ be the set of $X$ by $X$ matrices $M \maps X \times X \to R$. $\RMat(X)$ is equipped with the partial order where $M \leq N$ if and only if $M(i,j) \leq N(i,j)$ for all $i,j \in X$.
\end{defn}
\begin{prop}
$\RMat (X)$ is a quantale with
\begin{itemize}
    \item join given by pointwise sum of matrices,
    \item and multiplication given by matrix product.
\end{itemize}
\end{prop}
\noindent The proof of this proposition is left to the reader. All the required properties of $\RMat(X)$ follow from the analogous properties in $R$.

A square matrix $M\maps X \times X \to R$ represents a complete $R$-weighted graph whose vertex set is given by $X$. 
\begin{defn} Let $M\maps X \times X \to R$ be a square matrix. A \define{vertex} of $M$ is an element $i\in X$. An \define{edge} of $M$ is a tuple of vertices $(a,b) \in X \times X$. A \define{path} in $M$ from $a_0$ to $a_n$ is a list of adjacent edges $p=\large((a_0,a_1),(a_1,a_2), \ldots,(a_{n-1},a_{n}), \large)$. The \define{weight} of $p$ is defined as the product
\[l(p) = \Pi_{i=0}^{n-1} M(a_i,a_{i+1}) \]
in $R$. For vertices $i,j \in X$, let
\[P_{ij}^M = \{\text{ paths in $M$ from $i$ to $j$ } \}  \]
\end{defn}Let $i$ and $j$ be vertices of a square matrix $M \maps X \times X \to R$. The algebraic path problem asks to compute the quantity
\[\sum_{p \in P_{ij}^M} l(p) \]
in the quantale $R$. If $R$ is the quantale $([0,\infty],\inf,+)$ then the weight of an edge $M_{ij}$ represents the distance between vertex $i$ and vertex $j$ and the weight of a path $l(p)$ represents the total distance traversed by $p$. Summing the weights of all paths between a pair of vertices corresponds to finding the path with the minimum weight. For example, the algebraic path problem asks to compute the length of the shortest path in the case when $R$ is $([0,\infty],\inf,+)$.

A more tractable framing of the algebraic path problem can be found by considering matrix powers. The entries of $M^2$ are given by 
\[M^2(i,j) = \sum_{l \in X} M(i,l)M(l,j)= \mathrm{inf}_{l \in X}\{ M(i,l)+ M(l,j)\}. \]
Because $M(i,l)$ and $M(l,j)$ represent the distance from $i$ to $l$ and from $l$ to $j$, this infimum computes the cheapest way to travel from $i$ to $j$ while stopping at some $l$ in between. More generally, the entries of $M^n$ for $n\geq 0$ represent the shortest paths between nodes of your graph that occur in exactly $n$ steps. To compute the shortest paths which can occur in any number of steps, we must take the infimum of the matrices $M^n$ over all $n \geq 0$. This pattern replicates for other choices of quantale. Therefore, the \define{algebraic path problem} seeks to compute 
\begin{equation}\label{pathproblem}
    F(M) = \sum_{n \geq 0} M^n
\end{equation}
where $M$ is an $R$-matrix. The following table summarizes some instances of the algebraic path problem for different choices of $R$. Fink provides an explanation of the algebraic path problems for $([0,\infty],\leq)$ and $\{T,F\}$ and Foote provides an explanation for the quantales $([0,1],\leq)$ and $(\mathcal{P}(\Sigma),\subseteq)$
\cite{fink1992survey,foote2015kleene}.
\smallskip
\begin{center}
\begin{tabular}{ |c|c|c|c| } 
\hline
\textbf{poset} & \textbf{join} & \textbf{multiplication} & \textbf{solution of path problem}\\
 \hline
 $([0,\infty],\geq)$ & $\inf$& $+$ & shortest paths in a weighted graph \\ 
 \hline
 $([0,\infty],\leq)$&$\sup$& $\inf$ & maximum capacity in the tunnel problem\\
 \hline
 $([0,1],\leq)$ & $\sup$ & $\times$ & most likely paths in a Markov process\\ 
 \hline
 $\{T,F\}$ & $\mathrm{OR}$ & $\mathrm{AND}$ & transitive closure of a directed graph \\ 
 \hline
 $(\mathcal{P}(\Sigma^*),\subseteq)$ & $\bigcup$ & concatenation & decidable language of a NFA \\
 \hline
\end{tabular}
\end{center}
\smallskip
Note that in this table, $\mathcal{P}(\Sigma^*)$ denotes the power set of the language generated by an alphabet $\Sigma$.
\section{The Algebraic Path Problem Functor}\label{freecatinternalsec}
Equation \ref{pathproblem} is known to category theorists by a different name: the free monoid on $M$. Framing it in this way gives a categorical proof of existence and uniqueness of $F(M)$. A classic result from \cite[\S V11]{maclane} gives a construction of free monoids. MacLane's construction is defined as an adjunction into a category of internal monoids.
\begin{defn}
Let $(C, \otimes, I)$ be a monoidal category. A \define{monoid internal to $C$} is an object $A$ of $C$ equipped with morphisms
\[m \maps A \otimes A \to A \text{ and } i \maps I \to M \]
satisfying the axioms of associativity and unitality expressed as commutative diagrams. A \define{monoid homomorphism} from a monoid $A$ to a monoid $B$ is a morphism $f: A \to B$ in $C$ which commutes with the maps $m$ and $i$ of each monoid. Let $\mathsf{Mon}(C)$ be the category where objects are monoids internal to $C$ and morphisms are their homomorphisms.
\end{defn}
\begin{prop}[MacLane]\label{Mac}Let $(C, \otimes, I)$ be a monoidal category with countable coproducts such that tensoring on both sides preserves these coproducts. Then there is an adjunction
\[\begin{tikzcd} \ \ \ C\ar[r,bend left,"F"] \ar[r,phantom,"\bot",pos=.6]& \ar[l,bend left,"U",pos=.45] \mathsf{Mon}(C)\end{tikzcd} \]
whose left adjoint is given by the countable coproduct
\begin{equation}\label{freemon}
F(X) = \sum_{n \geq 0} X^n.
\end{equation}
\end{prop}
The poset $\RMat(X)$ when viewed as a category satisfies the hypotheses of Proposition \ref{Mac} and admits a free monoid construction.
\begin{prop}\label{fixedmon}
There is an adjoint pair
\[\begin{tikzcd}\RMat(X) \ar[r,bend left,"F_{X}"] &\mathsf{Mon}(\RMat(X)) \ar[l,bend left,"U_{X}"] \end{tikzcd} \]
where $F_{X}$ is the monotone map which produces the solution to the algebraic path problem on a matrix and $U_{X}$ is the natural forgetful map.
\end{prop}
\begin{proof}
Because $\RMat(X)$ is a quantale, it can be regarded as a monoidal category with all coproducts such that tensoring distributes over these coproducts. The result follows from applying Proposition \ref{Mac} and noticing that Equation \ref{freemon} matches Equation \ref{pathproblem} in the case when $C=\RMat$.
\end{proof}

\noindent Monoids internal to $\RMat(X)$ are $R$-enriched categories.
\begin{defn}
An \define{$R$-enriched category} $C$ with object set $X$ consists of an element $C(x,y)$ in $R$ for every $x,y \in X$ such that
\begin{itemize}
    \item $1 \leq C(x,y)$ (the identity law),
    \item and $C(x,y) C(y,z) \leq C(x,z)$ (the composition law).
\end{itemize}
$R$-enriched categories will be referred to as $R$-categories. Let $\RCat(X)$ be the poset whose elements are $R$-categories with object set $X$. For $R$-categories $C$ and $D$ in $\RCat(X)$,
\[C \leq D \leftrightarrow C(i,j) \leq D(i,j) \quad \forall i,j \in X. \]
\end{defn}
\begin{prop}
$\mathsf{Mon}(\RMat(X))$ is isomorphic to $\RCat(X)$, the poset of categories enriched in $R$ with object set $X$.
\end{prop}
\begin{proof}
The isomorphism in question assigns a matrix $M \maps X \times X \to R$ to the $R$-category with $\hom(x,y)=M(x,y)$. The identity law follows from the inequality $1 \leq M$ and the inequality $M^2 \leq M$ implies that for all $y \in X$,
\[\sum_{y \in X} M(x,y) M(y,z) \leq M(x,z) \] The composition law follows from the fact that any element of $R$ is less than a join which contains it.
\end{proof}

\noindent Proposition \ref{fixedmon} says that each matrix valued in $R$ has a unique, universally characterized solution to the algebraic path problem: namely the free $R$-category on that matrix. This adjunction can be extended to matrices over an arbitrary set. 
\begin{defn}\label{pushforward}
Let $f: X\to Y$ be a function and let $M: X \times X \to R$ be an $R$-matrix. Then the \define{pushforward} of $M$ along $f$ is the matrix $f_* (M) \maps Y \times Y \to R$ defined by 
\[f_*(M) (y,y') = \sum_{(x,x') \in (f \times f)^{-1}(y,y')} M(x,x'). \]
\end{defn}

\begin{defn}\label{matr}
Let $\RMat$ be the category where objects are square matrices $M \maps X \times X \to R$ on a set $X$ and where a morphism from $M \maps X \times X \to R$ to $N \maps Y \times Y \to R$ is a function $f \maps X \to Y$ satisfying 
\[f_*(M) \leq N. \]
Let $\RCat$ be the full subcategory of $\RMat$ consisting of matrices satisfying the axioms of an $R$-category.
\end{defn}

\begin{thm}\label{Grothendieckenriched}
The free monoid construction of Proposition \ref{fixedmon} extends to an adjunction
\[\begin{tikzcd}\RMat \ar[r,bend left,"F"] \ar[r,phantom,"\bot"] & \RCat \ar[l,bend left,"U"] .\end{tikzcd} \]
\end{thm}

\begin{proof}
Let $A \maps \Set \to \Cat$ be the functor which sends a set $X$ to the poset $\RMat(X)$ regarded as a category and sends a function $f \maps X \to Y$ to the pushforward functor
\[f_* \maps \RMat(X) \to \RMat(Y).\]Analogously, let $B \maps \Set \to \Cat$ be the functor which sends a set $X$ to the poset $\RCat(X)$ and sends a function $f$ to its pushforward functor. The functors $F_{X}$ form the components of a natural transformation $\mathbf{F} \maps A \Rightarrow B$ and the functors $U_{X}$ form the components of a natural transformation $\mathbf{U} \maps B \Rightarrow A$. Furthermore, these natural transformations form an adjoint pair in the $2$-category $[\Set, \Cat]$ of functors $\Set \to \Cat$, natural transformations between them, and modifications. $\mathbf{F}$ and $\mathbf{U}$ are adjoint because an adjoint pair in $[\Set,\Cat]$ is a pair of natural transformations which are adjoint in each component. To summarize, we have a pair of adjoint natural transformations 
\[
\begin{tikzcd}
\Set\, \ar[r,bend left=70,"A"{name=U}] \ar[r,bend right=70,"B"{name=D},swap] & \ar[Rightarrow,from=U, to=D,shorten <= 1.7ex, shorten >= 1.7ex,bend right,"\mathbf{F}"description] \ar[Rightarrow,from=D, to=U,shorten <= 1.7ex, shorten >= 1.7ex,bend right,"\mathbf{U}"description]\Cat
\end{tikzcd}
\]\noindent A restriction of the Grothendieck construction \cite{Borceux} defines a 2-functor
\[\int \maps [\Set,\Cat] \to \CAT \]
where $\CAT$ is the 2-category of large categories. Because every 2-functor preserves adjunctions, the above diagram maps to an adjunction
\[ \begin{tikzcd}\int A \ar[r,bend left,"\int \mathbf{F}"] & \int B \ar[l,bend left,"\int \mathbf{U}"] .\end{tikzcd} \]
The result follows from the equivalences $\int A \cong \RMat$ and $\int B \cong \RCat$. The desired functors $F$ and $U$ are obtained by composing $\int \mathbf{F}$ and $\int \mathbf{U}$ with these equivalences.
\end{proof}
We conclude this chapter with a property of the above adjunction which will be useful in the next chapter.
\begin{prop}
\[\begin{tikzcd}\RMat \ar[r,bend left,"F"] \ar[r,"\bot",phantom] & \RCat \ar[l,bend left,"U"] .\end{tikzcd} \]
is an idempotent adjunction.
\end{prop}

\begin{proof}
Every adjunction between posets is idempotent. Therefore the smaller adjunctions $F_{X} \dashv U_{X}$ are idempotent. Because $F$ and $U$ are stitched together using these adjunctions, it is idempotent as well.
\end{proof}

\chapter{Compositionality of the Algebraic Path Problem}\label{compalgpathchap}In this chapter we show how the algebraic path problem extends to the syntax of open $R$-matrices, i.e., $R$-matrices equipped with input and output nodes. In Section \ref{openmat}, we define open $R$-matrices and construct a symmetric monoidal double category $\Open(\RMat)$ whose horizontal composition is gluing of open $R$-matrices. In Section \ref{openalg}, we show that the algebraic path problem functor
\[F \maps \RMat \to \RCat \]
of Theorem \ref{freecatinternal} lifts to a symmetric monoidal double functor
\[\Open(F) \maps \Open(\RMat) \to \Open(\RCat).\] This double functor describes the compositionality of solutions to the algebraic path problem with respect to gluing of open $R$-matrices. In Section \ref{functional}, we define the black-boxing $\blacksquare(C)$ of an open $R$-category $C \maps X \to Y$. $\blacksquare(C)$ is a $R$-matrix whose values contain only the entries of $C$ which go from input nodes to output nodes. In Theorem \ref{blackbox} we show that the black-boxing functor is in general laxly functorial. However, in Theorem \ref{functional} we show that black-boxing is strictly functorial on functional open $R$-matrices, i.e.\ $R$-matrices for which every input is a source and every output is a sink. This gives a useful expression for solving the algebraic path problem compositionally.

\section{Open $R$-Matrices}\label{openmat}
$R$-matrices are made open by designating some of their vertices to be either inputs or outputs. In this section we show how these open $R$-matrices are composed by gluing the output vertices of one to the input vertices of another and adding the $R$-matrices on the overlap. To define open $R$-matrices, we need a notion of a discrete weighted matrix on a set. The map sending a set to its discrete $R$-matrix is a functor and a left adjoint.
\begin{prop}\label{Lenriched}
Let $R \maps \RMat \to \Set$ be the functor which sends a weighted graph to its underlying set of vertices and sends a morphism to its underlying function. Then $R$ has a left adjoint
\[ 0 \maps \Set \to \RMat\]
which sends a set $X$ to the $R$-weighted graph
\[0_X \maps  X \times X \to Y \]
defined by $0_X (i,j) = 0$ for all $i$ and $j$ in $X$. $F$ sends a function $f \maps X \to Y$ to the morphism of $R$-matrices which has $f$ as its underlying function between vertices. 
\end{prop}
\begin{proof}
The natural isomorphism 
\[\phi\maps \RMat (0_X, G) \cong \Set (X, R(G))\]
is formed by noting that a morphism $ 0_X \to R(G)$ is uniquely determined by its underlying function on vertices and every such function obeys the inequality in Definition \ref{matr}.
\end{proof}
A weighted graph can be opened up to its environment by equipping it with inputs and outputs.
\begin{defn}Let $M: A \times A \to R$ be an $R$-matrix.
An \define{open $R$-matrix} $M \maps X \to Y$ is a cospan in $\RMat$ of the form
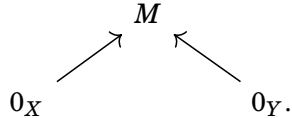
\[ \begin{tikzcd} & M & \\
 0_X \ar[ur] & & 0_Y. \ar[ul]\end{tikzcd}\]
\end{defn}
\noindent The idea is that the maps of this cospan point to input and output nodes of the matrix $M$. 

We can compose open $R$-matrices using the theory of structured cospans; concretely the composition has the following more elementary description. Let $M \maps X \to Y$ and $N \maps Y \to Z$
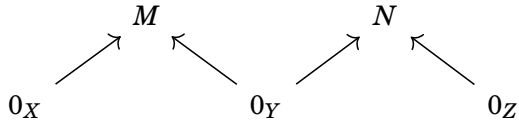
\[
\begin{tikzcd}
& M & & N & \\
0_X \ar[ur] & & 0_Y \ar[ul] \ar[ur] & & \ar[ul] 0_Z 
\end{tikzcd}
\]
be open $R$-matrices. The underlying sets of $M$ and $N$ form a diagram 
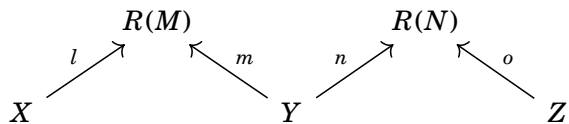
\[ \begin{tikzcd}
& R(M) & & R(N) & \\
X \ar[ur,"l"] & & Y \ar[ul,"m",swap] \ar[ur,"n"] & & \ar[ul,"o",swap] Z 
\end{tikzcd}\]
which generate a pushout
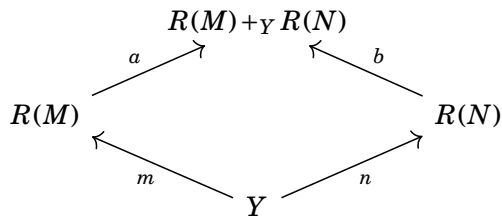
\[ 
\begin{tikzcd}
& R(M) +_Y R(N) & \\
R(M) \ar[ur,"a"] & & R(N) \ar[ul,"b",swap] \\
& Y \ar[ul,"m"] \ar[ur,"n",swap] & 
\end{tikzcd}
\]
The functions $a$ and $b$ of this pushout allow the matrices $M$ and $N$ to be compared on equal footing: the pushforwards $a_*(M)$ and $b_*(N)$ both have $R(M) +_Y R(N)$ as their underlying set. The matrices $a_*(M)$ and $b_*(N)$ are combined using pointwise sum.
\begin{defn}
For open $R$-matrices $M: X \to Y$ and $N \maps Y \to Z$ as defined above, their \define{composite} is defined by
\[N \circ M \maps X \to Z = \begin{tikzcd}&a_*(M) + b_*(N)& \\
LX \ar[ur,"\phi^{-1}(a \circ l)" ] & & LZ \ar[ul,"\phi^{-1}(b \circ r)",swap] \end{tikzcd} \]
where $\phi^{-1}$ gives the unique morphism out of a discrete $R$-matrix defined by a function on its underlying set.
\end{defn}

An $R$-matrix $M \maps X \times X \to R$ can represent a graph with vertex set $X$ weighted in $R$. Similarly, an open $R$-matrix, represents an $R$-weighted graph equipped with inputs and outputs. For example, let $R$ be the quantale $[0, \infty]$ where the addition is infimum and the multiplication is ordinary addition. Then the $[0,\infty]$-matrix
\[
\begin{bmatrix}
1 &2 & .1 \\
3 & 0 & .2 \\
\infty & 1 & .2 
\end{bmatrix}
\]
on the set $\{a,b,c\}$ can be regarded as on open $[0,\infty]$-matrix with left input set $\{1,2\}$ and right input set $\{3\}$. The mappings of the cospan are given by $1 \mapsto a, 2 \mapsto b$ and $3 \mapsto c$. This can be drawn as an open weighted graph
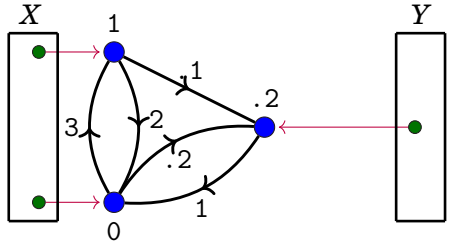
\[
\begin{tikzpicture}
	\begin{pgfonlayer}{nodelayer}
	    
		\node [style=empty] (X) at (-5.1, 2) {$X$};
		\node [style=none] (Xtr) at (-4.75, 1.75) {};
		\node [style=none] (Xbr) at (-4.75, -0.75) {};
		\node [style=none] (Xtl) at (-5.4, 1.75) {};
             \node [style=none] (Xbl) at (-5.4, -0.75) {};
	
		\node [style=inputdot] (1) at (-5, 1.5) {};
		\node [style=inputdot] (5) at (-5,-.5) {};

		\node [style=empty] (Y) at (0.1, 2) {$Y$};
		\node [style=none] (Ytr) at (.4, 1.75) {};
		\node [style=none] (Ytl) at (-.25, 1.75) {};
		\node [style=none] (Ybr) at (.4, -0.75) {};
		\node [style=none] (Ybl) at (-.25, -0.75) {};

		\node [style=inputdot] (3) at (0, 0.5) {};

		\node [style=node,label={\tt 1}] (A) at (-4, 1.5) {};
		\node [style=node,label=below:{\tt 0}] (A2) at (-4,-.5) {};
		\node [style=node,label={\tt .2}] (A3) at  (-2,.5) {};
		
	\end{pgfonlayer}
	\begin{pgfonlayer}{edgelayer}
		\draw [style=inputarrow] (1) to (A);

		\draw [style=inputarrow] (3) to (A3);
		
		\draw [style=inputarrow] (5) to (A2);
		\draw [style=arrow] (A) to [bend left] node [right,yshift=.1cm] {\tt 2} (A2);
		\draw [style=arrow] (A2) to [bend left] node [left] {\tt 3} (A); 
		\draw [style=arrow] (A) to node [above] {\tt .1} (A3);
		\draw [style=arrow] (A3) to [bend left] node [below] {\tt 1} (A2);
		\draw [style=arrow] (A2) to [bend left]  node [below] {\tt .2} (A3);
		\draw [style=simple] (Xtl.center) to (Xtr.center);
		\draw [style=simple] (Xtr.center) to (Xbr.center);
		\draw [style=simple] (Xbr.center) to (Xbl.center);
		\draw [style=simple] (Xbl.center) to (Xtl.center);
		\draw [style=simple] (Ytl.center) to (Ytr.center);
		\draw [style=simple] (Ytr.center) to (Ybr.center);
		\draw [style=simple] (Ybr.center) to (Ybl.center);
		\draw [style=simple] (Ybl.center) to (Ytl.center);
	\end{pgfonlayer}
\end{tikzpicture}
\]
where unlabeled edges are assumed to have a value of $\infty$. Similarly, we define an open $[0,\infty]$-matrix
on $\{d,e\}$ 
\[
\begin{bmatrix}
6 & \infty\\
0& 9
\end{bmatrix}
\]
with left input set given by $\{3\}$ and right input set given by $\{4\}$. The mappings in the cospan for this open $[0,\infty]$-matrix are given by the assignments $3 \mapsto d$ and $4 \mapsto e$. This open $[0,\infty]$-matrix is drawn as
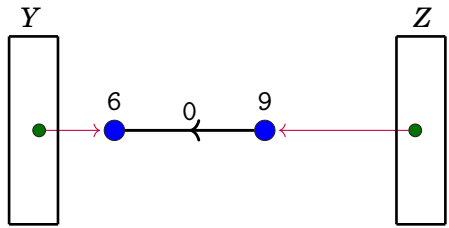
\[
\begin{tikzpicture}
	\begin{pgfonlayer}{nodelayer}
	    
		\node [style=empty] (X) at (-5.1, 2) {$Y$};
		\node [style=none] (Xtr) at (-4.75, 1.75) {};
		\node [style=none] (Xbr) at (-4.75, -0.75) {};
		\node [style=none] (Xtl) at (-5.4, 1.75) {};
             \node [style=none] (Xbl) at (-5.4, -0.75) {};
	
		\node [style=inputdot] (1) at (0, .5) {};

		\node [style=empty] (Y) at (0.1, 2) {$Z$};
		\node [style=none] (Ytr) at (.4, 1.75) {};
		\node [style=none] (Ytl) at (-.25, 1.75) {};
		\node [style=none] (Ybr) at (.4, -0.75) {};
		\node [style=none] (Ybl) at (-.25, -0.75) {};
		\node [style=inputdot] (3) at (-5, 0.5) {};
		\node [style=node,label={\tt 9}] (A) at (-2, .5) {};
		\node [style=node,label={\tt 6}] (A3) at  (-4,.5) {};

	\end{pgfonlayer}
	\begin{pgfonlayer}{edgelayer}
		\draw [style=inputarrow] (1) to (A);

		\draw [style=inputarrow] (3) to (A3);

		\draw[style=arrow] (A) to node [above] {\tt 0} (A3);

		\draw [style=simple] (Xtl.center) to (Xtr.center);
		\draw [style=simple] (Xtr.center) to (Xbr.center);
		\draw [style=simple] (Xbr.center) to (Xbl.center);
		\draw [style=simple] (Xbl.center) to (Xtl.center);
		\draw [style=simple] (Ytl.center) to (Ytr.center);
		\draw [style=simple] (Ytr.center) to (Ybr.center);
		\draw [style=simple] (Ybr.center) to (Ybl.center);
		\draw [style=simple] (Ybl.center) to (Ytl.center);
	\end{pgfonlayer}
\end{tikzpicture}
 \]
The composite of these two $[0,\infty]$-matrices is represented by
 
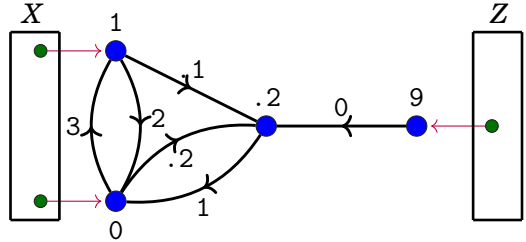
\[
\begin{tikzpicture}
	\begin{pgfonlayer}{nodelayer}
	    
		\node [style=empty] (X) at (-5.1, 2) {$X$};
		\node [style=none] (Xtr) at (-4.75, 1.75) {};
		\node [style=none] (Xbr) at (-4.75, -0.75) {};
		\node [style=none] (Xtl) at (-5.4, 1.75) {};
             \node [style=none] (Xbl) at (-5.4, -0.75) {};
	
		\node [style=inputdot] (1) at (-5, 1.5) {};
		\node [style=inputdot] (5) at (-5,-.5) {};
		\node [style=empty] (Y) at (1.1, 2) {$Z$};
		\node [style=none] (Ytr) at (1.4, 1.75) {};
		\node [style=none] (Ytl) at (.75, 1.75) {};
		\node [style=none] (Ybr) at (1.4, -0.75) {};
		\node [style=none] (Ybl) at (.75, -0.75) {};

		\node [style=inputdot] (3) at (1, 0.5) {};

		\node [style=node,label={\tt 1}] (A) at (-4, 1.5) {};
		\node [style=node,label=below:{\tt 0}] (A2) at (-4,-.5) {};
		\node [style=node,label={\tt .2}] (A3) at  (-2,.5) {};
		\node [style=node,label={\tt 9}] (A4) at  (0,.5) {};
	\end{pgfonlayer}
	\begin{pgfonlayer}{edgelayer}
		\draw [style=inputarrow] (1) to (A);

		\draw [style=inputarrow] (3) to (A4);
		
		\draw [style=inputarrow] (5) to (A2);
			\draw [style=arrow] (A) to [bend left] node [right,yshift=.1cm] {\tt 2} (A2);
		\draw [style=arrow] (A2) to [bend left] node [left] {\tt 3} (A); 
		\draw [style=arrow] (A) to node [above] {\tt .1} (A3);
		\draw [style=arrow] (A3) to [bend left] node [below] {\tt 1} (A2);
		\draw [style=arrow] (A2) to [bend left]  node [below] {\tt .2} (A3);
		\draw [style=arrow] (A4) to node [above] {\tt 0} (A3);
		\draw [style=simple] (Xtl.center) to (Xtr.center);
		\draw [style=simple] (Xtr.center) to (Xbr.center);
		\draw [style=simple] (Xbr.center) to (Xbl.center);
		\draw [style=simple] (Xbl.center) to (Xtl.center);
		\draw [style=simple] (Ytl.center) to (Ytr.center);
		\draw [style=simple] (Ytr.center) to (Ybr.center);
		\draw [style=simple] (Ybr.center) to (Ybl.center);
		\draw [style=simple] (Ybl.center) to (Ytl.center);
	\end{pgfonlayer}
\end{tikzpicture}
\]
where edges are omitted if their weight is infinite. The matrix on the apex of this composite is computed by pushing each component matrix forward to the pushout of their underlying sets and adding them together i.e.\ 
\[
\begin{bmatrix}
1 &2 & .1 & \infty \\
3 & 0 & .2 & \infty  \\
\infty & 1 & .2 & \infty \\
\infty & \infty & \infty & \infty 
\end{bmatrix}
+
\begin{bmatrix}
\infty & \infty & \infty & \infty \\
\infty & \infty & \infty & \infty \\
\infty & \infty & 6 & \infty\\
\infty & \infty & 0 & 9
\end{bmatrix}
=
\begin{bmatrix}
1 & 2 & .1 & \infty\\
3 & 0 & .2 & \infty \\
\infty& 1 & .2 & \infty \\
\infty &\infty & 0 & 9
\end{bmatrix}
\]
The entries of this matrix represent the shortest distance between pairs of vertices. Next we construct a double category where horizontal morphisms are open $R$-matrices.
\begin{thm}\label{openenriched}
For a quantale $R$, there is a symmetric monoidal double category $\Open(\RMat)$ where 
\begin{itemize}
    \item objects are sets $X$,$Y$,$Z \ldots$
    \item vertical morphisms are functions $f: X \to Y$, 
    \item a horizontal morphism $M \maps X \to Y$ is an open $R$-matrix
    
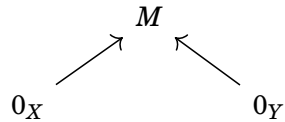
\[
    \begin{tikzcd}
     & M & \\\
    0_X \ar[ur] & & 0_Y \ar[ul]
    \end{tikzcd}
    \]
    \item vertical 2-morphisms are commutative rectangles
    
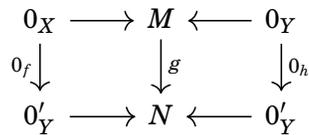
\[
    \begin{tikzcd}
    0_X\ar[d,"0_f",swap] \ar[r] & M \ar[d,"g"] & \ar[l] \ar[d,"0_h"] 0_Y \\
    0_Y' \ar[r] & N & \ar[l] 0_Y'
    \end{tikzcd}
    \]
    \item vertical composition is ordinary composition of functions,
    \item and horizontal composition is given by the composite operation defined above.
    
    
\end{itemize}
The symmetric monoidal structure is given by 
\begin{itemize}
    \item coproducts in $\Set$ on objects and vertical morphisms,
    \item pointwise coproducts on horizontal morphisms i.e.\ for open $R$-matrices,
    
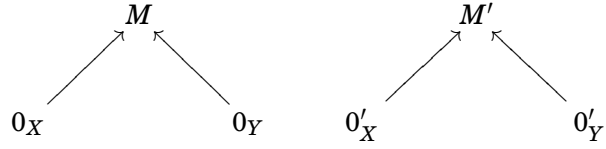
\[
    \xymatrix{ & M & & & M' & \\
	0_X \ar[ur] & & 0_Y \ar[ul] & 0_X' \ar[ur] & & 0_Y' \ar[ul] }
	\]
	their coproduct is
	
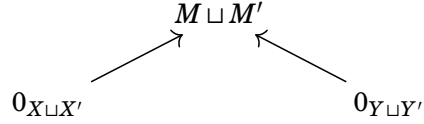
\[\begin{tikzcd}
	& M\sqcup M' & \\
	0_{X\sqcup X'} \ar[ur] & & \ar[ul] 0_{Y \sqcup Y'}
	\end{tikzcd} \]
	where $\sqcup$ indicates the coproduct in $\RMat$. For vertical 2-morphisms,
	\[
	\begin{tikzcd}
    0_X\ar[d,"0_f",swap] \ar[r] & M \ar[d,"g"] & \ar[l] \ar[d,"0_h"] 0_Y \\
    0_Z \ar[r] & N & \ar[l] 0_Q
    \end{tikzcd}
    \begin{tikzcd}
    0_X'\ar[d,"0_f'",swap] \ar[r] & M' \ar[d,"g'"] & \ar[l] \ar[d,"0_h'"] 0_Y' \\
    0_Z' \ar[r] & N' & \ar[l] 0_Q'
    \end{tikzcd}
    \]
    their coproduct is 
    \[\begin{tikzcd}
    0_{X\sqcup X'}\ar[d,"0_{f\sqcup f'}",swap] \ar[r] & M \sqcup M' \ar[d,"g\sqcup g'"] & \ar[l] \ar[d,"0_{h\sqcup h'}"] 0_{Y\sqcup Y'} \\
    0_{Z\sqcup Z'} \ar[r] & N\sqcup N' & \ar[l] 0_{Q\sqcup Q'}.
    \end{tikzcd} \]
\end{itemize}
\end{thm}

\begin{proof}

Lemma \ref{Courser} constructs this symmetric monoidal double category as long as
\begin{itemize}
    \item $\RMat$ has finite coproducts and pushouts,
    \item and $0 \maps \Set \to \RMat$ preserves pushouts and coproducts.
\end{itemize}Because $0$ is a left adjoint (Proposition \ref{Lenriched}) it preserves pushouts and coproducts when they exist so it suffices to prove the following lemma.
\end{proof}

\begin{lem}\label{matpush}
$\RMat$ has coproducts and pushouts.
\end{lem}
\begin{proof}
This is a consequence of Proposition 2.4 of \cite{wolff1974v} after noting that $\RMat$ is the category of $R$-graphs, the generating data for $R$-enriched categories. For concreteness and practicality, we offer an explicit construction of pushouts and coproducts here. 
 Let
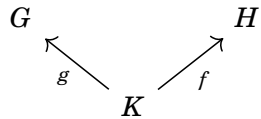
\[ \begin{tikzcd}G & & H \\
& \ar[ur,"f",swap] K \ar[ul,"g"] &\end{tikzcd}\]
be a diagram in $\RMat$ with
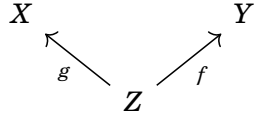
\[ \begin{tikzcd}X & & Y \\
& \ar[ur,"f",swap] Z \ar[ul,"g"] &\end{tikzcd}\]
as the underlying diagram of sets. To compute the pushout $G+_K H$ first we take the pushout of sets
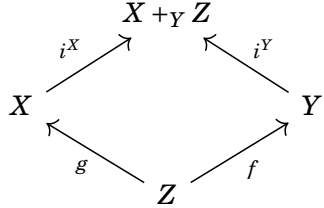
\[
\begin{tikzcd}
 & X+_Y Z & \\
X \ar[ur,"i^X"] & & Y \ar[ul,"i^Y",swap]\\
& \ar[ur,"f",swap] Z \ar[ul,"g"] &\end{tikzcd}
\]
push them forward to get matrices $i^X_*(G)$ and $i^Y_*(H)$ and add them together to get
\[G+_Y H \maps (X+_Y Z) \times (X +_Y Z) \to R = i^X_*(G) + i^Y_*(H). \]

This does indeed define a pushout in $\RMat$. Suppose we have a commutative diagram of $R$-matrices as follows:
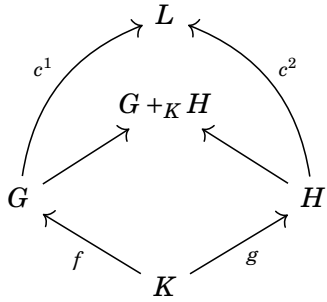
\[
\begin{tikzcd}
& L   &\\
& G+_K H   & \\
G \ar[uur,bend left,"c^1"]\ar[ur]  & & \ar[ul]\ar[uul,bend right,"c^2",swap] H \\
& K \ar[ul,"f"] \ar[ur,"g",swap] &
\end{tikzcd}
\]
then the underlying diagram of sets induces a unique function $u$
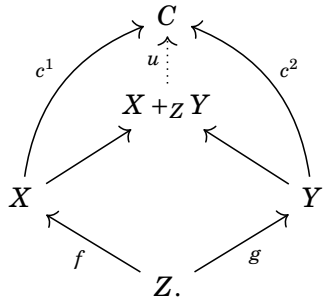
\[
\begin{tikzcd}
& C   &\\
& X+_Z Y  \ar[u,"u",dotted]  & \\
X \ar[ur]\ar[uur,bend left,"c^1"]  & & \ar[ul]\ar[uul,bend right,"c^2",swap] Y \\
& Z \ar[ul,"f"] \ar[ur,"g",swap]. &
\end{tikzcd}
\]
\noindent commuting suitably with $c^1$ and $c^2$. The map $u$ is certainly unique, it remains to show that it is well-defined i.e.\ it satisfies the inequality
\[u_*(G+_K H) \leq L. \]
Indeed, for $(x,y) \in C \times C$,
\begin{align*}
    u_*(G+_K H)(x,y) & = \sum_{(a,b) \in (u \times u)^{-1}(x,y)} (G+_K H)(a,b) \\
    & = \sum_{(a,b) \in (u \times u)^{-1}(x,y)} i^X_*(G)(a,b) + i^Y_*(H) (a,b) \\
    & = \sum_{(a,b) \in (u \times u)^{-1}(x,y)} i^X_*(G)(a,b) + \sum_{(a,b) \in (u \times u)^{-1}(x,y)} i^Y_*(H) (a,b) \\
    & = u_*(i^X_*(G))(x,y) + u_*(i^Y_*(H))(x,y)
\end{align*}
However, because 
\[u_*(i^X(G)) = c^1_*(G) \text{ and } u_*(i^Y(H)) = c^2_*(G) \]
the above expression is equal to
\[c^1_*(G)(x,y) + c^2_*(H) (x,y) \]
which is less than or equal to $L(x,y)$ because each term is and $+$ is the least upper bound.

For $R$-matrices $G \maps X \times X \to R$ and $H \maps Y \times Y \to R$, their coproduct is given by the pushout

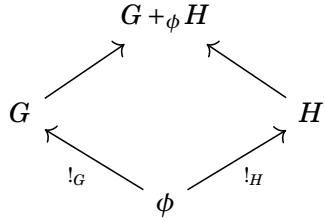
\[
\begin{tikzcd}
& G+_{\phi} H & \\
G \ar[ur]& & H\ar[ul] \\
 & \phi \ar[ul,"!_G"] \ar[ur,"!_H",swap]& 
\end{tikzcd}
\]
where $\phi$ is the initial object of $\RMat$, i.e., the unique $R$-matrix on the empty set, and $!_G$ and $!_H$ are the unique morphisms into $G$ and $H$ respectively.\end{proof}
\section{Compositional Semantics of the Algebraic Path Problem}\label{openalg} Computing the solution to the algebraic path problem on an $R$-matrix $G$ suffers from combinatorial explosion when the size of $G$ grows very large. Therefore, efficient strategies to compute the algebraic path problem must break down large matrices into small pieces, compute the algebraic path problem on each piece, and then combine those solutions together. This strategy can be understood using structured cospans. Suppose that an $R$-matrix $G$ is divided into open $R$-matrices
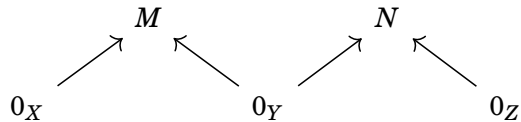
\[ \begin{tikzcd}& M & & N &\\
0_X \ar[ur] & & 0_Y \ar[ul] \ar[ur] & & 0_Z \ar[ul]\end{tikzcd} \]
sharing a common boundary $Y$. We may then apply the algebraic path problem functor $F$ to get two composable cospans of $R$-categories
\[ \begin{tikzcd}& F(M) & & F(N) &\\
1_X \ar[ur] & & 1_Y \ar[ul] \ar[ur] & & 1_Z. \ar[ul]\end{tikzcd} \] The pushout in $\RMat$, $UF(M) +_{1_Y} UF(N)$, is not equal to the solution $F(M+_{0_Y} N)$. The former optimizes over only paths that are the composite of a path in $M$ and a path in $N$. On the other hand, $F(M+_{0_Y} N)$ optimizes over paths that may zig-zag back and forth between $M$ and $N$ as many times as they like before arriving at their destination. Therefore, to construct $F(M+_{0_Y} N)$ from its components we turn to the pushout in $\RCat$.

\begin{prop}
$\RCat$ has pushouts and coproducts.
\end{prop}

\begin{proof}
More generally, $\RCat$ has all colimits by Corollary 2.14 of \cite{wolff1974v}. These colimits are constructed via the transfinite construction of free algebras \cite{kelly1980unified}. The idea behind the transfinite construction is that colimits in a category of monoids can be constructed by first taking the colimit of their underlying objects, taking the free monoid on that colimit, and then quotienting out by the equations in your original monoids. \end{proof}
Next we provide an explicit description of colimits in $\RCat$. 
\begin{prop}\label{colim}
For a diagram $D \maps C \to \RCat$, its colimit is given by the formula
\[\mathrm{colim}_{c \in C} D(c) \cong F( \mathrm{colim}_{c \in C} U(D(c)) ) \]
\end{prop}
\begin{proof}
It suffices to show that $F(\colim_{c\in C} U (D(c)))$ satisfies the universal property of $\colim_{c\in C} D(c)$. Let $\alpha \maps \Delta_d \Rightarrow D$ be a cocone from an object $d \in \RCat$ to our diagram $D$. Because $\alpha$ can be regarded as a cocone in $\RMat$, the universal property of colimits induces a unique map
\[ \colim_{c \in C} U (D(c)) \to U(d)\]
of $R$-matrices. Applying $F$ to this morphism gives a map
\[ F(\colim_{c \in C} U (D(c))) \to FU(d) = d\]
where the last equality follows either from elementary considerations or from the adjunction $F \dashv U$ being idempotent. The above map is a unique morphism satsifying the universal property for $\colim_{c\in C} D(c)$.
\end{proof}

\begin{cor}\label{comp}
For a diagram
\[
\begin{tikzcd}
M & & N \\
    & K \ar[ul] \ar[ur]& 
    \end{tikzcd}\]
    in $\RCat$, the pushout is given by
    \[M+_K N \cong F(U(M) +_{U(K)} U(N)) \]
    and the coproduct of $R$-categories is given by their coproduct in $\RMat$ i.e.
    \[M\sqcup N \cong U(M) \sqcup U(N).\]
\end{cor}
This pushout forms the horizontal composition of a double category of open $R$-categories. 
\begin{thm}\label{openRmatdouble}There is a symmetric monoidal double category $\Open(\RCat)$ where
\begin{itemize}
\item objects are sets,
\item vertical morphisms are functions,
\item horizontal morphisms are cospans 

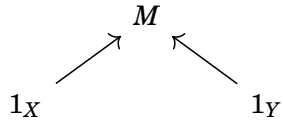
\[ \begin{tikzcd}
 & M & \\
   1_X \ar[ur] & & 1_Y \ar[ul]
\end{tikzcd}\]
where the apex $M$ satisfies the axioms of an $R$-category,
\item and vertical 2-morphisms are commuting rectangles
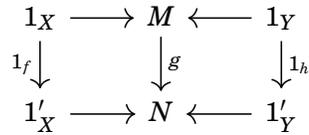
\[
\begin{tikzcd}
    1_X\ar[d,"1_f",swap] \ar[r] &  M \ar[d,"g"] & \ar[l] \ar[d,"1_h"] 1_Y \\
    1_X' \ar[r] & N & \ar[l] 1_Y'
    \end{tikzcd}
\]
\item The horizontal composition is given by pushout of open $R$-categories i.e.\ for open $R$-categories
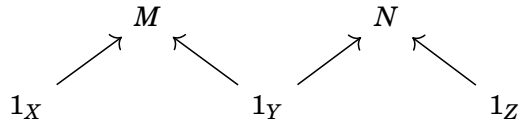
\[
\begin{tikzcd}
& M & & N & \\
1_X \ar[ur] & & 1_Y \ar[ul] \ar[ur] & &\ar[ul] 1_Z
\end{tikzcd}\]
their pushout is the cospan
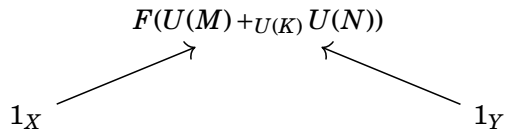
\[
\begin{tikzcd}
& F(U(M)+_{U(K)} U(N)) & \\
1_X \ar[ur] & & \ar[ul] 1_Y
\end{tikzcd}
\]
\end{itemize}
The symmetric monoidal structure of $\Open(\RCat)$ is given by
\begin{itemize}
    \item coproduct of sets and functions,
    \item pointwise coproduct on horizontal morphisms,
    \item and pointiwise coproduct on vertical 2-morphisms.
\end{itemize}

\end{thm}

\begin{proof}
To construct the desired symmetric monoidal double category, we apply Lemma \ref{Courser} to the composite left adjoint \[\begin{tikzcd}\Set \ar[r,"0"] & \RMat \ar[r,"F"] & \RCat. \end{tikzcd}\]
We write the above composite as \[1 \maps \Set \to \RCat \]
as it sends a set $X$ to the identity matrix $1_X$ on $X$.\end{proof}
\noindent So far we have the commutative diagram of functors
\[
\begin{tikzcd}
\RMat  \ar[rr,"F"] &  &\RCat \\
&\Set. \ar[ur,"1",swap] \ar[ul,"0"]&
\end{tikzcd}
\] The definition of $\Open$ is functorial with respect to this sort of diagram, i.e.\ it induces a symmetric monoidal double functor between the relevant double categories.
\begin{thm}\label{functor}
There is a symmetric monoidal double functor
\[\bigstar \maps \Open(\RMat) \to \Open(\RCat) \]
which is 
\begin{itemize}
\item the identity on objects and vertical morphisms,

\item an open $R$-matrix
\[
M \maps X \to Y = 
\begin{tikzcd}
& M & \\
0_X \ar[ur]& & \ar[ul] 0_Y
\end{tikzcd}
\]
is sent to the solution of its algebraic path problem
\[
\bigstar(M) \maps X \to Y =
\begin{tikzcd}
& FM & \\
1_X \ar[ur]& & \ar[ul] 1_Y,
\end{tikzcd}
\]
and
\item a vertical 2-morphism of open $R$-matrices 
\[
 \alpha \maps M \Rightarrow N =	\begin{tikzcd}
    0_X\ar[d,"0_f",swap] \ar[r] & M \ar[d,"g"] & \ar[l] \ar[d,"0_h"] 0_Y \\
    0_Z \ar[r] & N & \ar[l] 0_Q
    \end{tikzcd}
\] 
is sent to the 2-morphism given by pointwise application of $F$
\[ \bigstar(\alpha) \maps M \Rightarrow N =
	\begin{tikzcd}
    1_X\ar[d,"1_f",swap] \ar[r] & FM \ar[d,"Fg"] & \ar[l] \ar[d,"1_h"] 1_Y \\
    1_X' \ar[r] & FN & \ar[l] 1_Y'.
    \end{tikzcd}
\] 
\end{itemize}
\end{thm}

\begin{proof}
We apply Lemma \ref{openfunctoriality} to the square
\[\begin{tikzcd}\RMat \ar[r,"F"] & \RCat\\
\Set \ar[u,"0"] \ar[r,equals] & \Set \ar[u,"1",swap]
\end{tikzcd}
\]
to obtain the desired double functor. Because $F$ preserves pushouts, this double functor preserves horizontal composition and monoidal product up to isomorphism.
\end{proof}\noindent The definition of symmetric monoidal double functor packages up a lot of information very succinctly. In particular, it contains a coherent comparison isomorphism relating the solution of the algebraic path problem on a composite matrix to the solution on its components. For open $R$-matrices $M: X \to Y$ and $N : Y \to Z$, there is a composition comparison
\begin{equation}\label{funct}\phi_{MN} \maps \bigstar(M)\circ \bigstar(N)  \xrightarrow{\sim} \bigstar(M \circ N)\end{equation}
and monoidal comparison
\begin{equation}\label{functmon}
    \psi_{MM'} \maps \bigstar(M + M') \xrightarrow{\sim} \bigstar(M) + \bigstar(M')
\end{equation}
%
giving recipes to break solutions to the algebraic path problem into their components. In other words, the left-hand side of each comparison is computed to determine the right-hand side

Pouly and Kohlas present a similar relationship in the context of valuation algebras \cite[\S 6.7]{pouly2012generic}. For matrices $M$ and $N$ representing weighted graphs on vertex sets $s$ and $t$ respectively, the solution to the algebraic path problem on the union of their vertex sets is given by
\[F(M) \otimes F(N) = F\left( F(M)^{\uparrow s \cup t} + F(N)^{\uparrow s \cup t} \right)\]
In this formula, $\uparrow s \cup t$ indicates that the matrix is trivially extended to the union of the vertex sets. This formula is less general than comparison \ref{funct}: it corresponds to the special case when the legs of the open $R$-matrices are inclusions.

A typical algorithm for the algebraic path problem has spatial complexity $\Theta(n^3)$ where $n$ is the number of vertices in your weighted graph \cite{hofner2012dijkstra}. The comparisons \ref{funct} and \ref{functmon} suggests a faulty strategy for computing the solution to the algebraic path problem which reduces this complexity. First cut your weighted graph into smaller chunks, compute the solution to the algebraic path problem on those chunks, then combine their solutions using \ref{funct} and \ref{functmon}. Unfortunately, this strategy will in general to take \emph{more} time to compute the solution to the algebraic path problem on a composite because the right hand side of comparison \ref{funct} requires three applications of the functor $F$. However, the situation improves if the open $R$-matrices are functional. 

\section{Functional Open Matrices}\label{functional}
In this section we define functional open $R$-matrices, a class of open $R$-matrices for which the composition comparison
\[\phi_{MN} \maps \bigstar(M) \circ \bigstar(N) \xrightarrow{\sim} \bigstar(M \circ N) \]
can be expressed in terms of matrix multiplication. The one caveat is that this expression requires that the open matrices be restricted to their inputs and outputs as follows.
\begin{defn}
Let $M \maps X \to Y$ be the open $R$-category
\[\begin{tikzcd}
& M & \\
1_X \ar[ur,"i"] & & \ar[ul,"o",swap] 1_Y. 
\end{tikzcd} \] Then the \define{black-boxing} of $M$ is the matrix
\[\blacksquare (M) \maps X \times Y \to R \]
given by 
\[\blacksquare(M) ( x,y) = M(i(x),o(y)). \]
\end{defn}\noindent At first the relationship between this black-boxing and the black-boxing of Theorem \ref{blackboxgraph} and Theorem \ref{blackboxqnet} may be opaque. In these theorems we considered black-boxings as profunctors 
\[P \maps A \times B \to \Set \]
where $A$ and $B$ are discrete categories containing only identity morphisms. When enriching in the quantale $R$, we replace $\Set$ with $R$ to obtain a function 
\[P \maps A \times B \to R. \] In general $R$-enriched profunctors must satisfy axioms expressing compatibility with $R$-category structure of $A$ and $B$. However, because $A$ and $B$ are discrete, these axioms become trivial and the above $R$-enriched profunctor is exactly the same as a matrix valued in $R$. Furthermore, as shown in \cite{matricesareenriched}, this correspondence embeds the category of $R$-matrices into the category $R$-enriched profunctors. In this thesis we require a double category of $R$-matrices to match the double categories defined earlier in this chapter.
\begin{defn}
Let $\Mat_R$ be the double category where
\begin{itemize}
    \item an object is a set $X$,$Y$,$Z$,$\ldots$
    \item a vertical morphism is a function $f\maps X \to Y$,
    \item a horizontal morphism $M\maps X \to Y$ is a matrix $M \maps X \times Y \to R$,
    \item a vertical 2-morphism from $M \maps X \to Y$ to $N \maps X' \to Y'$ is a square
    \[
    \begin{tikzcd}
    X \ar[d,"f",swap] \ar[r,"M"] & Y \ar[d,"g"] \\
    X' \ar[r,"N",swap] & Y'
    \end{tikzcd}
    \]
    such that 
    \[ \sum_{x \in f^{-1}(x'),\, y \in g^{-1}(y')} M(x,y) \leq N(x',y')  \]
    for all $x' \in X'$ and $y' \in Y'$,
    \item vertical composition is function composition, and
    \item horizontal composition is given by matrix multiplication.
\end{itemize}
\end{defn}\noindent In this double category, the composite of matrices $M$ and $N$ is written as the juxtaposition $MN$. Black-boxing is extended to the double category of open $R$-categories.
\begin{thm}\label{blackboxenriched}
There is a lax double functor 
\[\blacksquare \maps \Open(\RCat) \to \Mat_R \]
which  
\begin{itemize}
\item is the identity on objects,
\item sends an open $R$-category $M\maps X \to Y$ to its black-boxing $\blacksquare(M)$, and
\item sends a vertical 2-cell
\[
\begin{tikzcd}
1_X \ar[r] \ar[d,"1_f",swap]  &\ar[d,"g"] M & 1_Y\ar[d,"1_h"] \ar[l] \\
1_{X'} \ar[r] & N &\ar[l] 1_{Y'}
\end{tikzcd}
\]
to the vertical 2-cell
\[
\begin{tikzcd}
X \ar[d,"f",swap]\ar[r,"\blacksquare(M)"] & Y \ar[d,"g"] \\
X' \ar[r,"\blacksquare(N)",swap] & Y'.
\end{tikzcd}
\]
\end{itemize}
\end{thm}

\begin{proof}
First observe that this lax double functor is well-defined on 2-cells. This amounts to showing that the inequality 
\begin{equation}\label{wts} \sum_{x \in f^{-1}(x'),\, y \in h^{-1}(y')} M(i(x),j(y)) \leq N(i'(x'),j'(y'))\end{equation}
holds. Because $g$ is a morphism of $R$-matrices, we have that
\begin{equation}\label{step} \sum_{a \in g^{-1}(i'(x')),\,b\in g^{-1}(j'(y'))} M(a,b) \leq N(i'(x'),j'(y')) \end{equation}Let $M(i(x),j(y))$ be a term on the left hand side of inequality \ref{wts}.
Then by definition, $x'=f(x)$ and $y'=h(y)$ so $ a \in g^{-1}(i'(f(x))$ and $b \in g^{-1}(j'(h(y))$. However, because we started with a 2-cell in $\Open(\RCat)$, $i' \circ f = g \circ i$ and $j' \circ h = g \circ j$ so we can rewrite inequality \ref{step} as 
\[ \sum_{a \in g^{-1}(g \circ i(x)),\, b \in g^{-1} (g \circ j(y))} M(a,b) \leq N(i'(x'), j'(y'))\]
The term $M(i(x),j(y))$ of the left hand side of inequality \ref{wts} is also a term of the left hand side of inequality \ref{step} so we have that 
\[M(i(x),i(y)) \leq \sum_{a \in g^{-1}(g \circ i(x)),\, b \in g^{-1} (g \circ j(y))} M(a,b) \leq N(i'(x'), j'(y'))  \]
Because each term on the left hand side of \ref{wts} is less than the desired quantity, the join of all the terms will be as well. Therefore the lax double functor is well-defined on 2-cells. Note that $\Mat_R$ is locally posetal, i.e.\ for every square 
 \[
    \begin{tikzcd}
    X \ar[d,"f",swap] \ar[r,"M"] & Y \ar[d,"g"] \\
    X' \ar[r,"N",swap] & Y'
    \end{tikzcd}
    \]
    there is at most one 2-cell filling it.
This property makes it so many of the axioms in the definition of lax double functor are satisfied trivially. It suffices to show that the globular composition and identity comparisons exist. The identity morphism in $\Open(\RCat)$ on a set $X$ is the cospan
\[\begin{tikzcd}
& 1_X & \\
1_X \ar[equals,ur] &  & \ar[equals,ul,,swap]1_X\end{tikzcd}  \]
The black-box of this cospan is equal to the identity matrix on $X$, so the identity comparison is the identity.
The composition comparison
\[\blacksquare(M) \blacksquare(N) \leq \blacksquare(M \circ N) \]
follows from the chain of inequalities
\begin{align*}
   \blacksquare(M) \blacksquare(N)& = \sum_{y \in Y} \blacksquare(M)(x,y) \blacksquare(N)(y,z) \\
   &= \sum_{y \in Y} M(i(x),j(y))N(i'(y),j'(z)) \\
   &= (M +_{1_Y} N)^2 \\
   & \leq \sum_{n geq 0} (M +_{1_Y} N)^n (i(x),j'(z)) \\
   & = \blacksquare ( M \circ N) (x,z).
\end{align*}\qedhere \end{proof}
We can compose the black-boxing operation with the algebraic path problem functor to get a lax symmetric monoidal double functor
\[\Open(\RMat) \xrightarrow{\bigstar} \Open(\RCat) \xrightarrow{\blacksquare} \Mat_R \]
This lax symmetric monoidal double functor gives the solution to the algebraic path problem on an open $R$-matrix \emph{restricted to its boundaries}. It is natural to ask when this mapping is strictly functorial, as this yields a very simple compositional formula for the algebraic path problem:
\[ \blacksquare( \bigstar(M \circ N)) = \blacksquare(\bigstar(M)) \blacksquare(\bigstar(N)).\]
The double functor $\blacksquare \circ \bigstar$ is strictly functorial on ``functional" open $R$-matrices. 
\begin{defn}\label{functmat} Let $M \maps A \times A \to R$ be an $R$-matrix. An element $a \in X$ is a \define{source} if for every $b \in X$, $M(b,a)=0$ and a \define{sink} if $M(a,b)=0$. A \define{functional open $R$-matrix} is an open $R$-matrix 
\[
\begin{tikzcd}
& M & \\
0_X \ar[ur,"l"] & & \ar[ul,"r",swap] 0_Y
\end{tikzcd}
\]
such that for every $x \in X$, $l(x)$ is a source and for every $y \in Y$, $r(y)$ is a sink.
\end{defn}

 

\noindent Because the composite of functional open $R$-matrices is also functional, we can form the following sub-double category.
\begin{defn}
Let $\Open(\RMat)_{fxn}$ be the full sub-symmetric monoidal double category generated by the open $R$-matrices which are functional.
\end{defn}

\begin{thm}\label{strict}
The composite $\blacksquare \circ \bigstar$ restricts to a strict double functor
\[\blacksquare \circ \bigstar_{fxn} \maps \Open(\RMat)_{fxn} \to \Mat_R \]
\end{thm}\noindent The proof of this theorem relies on a lemma which resembles the the binomial expansion of $(a+b)^n$ in the case when $ba=0$. If $a$ and $b$ represent black-boxes of functional open matrices, then the identity $ba=0$ indicates that there are no paths which go backwards.
\begin{lem}\label{binomial}
For functional open $R$-matrices $M \maps X \to Y$ and $N \maps Y \to Z$ we have that
\[\blacksquare (M +_{1_Y} N)^n = \sum_{i+j=n} \blacksquare(M^i)\blacksquare(N^j).\]
\end{lem}
\begin{proof}
The entries of the left hand side are expanded as 
\[\blacksquare( (M +_{1_Y} N)^n)(a_0,a_{n}) = \sum_{a_1,a_2,\ldots,a_{n-1}} (M +_{1_Y} N)(a_0,a_1) (M +_{1_Y} N)(a_1,a_2) \ldots (M +_{1_Y} N)(a_{n-1},a_n)  \]
where the $a_i$ are equivalence classes in $RM +_Y RN$. For a particular term of this sum, let $1 \leq k \leq n$ be the first natural number such that $a_k$ contains an element of $RN$. Because $M$ and $N$ are functional, for $k \leq i \leq n$ the equivalence classes $a_i$ must also contain an element of $RN$ if our term is nonzero. Therefore for a fixed $k$ the contribution to the above sum is given by
\[
\sum M(a_0,a_1) \ldots M(a_{k-1},a_k) N(a_k,a_{k+1}) \ldots N(a_{n-1},a_n)
\]
which simplifies to 
\[\blacksquare(M^k)\blacksquare( N^{n-k}) (a_0,a_n). \]
Because $k$ can occur in any entry we have that
\begin{align*}
    \blacksquare((M+_{1_Y} N)^n) & = \sum_{k \leq n} \blacksquare(M^k) \blacksquare(N^{n-k}) \\
    &= \sum_{i+j=n} \blacksquare(M^i) \blacksquare(N^j)
\end{align*}
\end{proof}
\noindent \textbf{Proof of Theorem \ref{strict}:} It suffices to prove that for functional open matrices \[\begin{tikzcd}0_X \ar[r] & M & \ar[l] 0_Y \end{tikzcd}\] and \[\begin{tikzcd} 0_Y \ar[r] & N & \ar[l] 0_Z \end{tikzcd}\] the equation
\[\blacksquare( \bigstar(M \circ N)) = \blacksquare(\bigstar(M)) \blacksquare(\bigstar(N)) \]
holds. Consider the left-hand side:
\begin{align*}
   \blacksquare (\bigstar (M \circ N)) & = \blacksquare \sum_{n \geq 0} (M \circ N)^n \\
   &= \sum_{n \geq 0} \blacksquare (M \circ N)^n \\
   &=\sum_{n \geq 0} \sum_{i+j=n} \blacksquare(M^i)\blacksquare(N^j). \\
\end{align*}
where the third step uses Lemma \ref{binomial}. On the other hand, 
\begin{align*}
    \blacksquare(\bigstar(M))\blacksquare(\bigstar(N)) &= \sum_{i \geq 0} \blacksquare(M^i) \sum_{j \geq 0} \blacksquare(N^j) \\
    &= \sum_{i,j \geq 0} \blacksquare(M^i) \blacksquare(M^j)
\end{align*}
Both sums contain the term $\blacksquare(M^i)\blacksquare(N^j)$ for every value of $i$ and $j$, but the left hand side may contain repeated terms. However, because addition is idempotent, repeated terms don't contribute to the sum and the two sides are the same.\hfill $\square$
\\

The functoriality of Theorem \ref{strict} might not be surprising. It says that if your open matrices are joined together directionally along bottlenecks, then the computation of the algebraic path problem can be reduced to a computation on components. This strategy has already proven sucessful. In \cite{sairam1992divide}, Sairam, Tamassia, and Vitter show how choosing \emph{one way separators} as cuts in a graph, allow for an efficient divide and conquer parallel algorithm for computing shortest paths. In \cite{rathke2014compositional} Rathke, Sobocinksi, and Stephens show how the reachability problem on a 1-safe Petri net can be computed more efficiently by cutting it up into more manageable pieces. Theorem \ref{functor} provides a framework for compositional formulas of this type. In future work we plan on extending the construction of this theorem to many other sorts of discrete event dynamic systems.

Lemma \ref{binomial} also holds independent computational interest. The equation given there gives a novel compositional formula for computing the solution to the algebraic path problem. The author has implemented this formula for the special case of Markov processes \cite{compmarkov}. We hope that this is the start of a more extensive library, made faster and more reliable by the mathematics developed in this chapter.

 \chapter{Conclusion}
There are a few directions of research which would make this thesis more complete:
\begin{itemize}
    \item Enriched graphs are only considered when the enriching category is a quantale. More generally, the theory of this thesis could be developed for graphs enriched in a monoidal closed category $(V, \otimes)$ with all colimits. Quantales are a particularly simple example of these, and we are excited about the possibility of enriching in categories that are not posets. For example, we may enrich in the category $(\Grph, \times)$. In this case a graph enriched in $\Grph$ with vertex set $X$ may be regarded as a function 
    \[M \maps X \times X \to \Grph. \]
    For a pair of vertices $x,y \in X$, $M(x,y)$ is a graph whose vertices represent different ways of turning $x$ into $y$. The edges of $M(x,y)$ may represent higher order relationships between the vertices. An operational semantics and compositional theory may be developed for $\Grph$-enriched graphs which is similar the theory developed in this thesis.
    \item To define $\Q$-nets we used the finitary monad $\Q$ induces on $\Set$. This could be made more general by considering an arbitrary monad on an arbitrary category. For example, let $\mathsf{Meas}$ be the category where objects are sets equipped with a $\sigma$-algebra and morphisms are $\sigma$-algebra preserving functions. Let $G \maps \mathsf{Meas} \to \mathsf{Meas}$ be the Giry monad defined in \cite{Giry}. For a measurable set $(X,\Sigma)$, $G((X,\Sigma))$ is the measurable space of probability measures $\mu \maps \Sigma \to [0,1]$ on $(X, \Sigma)$. $G$ is not a finitary monad and therefore does not come from a Lawvere theory. However, we may still define a $G$-net to be a pair of functions
    \[\begin{tikzcd} T \ar[r,shift left=.5ex,"s"] \ar[r,shift right=.5ex,"t",swap] & G(S) \end{tikzcd} \]
    where $T$ and $S$ are measurable spaces. For a transition $\tau \in T$, the probability distributions $s(\tau)$ and $t(\tau)$ may represent uncertainty about the pre- and post-conditions of the event represented by $\tau$. We may attempt to develop an operational semantics for $G$-nets by turning them into free categories internal to the category of $G$-algebras. However, it remains to be seen whether or not this operational semantics is relevant to their natural interpretation.
    \item In Section \ref{QNet}, the functoriality of the definition of $\Q$-net is explored in detail to understand relationships between different variants of $\Q$-nets. The same could be done for $R$-matrices. A morphism of quantales $f:R \to S$ lifts to a functor between categories $\Mat_R \to \Mat_S$ giving a functorial way to translate between matrices with different weights. There is work to be done to understand how these change of enrichment functors relate their corresponding solutions to the algebraic path problem.
    \item The black-boxing double functors of Theorems \ref{blackbox}, \ref{blackboxqnet}, and \ref{blackboxenriched} may be upgraded to symmetric monoidal double functors. This reflects the fact that black-boxing commutes with placing open networks in parallel.
\end{itemize}
In general, this thesis aims to provide a setting for reasoning about the compositionality of networks, but leaves most of that reasoning to future work. As shown in Example \ref{loop}, an infinitely large operational operational semantics may arise from the composite of very small networks. We believe that the best way to study this sort of emergence is to build up to it slowly. Functional open networks are intended to start at the bottom, i.e.\ they are open networks for which the behavior on a composite may be entirely derived from the behavior on its components as shown in Theorems \ref{functionalgraph}, \ref{blackboxfunctqnet}, and \ref{strict}.
We may attempt to prove similar theorems for networks which do exhibit emergent behavior when composed. This may be easier when we choose an operational semantics which is bounded or restricted in some way. For example, we may consider an operational semantics of networks containing paths with a length at most $n$. It is an open question whether or not this operational semantics gives a double functor whose domain is a structured cospan double category of open networks.
\appendix
\setcounter{chapter}{0}
\chapter{Double Categories}\label{appendixdouble}

What follows is a brief introduction to double categories. A more detailed exposition can be found in the work of Grandis and Par\'e \cite{GP1,GP2}, and for monoidal double categories the work of Shulman \cite{Shulman2}.  We use `double category' to mean what earlier authors called a `pseudo
double category'.

\begin{defn}
\label{defn:double_category}
A \textbf{double category} is a category weakly internal to $\Cat$. More explicitly, a double category $\lD$ consists of:
\begin{itemize}
\item a \define{category of objects} $\lD_0$ and a \define{category of arrows} $\lD_1$,
\item  \define{source} and \define{target} functors
\[  S,T \colon \lD_1 \to \lD_0 ,\]
an \define{identity-assigning} functor
\[  U\colon \lD_0 \to \lD_1 ,\]
and a \define{composition} functor
\[ \circ \colon \lD_1 \times_{\lD_0} \lD_1 \to \lD_1 \]
where the pullback is taken over $\lD_1 \xrightarrow[]{T} \lD_0 \xleftarrow[]{S} \lD_1$,
such that
\[  S(U_{A})=A=T(U_{A}) , \quad
	S(M \circ N)=SN, \quad
   T(M \circ N)=TM, \]
\item natural isomorphisms called the \define{associator}
\[ \alpha_{N,N',N''} \maps (N \circ N') \circ N'' \xrightarrow{\sim} N \circ (N' \circ N'') , \]
the \define{left unitor}
\[		\lambda_N \maps U_{T(N)} \circ N \xrightarrow{\sim} N, \]
and the \define{right unitor}
\[  \rho_N \maps N \circ U_{S(N)} \xrightarrow{\sim} N \]
such that $S(\alpha), S(\lambda), S(\rho), T(\alpha), T(\lambda)$ and $T(\rho)$ are all identities and such that the standard coherence axioms hold: the pentagon identity for the 
associator and the triangle identity for the left and right unitor \cite[Sec.\ VII.1]{maclane}.
\end{itemize}
If $\alpha$, $\lambda$ and $\rho$ are identities, we call $\lD$ a \define{strict} double category.
\end{defn}

Objects of $\lD_0$ are called \define{objects} and morphisms in $\lD_0$ are called \define{vertical 1-morphisms}.  Objects of $\lD_1$ are called \define{horizontal 1-cells} of $\lD$ and morphisms in $\lD_1$ are called \define{2-morphisms}.   A morphism $\alpha \maps M \to N$ in $\lD_1$ can be drawn as a square:
\[
\begin{tikzpicture}[scale=1]
\node (D) at (-4,0.5) {$A$};
\node (E) at (-2,0.5) {$B$};
\node (F) at (-4,-1) {$C$};
\node (A) at (-2,-1) {$D$};
\node (B) at (-3,-0.25) {$\Downarrow \alpha$};
\path[->,font=\scriptsize,>=angle 90]
(D) edge node [above]{$M$}(E)
(E) edge node [right]{$g$}(A)
(D) edge node [left]{$f$}(F)
(F) edge node [above]{$N$} (A);
\end{tikzpicture}
\]
where $f = S\alpha$ and $g = T\alpha$.  If $f$ and $g$ are identities we call $\alpha$ a \textbf{globular 2-morphism}.  These give rise to a bicategory:

\begin{defn}
\label{defn:horizontal}
Let $\lD$ be a double category. Then the $\textbf{horizontal bicategory}$ of $\lD$, denoted $H(\lD)$, is the bicategory consisting of objects, horizontal 1-cells and globular 2-morphisms of $\lD$.
\end{defn}

We have maps between double categories, and also transformations between maps:

\begin{defn}
\label{defn:double_functor}
Let $\lA$ and $\lB$ be double categories. A \textbf{double functor} $F \maps \lA \to \lB$ consists of:
\begin{itemize}
\item functors $F_0 \maps \lA_0 \to \lB_0$ and $F_1 \maps \lA_1 \to \lB_1$ obeying the following
equations: 
\[S \circ F_1 = F_0 \circ S, \qquad T \circ F_1 = F_0 \circ T,\]
\item natural isomorphisms called the \define{composition comparison}: 
\[   \phi(N,N') \maps F_1(N) \circ F_1(N') \stackrel{\sim}{\longrightarrow} F_1(N \circ N') \]
and the \define{identity comparison}:
\[  \phi_{A} \maps U_{F_0 (A)} \stackrel{\sim}{\longrightarrow} F_1(U_A) \]
whose components are globular 2-morphisms, 
\end{itemize}
such that the following diagram commmute:
\begin{itemize} 
\item a diagram expressing compatibility with the associator: 
\[\xymatrix{ 	(F_1(N) \circ F_1(N')) \circ F_1(N'') \ar[d]_{\phi (N,N') \circ 1} \ar[r]^{\alpha} & F_1(N) \circ (F_1(N') \circ F_1(N'')) \ar[d]^{1 \circ \phi(N',N'')} \\
			F_1(N \circ N') \circ F_1(N'') \ar[d]_{\phi(N \circ N', N'')} & F_1(N) \circ F_1(N' \circ N'') \ar[d]^{\phi(N, N'\circ N'')}\\
F_1((N \circ N') \circ N'') \ar[r]^{F_1(\alpha)} & F_1(N \circ (N' \circ N'')) }	\]
\item two diagrams expressing compatibility with the left and right unitors:
	\[
	\begin{tikzpicture}[scale=1.5]
	\node (A) at (1,1) {$F_1(N) \circ U_{F_0(A)}$};
	\node (A') at (1,0) {$F_1(N) \circ F_1(U_{A})$};
	\node (C) at (3.5,1) {$F_1(N)$};
	\node (C') at (3.5,0) {$F_1(N \circ U_A)$};
	\path[->,font=\scriptsize,>=angle 90]
	(A) edge node[left]{$1 \circ \phi_{A}$} (A')
	(C') edge node[right]{$F_1(\rho_N)$} (C)
	(A) edge node[above]{$\rho_{F_1(N)}$} (C)
	(A') edge node[above]{$\phi(N,U_{A})$} (C');
	\end{tikzpicture}
	\]
	\[
	\begin{tikzpicture}[scale=1.5]
	\node (B) at (5.5,1) {$U_{F_0(B)} \circ F_1(N)$};
	\node (B') at (5.5,0) {$F_1(U_{B}) \circ F_1(N)$};
	\node (D) at (8,1) {$F_1(N)$};
	\node (D') at (8,0) {$F_1(U_{B} \circ N).$};
		\path[->,font=\scriptsize,>=angle 90]
		(B) edge node[left]{$\phi_{B} \circ 1$} (B')
	(B') edge node[above]{$\phi(U_{B},N)$} (D')
	(B) edge node[above]{$\lambda_{F_1(N)}$} (D)
	(D') edge node[right]{$F_1(\lambda_{N})$} (D);
	\end{tikzpicture}
	\]
\end{itemize}
If the 2-morphisms $\phi(N,N')$ and $\phi_A$ are identities for all $N,N' \in \lA_1$ and 
$A \in \lA_0$, we say $F \maps \lA \to \lB$ is a \define{strict} double functor.  If on the other hand we drop the requirement that these 2-morphisms be invertible, we call $F$ a \define{lax} double
functor.
\end{defn}
	
\begin{defn}
Let $F \maps \lA \to \lB$ and $G \maps \lA \to \lB$ be lax double functors. A \define{transformation} $\beta \maps F \Rightarrow G$ consists of natural transformations $\beta_0 \maps F_0 \Rightarrow G_0$ and $\beta_1 \maps F_1 \Rightarrow G_1$ (both usually written as $\beta$) such that 
		\begin{itemize}
			\item $S( \beta_M) = \beta_{SM}$ and $T(\beta_M) = \beta_{TM}$ for any object $M \in \mathbb{A}_1$, 
			\item $\beta$ commutes with the composition comparison, and
			\item $\beta$ commutes with the identity comparison.
		\end{itemize}
\end{defn}
	
Shulman defines a 2-category $\mathbf{Dbl}$ of double categories, double functors, and transformations \cite{Shulman2}.  This has finite products.  In any 2-category with finite products we can define a pseudomonoid \cite{DayStreet}, which is a categorification of the concept of monoid.  For example, a pseudomonoid in $\mathsf{Cat}$ is a monoidal category.
	
\begin{defn}
	\label{defn:monoidal_double_category}
A \textbf{monoidal double category} is a pseudomonoid in $\mathbf{Dbl}$. Explicitly, a monoidal double category is a double category equipped with double functors $\otimes \maps \lD \times \lD \to \lD$ and $I \maps * \to \lD$ where $*$ is the terminal double category, along with invertible transformations called the \define{associator}:
\[  A \maps \otimes \, \circ \; (1_{\lD} \times \otimes ) \Rightarrow \otimes \; \circ \; (\otimes \times 1_{\lD}) ,\]
\define{left unitor}:
\[ L\maps \otimes \, \circ \; (1_{\lD} \times I) \Rightarrow 1_{\lD} ,\]
and \define{right unitor}:
\[ R \maps \otimes \,\circ\; (I \times 1_{\lD}) \Rightarrow 1_{\lD} \]
satisfying the pentagon axiom and triangle axioms.
\end{defn}

This definition neatly packages a large quantity of information.   Namely:
\begin{itemize}
\item $\lD_0$ and $\lD_1$ are both monoidal categories.
\item If $I$ is the monoidal unit of $\lD_0$, then $U_I$ is the
monoidal unit of $\lD_1$.
\item The functors $S$ and $T$ are strict monoidal.
\item $\otimes$ is equipped with composition and identity comparisons
\[ \chi \maps (M_1\ten N_1)\circ (M_2\ten N_2) \stackrel{\sim}{\longrightarrow}
(M_1\circ M_2)\ten (N_1\circ N_2)\]
\[ \mu \maps U_{A\ten B} \stackrel{\sim}{\longrightarrow} (U_A \ten U_B)\]
making three diagrams commute as in Def.\ \ref{defn:double_functor}.

\item The associativity isomorphism for $\ten$ is a transformation between double functors.
		\item The unit isomorphisms are transformations
between double functors.
	\end{itemize}

	\begin{defn}
	\label{defn:symmetric_monoidal_double_category}
A \define{braided monoidal double category} is a monoidal double
category equipped with an invertible transformation
\[ \beta \maps \otimes \Rightarrow \otimes \circ \tau \]
called the \define{braiding}, where $\tau \maps \lD \times \lD \to \lD \times \lD$ is the twist double functor sending pairs in the object and arrow categories to the same pairs in the opposite order. The braiding is required to satisfy the usual two hexagon identities \cite[Sec.\ XI.1]{maclane}.  If the braiding is self-inverse we say that $\lD$ is a \define{symmetric monoidal double category}.
	\end{defn}
	
In other words:
\begin{itemize}
		\item $\lD_0$ and $\lD_1$ are braided (resp. symmetric) monoidal categories,
		\item the functors $S$ and $T$ are strict braided monoidal functors, and
		\item the braiding is a transformation between double functors.
\end{itemize}

\begin{defn}
\label{defn:monoidal_double_functor}
A \define{monoidal lax double functor} $F \colon \lC \to \lD$ between monoidal double categories $\lC$ and $\lD$ is a lax double functor $F \maps \lC \to \lD$ such that
	\begin{itemize}
		\item $F_0$ and $F_1$ are monoidal functors,
		\item $SF_1= F_0S$ and $TF_1 = F_0T$ are equations between monoidal functors, and
		\item the composition and unit comparisons $\phi(N_1,N_2) \maps F_1(N_1) \circ F_1(N_2) \to F_1(N_1\circ N_2)$ and $\phi_A \maps U_{F_0 (A)} \to F_1(U_A)$ are monoidal natural transformations.
	\end{itemize}
The monoidal lax double functor is \define{braided} if $F_0$ and $F_1$ are braided monoidal functors and \define{symmetric} if they are symmetric monoidal functors. 
\end{defn}

\chapter{Lawvere Theories}\label{law}
Introduced by Lawvere in his landmark thesis \citep{Lawvere}, Lawvere theories are a general framework for reasoning about algebraic structures \citep{ttt,buckley}.

\begin{defn}
A Lawvere theory $\law{Q}$ is a small category with finite products such that every object is isomorphic to the iterated finite product $x^n = x \times \ldots \times x$ for a \define{generic object} $x$ and natural number $n$. Equivalently, Lawvere theories can be thought of as categories whose objects are given by natural numbers $n \in \N$ and with cartesian product given by $+$. The morphisms in a Lawvere theory are called \define{operations}.
\end{defn}
\noindent The idea is that a Lawvere theory represents the platonic embodiment of an algebraic gadget.

\begin{expl}\label{mon}
A canonical example is the Lawvere theory $\law{MON}$ of monoids. Like all Lawvere theories, the objects of $\law{MON}$ are given by natural numbers. In addition $\law{MON}$ contains the morphisms 
\[ m \maps 2 \to 1 \text{ and } e \maps 0 \to 1\]
For a monoid $M$, this represents the multiplication map
\[ \boldsymbol{\cdot} \maps M \times M \to M \]
and the map
\[ e \maps \{*\} \to M\]
which picks out the identity element of $M$. 
These maps are required to satisfy the associative law
\[
\begin{tikzcd}
3 \ar[r, "\mathrm{id} \times m"] \ar[d, "m \times \mathrm{id}",swap] & 2 \ar[d, "m"]\\
2 \ar[r, "m",swap] & 1
\end{tikzcd}
\]
and the unital laws for monoids.
\[
\begin{tikzcd}
1 \ar[r, "\mathrm{id} \times e"] \ar[dr,"\mathrm{id}",swap]& 2\ar[d, "m"] & 1 \ar[l,"e \times \mathrm{id}",swap] \ar[dl,"\mathrm{id}"] \\
 & 1 &
\end{tikzcd}
\]
$\law{MON}$ also contains all composites, tensor products, and maps necessary to make $n$ into the product $x^n$ induced by the maps $m$ and $e$.
\end{expl}
\noindent Like all good things, Lawvere theories form a category. 
\begin{defn}
 Let $\cat{Law}$ be the category where objects are Lawvere theories and morphisms are product preserving functors.
\end{defn}
Note that because morphisms of Lawvere theories preserve products, they must send the generic object of their source to the generic object of their target. Therefore to specify a morphism of Lawvere theories, it suffices to make an assignment of the morphisms which are not part of the product structure.

Let $\law{Q}$ be a Lawvere theory and $C$ a category with finite products. We can impose the axioms and operations of $\law{Q}$ onto an object in $C$ via a product preserving functor $F \maps \law{Q} \to  C$. The image $F(1)$ of the generating object $1$ gives the underlying object of $F$ and for an operation $o \maps n \to k$ in $\law{Q}$, $F(o) \maps F(x)^n \to F(x)^k$ gives a specific instance of the algebraic operation represented by $o$. There is a natural way to make a category of these functors.
\begin{defn}\label{models}
Let $\law{Q}$ be a Lawvere theory and $C$ a category with finite products. Then there is a category $\Mod(\law{Q},C)$ where 
	    \begin{itemize}
			\item objects are product preserving functors $F\maps \law{Q} \to C$ and,
			\item morphisms are natural transformations between these functors.
		\end{itemize}
When $\Mod(\law{Q})$ is written without the second argument, it is assumed to be $\Set$. We will refer to objects in $\Mod(\law{Q})$ as $\law{Q}$\define{-models} and morphisms in $\Mod(\law{Q})$ as $\law{Q}$\define{-model homomorphisms}. When $C=\Cat$, we will refer to these objects as $\Q$\define{-categories}.
\end{defn}
\noindent When the category of models is $\Set$ then there is a forgetful functor 
\[ R_{\law{Q}} \maps \Mod(\law{Q}) \to \Set \] 
which sends a product preserving functor $F \maps \law{Q} \to \Set$ to  image on the generating object $F(1)$ and a natural transformation to  component on the object $1$.
A classical result says that $R_{\law{Q}}$ \emph{always} has a left adjoint
\[ L_{\law{Q}} \maps \Set \to \Mod(\law{Q})\]
which for a set $X$, $L_{\law{Q}} X$ is referred to as the \define{free model of }$\law{Q}$\define{ on }$X$. In fact, this construction extends to fully faithful functor 
\[ \Law \to \Mnd\]
which sends a Lawvere theory $\law{Q}$ to the monad $R_{\law{Q}} \circ L_{\law{Q}} \maps \Set \to  \Set$ and where $\Mnd$ is the category of monads on $\Set$ \citep{linton}. For a Lawvere theory $\law{Q}$ we will denote the monad it induces via this functor by $M_{\law{Q}} \maps \Set \to \Set$.

For $\law{Q}= \law{MON}$, $\Mod(\law{MON}, \Set)$ is equivalent to the category $\cat{Mon}$ of monoids and monoid homomorphisms. In this case the functor $R_{\law{MON}} \maps \cat{Mon} \to \Set$ turns monoids and monoid homomorphisms into their underlying sets and functions. $R_{\law{MON}}$ has a left adjoint 
\[L_{\law{MON}} \maps \Set \to \Mon\]
which sends a set $X$ to the free monoid $L_{\law{MON}} X$. For a function $f \maps X \to Y$, $L_{\law{MON}} f$ is the unique multiplication preserving extension of $f$ to $L_\law{MON} X$.


\end{document}